\documentclass[11pt]{article}
\usepackage[margin=1in]{geometry}
\usepackage[T1]{fontenc}
\usepackage{mathtools,amsthm}
\usepackage{newtxtext,newtxmath}
\usepackage{bm}
\usepackage{microtype}
\usepackage[numbers,sort&compress]{natbib}
\usepackage{algorithm,algpseudocode}
\usepackage{booktabs,array}
\usepackage{graphicx,placeins,needspace,xurl}
\usepackage{enumitem}
\usepackage{xcolor}
\usepackage{aliascnt}
\usepackage{hyperref}
\hypersetup{colorlinks=true,linkcolor=blue!40!black,citecolor=blue!40!black,urlcolor=blue!40!black,pdftitle={Log-optimal two-sample testing by betting},pdfauthor={Aaditya Ramdas},pdfsubject={Sequential two-sample testing and log-optimal e-processes},pdfkeywords={sequential testing, e-process, Jensen--Shannon divergence, characteristic kernel, universal portfolio}}
\usepackage[nameinlink,noabbrev]{cleveref}
\setlist{topsep=3pt,itemsep=2pt}
\numberwithin{equation}{section}
\newtheorem{theorem}{Theorem}[section]
\newaliascnt{proposition}{theorem}
\newtheorem{proposition}[proposition]{Proposition}
\aliascntresetthe{proposition}
\newaliascnt{lemma}{theorem}
\newtheorem{lemma}[lemma]{Lemma}
\aliascntresetthe{lemma}
\newaliascnt{corollary}{theorem}
\newtheorem{corollary}[corollary]{Corollary}
\aliascntresetthe{corollary}
\theoremstyle{definition}
\newaliascnt{definition}{theorem}
\newtheorem{definition}[definition]{Definition}
\aliascntresetthe{definition}
\newaliascnt{assumption}{theorem}

\aliascntresetthe{assumption}
\newaliascnt{example}{theorem}
\newtheorem{example}[example]{Example}
\aliascntresetthe{example}
\theoremstyle{remark}
\newaliascnt{remark}{theorem}
\newtheorem{remark}[remark]{Remark}
\aliascntresetthe{remark}
\newcommand{\E}{\mathbb E}
\newcommand{\PP}{\mathbb P}
\newcommand{\R}{\mathbb R}
\newcommand{\N}{\mathbb N}
\newcommand{\X}{\mathcal X}
\newcommand{\F}{\mathcal F}
\newcommand{\B}{\mathcal B}
\newcommand{\HH}{\mathcal H}
\newcommand{\one}{\mathbf 1}
\newcommand{\KL}{D_{\mathrm{KL}}}
\newcommand{\JS}{\operatorname{JS}}
\newcommand{\kl}{\operatorname{kl}}
\newcommand{\spplus}{\operatorname{sp}}
\DeclareMathOperator{\logaddexp}{logaddexp}
\newcommand{\clip}{\operatorname{clip}}
\newcommand{\sig}{\sigma}
\newcommand{\as}{\quad\text{a.s.}}
\newcommand{\norm}[1]{\left\lVert#1\right\rVert}
\newcommand{\inner}[2]{\left\langle#1,#2\right\rangle}

\title{\textbf{Log-optimal two-sample testing by betting}}
\author{Aaditya Ramdas\\
\texttt{aramdas@stanford.edu}
\\[3pt]\small Department of Statistics, Stanford University}
\date{November 12, 2026}
\begin{document}
\maketitle
\vspace{-0.5em}
\begin{abstract}
For each standard Borel observation space, we construct a single anytime-valid two-sample procedure that is asymptotically log-optimal for every pair of unknown probability laws and every limiting ratio of diverging sample sizes, including zero and infinity. For independent balanced pairs, the optimal log-wealth rate is $2\JS(P,Q)$, where $\JS$ is the Jensen--Shannon divergence. Unequal sizes yield weighted Jensen--Shannon information and the one-sample limit $\log W(n_1,n_2)/n_2\to\KL(Q\|P)$ when $n_1/n_2\to\infty$.  Refining finite partitions and logarithmic-loss learning with measurable characteristic kernels provide universal constructions. 
Our main contribution is the generality of our construction, and its connections to classical batch efficiency notions. Practical implementations allow continuous classifier fine-tuning, as demonstrated on finite alphabets, Gaussian data, and handwritten digits.
\end{abstract}

\section{Introduction}
Sequential two-sample testing asks whether two data sources have the same law while allowing the analyst to inspect accumulating evidence and stop adaptively. A useful way to construct such evidence is to regard a nonnegative test statistic as the wealth of a bettor against the null. Validity restricts the bettor's conditional expected return under the null; power depends on how quickly wealth grows under an alternative. We characterize the largest achievable first-order logarithmic growth rate and construct procedures that attain it.

For independent streams observed one pair at a time, we establish the benchmark
\begin{equation}\label{eq:intro-rate}
G(P,Q)=2\JS(P,Q),\qquad W_n=\exp\{nG(P,Q)+o(n)\}.
\end{equation}
Here $\JS$ is the Jensen--Shannon divergence, defined later. The normalization is per \emph{pair}; per individual observation it is $\JS(P,Q)$. For exchangeable pairs, which may be dependent, the benchmark is $\JS(P\otimes Q,Q\otimes P)$. 
We give separate constructions in general Borel spaces attaining the independent-pair and exchangeable-pair optima; this generality is our main contribution\footnote{A companion paper addresses the related but different problem of independence testing~\cite{ramdas2026independence}.}.

\paragraph{Sequential evidence and logarithmic growth.}
The connection between information and multiplicative betting goes back to \citet{kelly1956}. Modern e-value and test-supermartingale formulations of optional-stopping-safe inference are surveyed by \citet{ramdas2023savi}; growth-rate-optimal e-values and their information-projection interpretation are developed by \citet{grunwald2024safe}. In the two-stream setting, \citet{turner2024generic} develop blockwise e-variables for specified families represented by densities or mass functions, and establish oracle expected-log optimality under convexity. We give a measure-theoretic formulation and construct one data-driven procedure attaining the oracle rate for every pair of probability laws on a standard Borel space, without a prespecified dominating measure. Refining partitions and kernel logarithmic-loss learning supply the universal attainment guarantee.

\paragraph{Nonparametric and predictive two-sample testing.}
The batch kernel-MMD framework of \citet{gretton2012kernel} compares probability laws through their RKHS mean embeddings. \citet{balsubramani2016} develop sequential nonparametric tests using time-uniform martingale bounds. \citet{shekhar2024betting} construct a general two-sample betting framework based on predictably learned integral-probability-metric witnesses, with a practical kernel-MMD instantiation. We analyze their additive payoff family and its relation to the information bounds for independent and exchangeable pairs.

Predictive reductions are also established. \citet{lheritier2018} turn sequential probabilistic prediction into two-sample tests, and \citet{lheritier2019} develop a computationally efficient pointwise-universal predictor for multivariate continuous features. \citet{podkopaev2023predictive} study sequential predictive two-sample and independence tests. \citet{pandeva2024} develop classifier-based e-value tests for high-dimensional data. None of these past works proved log-optimality, especially in such generality as our current work. We use learned class probabilities to attain the information bound for the full common-distribution null. The exact growth shortfall equals the excess logarithmic prediction loss, which determines the learning objective.

\paragraph{Kernels and portfolios.}
The distinction between characteristic and universal kernels is studied by \citet{sriperumbudur2011}. We use a short signed-measure argument to show that adding a constant feature to a bounded characteristic kernel provides the $L^2$ approximation property needed for our logistic construction. The online optimization component is a standard projected-gradient mechanism; see \citet{orabona2019} for the broader online-learning framework. Our analysis gives an explicit regret constant in the two-sample normalization. \citet{cover1991} introduced universal portfolios. \citet{waudbysmith2026} establish instance-log-optimality and rejection-time results from portfolio regret for a broad family of testing problems, including an analysis of betting-fraction restrictions in two-sample witness procedures. Our scalar-portfolio bound applies pathwise to adaptively learned predictors.

\paragraph{A single procedure for all allocation regimes.}
The finite-sample information bound is a weighted Jensen--Shannon divergence. Each construction gives one e-process, fixed before observing the data, that attains the corresponding first-order rate for every limiting sample-size ratio. The ratio may be any positive constant, zero, or infinity. Under extreme imbalance, the rate per minority observation converges to the appropriate directed Kullback--Leibler divergence. \Cref{cor:all-ratios} states the result for a single rule with no knowledge of the limiting allocation. 

\paragraph{Classical fixed-sample efficiency.}
The same wealth processes provide fixed-endpoint p-values with the largest possible Bahadur slope. For the finite-alphabet code, an explicit logarithmic comparison with the homogeneity likelihood ratio yields optimal Stein and Hoeffding error exponents, an alternative-dependent Chernoff calibration, and their minority-sample limits under arbitrary imbalance. At a fixed level, the miss exponent per observation is $-\log\sum_i p_i^\rho q_i^{1-\rho}$, which differs from the Jensen--Shannon significance exponent. Exact permutation calibration of the code is the probability-ordered Fisher--Freeman--Halton test: it preserves the error-exponent guarantees and recovers the classical local likelihood-ratio power. (The direct e-value threshold pays a logarithmic prediction cost and has vanishing power under root-sample-size local alternatives.) \Cref{sec:batch-efficiency} states the results, and \cref{app:batch-proofs} gives the proofs.

\paragraph{Classical antecedents and the scope of the contribution.}
Optimal finite-alphabet homogeneity tests were developed by \citet{gutman1989} and analyzed by \citet{unnikrishnan2012}; \citet{zhu2021} established optimal fixed-level miss exponents for kernel-MMD tests. Recent work gives further finite-alphabet and second-order analyses \citep{harsha2026,grootveld2026}. Exact conditional tables and their likelihood-ratio limits belong to the classical theory \citep{freeman1951,wilks1938}. Our constructions connect these efficiency properties to a single anytime-valid rule with general-space learning and arbitrary limiting sample allocations. Universal-predictor/null-MLE optimality on finite alphabets is also part of \citet{shekhar2025optimal}. General-space learning supplies Bahadur efficiency; the sharper miss-exponent and local-power results use the finite-alphabet coding structure.

\paragraph{Generality and practical implementation.}
The main contribution is universal attainment over arbitrary laws and allocation regimes, with constructions that use only measurable partitions or kernel evaluations. The finite-sample performance depends on how these representations, complexity budgets, and optimization methods are instantiated. Efficient tests for particular spaces and alternatives may benefit from additional problem-specific algorithm design. \Cref{sec:experiments} reports reproducible comparisons of finite implementations with the kernel-MMD betting test of \citet{shekhar2024betting}. The experiments examine both evidence growth and sample efficiency and include settings in which the earlier test or its full-range betting variant stops sooner. Source classifiers also provide practical tests for images. We allow the encoder and source head to adapt throughout the test, using each pair for training only after its bet has been recorded. \Cref{sec:finetuning,sec:image-results} give the protocol and controlled digit-image experiments; Appendices~\ref{app:adaptive-proofs} and~\ref{app:image-experiments} provide the learning analysis and experimental details.

\section{Model, validity, and the optimal rate}\label{sec:model}
Let $(\X,\B)$ be a measurable space. Under an alternative, pairs are independent and identically distributed with law
\[
A=P\otimes Q,
\qquad \F_n=\sigma((X_t,Y_t):1\le t\le n).
\]
Any algorithmic randomization is independent of the data and included in the filtration before it is used. The information bound holds on a measurable space. Universal attainment is established on standard Borel spaces for the partition method and under the measurable-kernel assumptions of Section~\ref{sec:kernel} for the kernel method.

The basic null is the family of i.i.d. pair laws
\begin{equation}\label{eq:null}
\HH_0=\{R\otimes R:R\text{ a probability law on }(\X,\B)\}.
\end{equation}
The constructed processes satisfy a stronger, conditional version: under the null, given $\F_{t-1}$, the two observations are independent with a common, possibly past-dependent law $R_t$.

\begin{definition}
A test supermartingale is a nonnegative adapted process with $W_0=1$ and $\E_0[W_t\mid\F_{t-1}]\le W_{t-1}$ under every null law in the specified model. An e-process is a nonnegative adapted process satisfying $\E_0 W_\tau\le1$ for every bounded stopping time $\tau$ and every null law. The paired constructions before Section~\ref{sec:imbalance} are test supermartingales. The unequal-sample constructions in Section~\ref{sec:imbalance} satisfy the e-process guarantee under sequential sampling. For the upper bounds, it suffices that $\E_{(R\otimes R)^{\otimes n}}W_n\le1$ at every deterministic $n$.
\end{definition}
For $0<\alpha<1$, reject at
\begin{equation}\label{eq:tau}
\tau_\alpha=\inf\{n\ge1:W_n\ge\alpha^{-1}\}.
\end{equation}
The convention is $\inf\varnothing=\infty$. To verify anytime validity directly, fix $N$ and stop a test supermartingale at $\tau_\alpha\wedge N$. Its expectation is at most one, while its value is at least $\alpha^{-1}$ on $\{\tau_\alpha\le N\}$. Consequently $\PP_0(\tau_\alpha\le N)\le\alpha$, and letting $N\to\infty$ gives $\PP_0(\tau_\alpha<\infty)\le\alpha$.

Write $M=(P+Q)/2$. Relative entropy uses the convention $\KL(P\|R)=\infty$ when $P\not\ll R$, and $0\log0=0$. Define
\begin{equation}\label{eq:G}
G(P,Q)=\KL(P\|M)+\KL(Q\|M)=2\JS(P,Q).
\end{equation}
Since $dP/dM,dQ/dM\le2$, $0\le G(P,Q)\le\log4$, with equality at zero exactly when $P=Q$.

\begin{lemma}[Likelihood comparison]\label{lem:comparison}
Let $A\ll B$ be probability laws, let $L=dA/dB$, and let $Z\ge0$ satisfy $\E_B Z\le1$. Then
\[
\E_A[Z/L]\le1,
\qquad \PP_A(\log Z-\log L\ge u)\le e^{-u}\quad(u>0).
\]
The ratio is needed only $A$-almost surely. If $\log L$ is $A$-integrable, then $\E_A\log Z\le\E_A\log L$, allowing $-\infty$ on the left.
\end{lemma}
\begin{proof}
Integration on $\{L>0\}$ gives $\E_A[Z/L]=\int_{\{L>0\}}Z\,dB\le1$. Markov's inequality proves the probability bound, and Jensen's inequality proves $\E_A\log(Z/L)\le0$. The positive part of $\log(Z/L)$ is integrable because $\log^+v\le v$. The expectation statement therefore also holds when the negative part is infinite.
\end{proof}

\begin{theorem}[Unrestricted optimal log-growth]\label{thm:optimal}
For every process valid at deterministic times under \eqref{eq:null},
\begin{align}
\E_{A^{\otimes n}}\log W_n&\le nG(P,Q),\label{eq:expected-upper}\\
\limsup_{n\to\infty}\frac{\log W_n}{n}&\le G(P,Q)\quad A^{\otimes\infty}\text{-a.s.}\label{eq:as-upper}
\end{align}
Moreover,
\begin{equation}\label{eq:projection}
\inf_R\{\KL(P\|R)+\KL(Q\|R)\}=G(P,Q).
\end{equation}
The oracle product process based on
\begin{equation}\label{eq:oracle}
\eta^*(z)=\frac{dP}{d(P+Q)}(z),\qquad
 e^*(x,y)=4\eta^*(x)\{1-\eta^*(y)\}
\end{equation}
is a test supermartingale under the conditional independent-pair null and satisfies $n^{-1}\log\prod_{t=1}^n e^*(X_t,Y_t)\to G(P,Q)$ almost surely.
\end{theorem}
\begin{proof}
Let $B=M\otimes M$. The likelihood ratio of $A$ to $B$ is $e^*$, and that of $A^{\otimes n}$ to $B^{\otimes n}$ is $L_n=\prod_{t=1}^n e^*(X_t,Y_t)$. Its log increment is integrable. Indeed, for $r=dP/dM\le2$, the positive part of $\log r$ under $P$ is bounded by $\log2$, while
\[
\int_{\{r<1\}}|\log r|\,dP
=\int_{\{r<1\}}r|\log r|\,dM\le1/e;
\]
the same argument applies to $Q$. The strong law gives $n^{-1}\log L_n\to G(P,Q)$.

Apply \cref{lem:comparison} at each $n$ to obtain \eqref{eq:expected-upper} and
\[
\PP_A(\log W_n-\log L_n\ge n\varepsilon)\le e^{-n\varepsilon}.
\]
The sum over $n$ is finite. Borel--Cantelli, followed by a countable sequence $\varepsilon\downarrow0$, proves \eqref{eq:as-upper}.

For any $R$, decomposition of the log density ratios gives, in the extended sense,
\begin{equation}\label{eq:pythagorean}
\KL(P\|R)+\KL(Q\|R)
=G(P,Q)+2\KL(M\|R).
\end{equation}
If $P$ or $Q$ is not absolutely continuous with respect to $R$, both the relevant left side and $\KL(M\|R)$ are infinite. Otherwise split $\log(dP/dR)$ into $\log(dP/dM)+\log(dM/dR)$, do the same for $Q$, and use $P+Q=2M$. The first two terms are finite, so this also justifies the identity when the remaining divergence is infinite. Taking $R=M$ proves \eqref{eq:projection}.

Finally, choose a measurable version of $\eta^*$ in $[0,1]$, setting it to $1/2$ on any exceptional $(P+Q)$-null set. For any common null law $R$, let $a=\int\eta^*\,dR$. Independence gives
\begin{equation}\label{eq:foura}
\E_{R\otimes R}e^*(X,Y)=4a(1-a)\le1.
\end{equation}
This calculation is conditional on the past as well. The oracle product is therefore valid under every conditional null law, and its asserted growth was already obtained from the strong law.
\end{proof}

\begin{remark}[Measurable-space formulation]
The dominating measure $P+Q$ allows arbitrary distributions, including singular and discrete laws. Equation~\eqref{eq:foura} uses independence within each null pair. Section~\ref{sec:mmd} treats exchangeable pairs and gives their corresponding information bound.
\end{remark}

\section{Predictive bets and aggregation}\label{sec:prediction}
\subsection{Exact reduction to logarithmic prediction loss}
Let $\eta_t:\X\to(0,1)$ be a measurable predictor selected using only $\F_{t-1}$. Both evaluations on pair $t$ use this same frozen predictor. Define
\begin{equation}\label{eq:predictive}
E_t=4\eta_t(X_t)\{1-\eta_t(Y_t)\},\qquad W_n^{\eta}=\prod_{t=1}^n E_t.
\end{equation}
The conditional version of \eqref{eq:foura} proves validity for any such learning rule. Let
\[
\ell_t(\eta)=-\log\eta(X_t)-\log\{1-\eta(Y_t)\},\qquad
\mathcal R(\eta)=\E\ell_t(\eta)
\]
for a fixed predictor. Then
\begin{equation}\label{eq:wealth-loss}
\log W_n^{\eta}=n\log4-\sum_{t=1}^n\ell_t(\eta_t).
\end{equation}

\begin{proposition}[Calibration identity]\label{prop:calibration}
For any fixed predictor $\eta:\X\to(0,1)$,
\begin{equation}\label{eq:calibration}
G(P,Q)-\E_A\log\bigl[4\eta(X)\{1-\eta(Y)\}\bigr]
=2\int\kl(\eta^*(z)\|\eta(z))\,dM(z).
\end{equation}
The same equality holds conditionally on $\F_{t-1}$ for $\eta_t$ under the i.i.d. alternative. If $H(u)=-u\log u-(1-u)\log(1-u)$, then
\begin{equation}\label{eq:bayes-risk}
\inf_\eta\mathcal R(\eta)=\mathcal R^*=2\int H(\eta^*)\,dM,
\qquad G(P,Q)=\log4-\mathcal R^*.
\end{equation}
\end{proposition}
\begin{proof}
Use $dP=2\eta^*dM$ and $dQ=2(1-\eta^*)dM$ in the two log-loss terms. Subtracting the Bayes loss gives twice the Bernoulli relative entropy, pointwise. Integration proves \eqref{eq:calibration}. The Bayes infimum over open-interval predictors follows by clipping $\eta^*$ away from zero and one and then taking the clipping level to zero. Conditional independence of the next pair from the past under the alternative proves the conditional statement.
\end{proof}

\subsection{Universal portfolios over the betting fraction}\label{sec:portfolio}
For any conditional base e-multiplier $E_t\ge0$, a fraction $\lambda\in[0,1]$ gives the valid multiplier $(1-\lambda)+\lambda E_t$. Its uniform-portfolio wealth is
\begin{equation}\label{eq:UP}
U_n=\int_0^1\prod_{t=1}^n\{(1-\lambda)+\lambda E_t\}\,d\lambda.
\end{equation}
A fixed mixture over valid processes is valid by conditional Tonelli. The next statement also supplies an exact implementation.

\begin{lemma}[Two-asset portfolio bound and recurrence]\label{lem:UP}
For nonnegative finite numbers $a_t,b_t$, define
\[
F_n(\lambda)=\prod_{t=1}^n\{(1-\lambda)a_t+\lambda b_t\}.
\]
Then
\begin{equation}\label{eq:UP-regret}
\int_0^1F_n(\lambda)\,d\lambda\ge\frac{\max_{\lambda\in[0,1]}F_n(\lambda)}{n+1}.
\end{equation}
Write $F_n(\lambda)=\sum_{j=0}^n d_{n,j}\binom nj\lambda^j(1-\lambda)^{n-j}$. With $d_{0,0}=1$ and nonexistent entries treated as zero,
\begin{align}
d_{n,j}&=\frac{n-j}{n}a_n d_{n-1,j}+\frac{j}{n}b_n d_{n-1,j-1},\label{eq:bernstein-update}\\
\int_0^1F_n(\lambda)\,d\lambda&=\frac{1}{n+1}\sum_{j=0}^n d_{n,j}.\label{eq:bernstein-int}
\end{align}
In particular, \eqref{eq:UP} loses at most $\log(n+1)$ log wealth relative to any fixed fraction, including full betting $\lambda=1$.
\end{lemma}
\begin{proof}
Multiplication by the next affine factor and the identities
$\binom{n-1}{j}/\binom nj=(n-j)/n$ and
$\binom{n-1}{j-1}/\binom nj=j/n$
give \eqref{eq:bernstein-update}; all coefficients are nonnegative. Each Bernstein basis polynomial integrates to $1/(n+1)$, proving \eqref{eq:bernstein-int}. The basis polynomials are nonnegative and sum to one, so $F_n(\lambda)\le\max_j d_{n,j}\le\sum_j d_{n,j}$. This proves \eqref{eq:UP-regret}.
\end{proof}
The actual predictable betting fraction in \eqref{eq:UP} is the wealth-weighted mean of $\lambda$ under the previous portfolio. For degree $m$ coefficients it equals
\begin{equation}\label{eq:posterior-stake}
\widehat\lambda_{m+1}
=\frac{\sum_{j=0}^m d_{m,j}(j+1)/(m+2)}{\sum_{j=0}^m d_{m,j}}.
\end{equation}
Equation~\eqref{eq:posterior-stake} expresses the universal portfolio as a self-financing predictable betting strategy.

\subsection{Countable model mixtures with finite computation}\label{sec:activation}
Index base learners by $j\ge1$. Fix deterministic activation times $s_j\uparrow\infty$, with finitely many $s_j\le n$ for every $n$, and prior weights $\pi_j>0$. Include a cash weight $\pi_0>0$, with $\pi_0+\sum_{j\ge1}\pi_j=1$. Expert $j$ starts from its prescribed initial state at $s_j$ and earns unit wealth before then. Let $B_{j,n}$ denote its cumulative base wealth on pairs $s_j,\ldots,n$, or one when $n<s_j$; let $U_{j,n}$ be its optional stake portfolio. Define
\begin{equation}\label{eq:outer}
W_n=\pi_0+\sum_{j\ge1}\pi_j U_{j,n}.
\end{equation}
Using $U_{j,n}=B_{j,n}$ omits fraction adaptation. The explicit choice used in the growing reference implementation is
\begin{equation}\label{eq:schedule}
s_j=2^{j-1},\qquad \pi_0=\frac12,\qquad
\pi_j=\frac{1}{2j(j+1)}.
\end{equation}
When $J$ experts are active, all inactive capital plus permanent cash equals $1/2+1/\{2(J+1)\}$. The inactive reserve retains the remaining prior mass exactly.

\begin{proposition}[Aggregation guarantees]\label{prop:aggregation}
The process \eqref{eq:outer} is a test supermartingale. For $n\ge s_j$ and $m=n-s_j+1$, the continuous-fraction version satisfies
\begin{equation}\label{eq:outer-oracle}
\log W_n\ge\log B_{j,n}-\log(m+1)-\log(1/\pi_j).
\end{equation}
Without fraction adaptation, omit $\log(m+1)$. With a fixed discrete fraction mixture assigning weight $w_1>0$ to $\lambda=1$, replace it by $\log(1/w_1)$. Delaying any fixed expert by a finite deterministic time has no first-order effect on its i.i.d. growth guarantee.
\end{proposition}
\begin{proof}
Each expert is a test supermartingale, including the constant pre-activation part. Conditional Tonelli proves validity of the countable sum. The inequality $W_n\ge\pi_j U_{j,n}$ and \cref{lem:UP} give \eqref{eq:outer-oracle}; the discrete case uses $U_{j,n}\ge w_1 B_{j,n}$. A fixed omitted prefix has finite length, and the strong law applies to the remaining i.i.d. tail.
\end{proof}

\section{A universal non-kernel construction}\label{sec:partition}
\subsection{Simultaneous-pair Laplace prediction}
Fix a finite measurable partition $\mathcal P=\{C_1,\ldots,C_d\}$. An expert starts with zero counts. Before its local pair $t$, let $A_{t-1,c}$ and $B_{t-1,c}$ count earlier $X$- and $Y$-observations in cell $C_c$. Predict
\begin{equation}\label{eq:laplace}
\eta_t(z)=\frac{A_{t-1,c}+1}{A_{t-1,c}+B_{t-1,c}+2},\qquad z\in C_c,
\end{equation}
and use \eqref{eq:predictive}. Only after computing both terms are both counts updated. This is Laplace smoothing with prediction frozen across each pair.

\begin{algorithm}[t]
\caption{One partition expert (local time)}\label{alg:partition}
\begin{algorithmic}[1]
\Require Cell map $c:\X\to\{1,\ldots,d\}$; optional stake wrapper
\State Initialize $A_c=B_c=0$ for all cells; $\log B=0$
\For{each new pair $(x,y)$}
\State $i\gets c(x)$; $j\gets c(y)$
\State $u\gets(A_i+1)/(A_i+B_i+2)$; $v\gets(A_j+1)/(A_j+B_j+2)$
\State $z\gets\log4+\log u+\log(1-v)$
\State Add $z$ to $\log B$, or pass $e^z$ to the stake wrapper
\State $A_i\gets A_i+1$; $B_j\gets B_j+1$ \Comment{after both predictions}
\EndFor
\end{algorithmic}
\end{algorithm}

\begin{theorem}[Finite-partition pathwise regret]\label{thm:partition-regret}
Let $B_{\mathcal P,m}$ be the base wealth after $m$ local pairs, and let $N_c=A_{m,c}+B_{m,c}$. For every vector $u=(u_1,\ldots,u_d)\in(0,1)^d$, viewed as a cellwise-constant predictor,
\begin{align}
\log B_{\mathcal P,m}
&\ge m\log4-\sum_{t=1}^m\ell_t(u)-\sum_{c=1}^d\log(N_c+1)\label{eq:hist-regret}\\
&\ge m\log4-\sum_{t=1}^m\ell_t(u)
-d\log(1+2m/d).\label{eq:hist-regret-coarse}
\end{align}
The inequalities hold for every data sequence.
\end{theorem}
\begin{proof}
Compare the product of frozen-pair label probabilities with a Bayesian sequential label predictor having independent uniform $\operatorname{Beta}(1,1)$ priors in the cells. Process the $X$ label first and the $Y$ label second within each pair. If the two observations lie in different cells, the Bayesian and frozen-pair probabilities for that pair agree. If both lie in a cell with previous beta parameters $a=A+1$, $b=B+1$ and $s=a+b$, the frozen product is
\[
\frac{a}{s}\frac{b}{s},
\]
whereas the Bayesian sequential product is $\frac{a}{s}\frac{b}{s+1}$. The former is at least the latter. Both methods have the same counts after the pair, so multiplication over pairs gives
\begin{equation}\label{eq:beta-dominance}
\prod_{t=1}^m\eta_t(X_t)\{1-\eta_t(Y_t)\}
\ge\prod_{c=1}^d\int_0^1v^{A_{m,c}}(1-v)^{B_{m,c}}\,dv.
\end{equation}
For nonnegative integers $a,b$ and $N=a+b$,
\[
\int_0^1v^a(1-v)^b\,dv
=\frac{1}{(N+1)\binom Na}
\ge\frac{u^a(1-u)^b}{N+1},
\]
because $\binom Na u^a(1-u)^b$ is a binomial probability, hence at most one. Substitution in \eqref{eq:beta-dominance}, multiplication by $4^m$, and logarithms prove \eqref{eq:hist-regret}. Finally $\sum_cN_c=2m$, so concavity of $\log$ gives $\sum_c\log(N_c+1)\le d\log(1+2m/d)$.
\end{proof}

\subsection{Information recovery through refining partitions}
On a standard Borel space choose increasing finite partitions $(\mathcal P_j)_{j\ge1}$ whose union generates $\B$. For example, enumerate a countable generator of $\B$ and take successive finite atom partitions. Empty cells can be discarded. Write $P_j,Q_j,M_j$ for the induced probability vectors and $G_j=2\JS(P_j,Q_j)$.

\begin{lemma}[Information recovered by a generating partition sequence]\label{lem:partition-info}
For increasing finite generating partitions, $G_j\uparrow G(P,Q)$.
\end{lemma}
\begin{proof}
Let $\B_j=\sigma(\mathcal P_j)$ and $\eta_j^*=\E_M[\eta^*\mid\B_j]$. On a cell of positive $M$ mass this equals $P(C)/(P(C)+Q(C))$. Thus
\[
G_j=\log4-2\int H(\eta_j^*)\,dM.
\]
The bounded martingale convergence theorem gives $\eta_j^*\to\eta^*$ $M$-almost surely, since $\sigma(\bigcup_j\B_j)=\B$. Entropy is continuous and bounded on $[0,1]$, so dominated convergence gives $G_j\to G$. Conditional Jensen applied to the concave function $H$ shows that the entropy integrals decrease with $j$, proving monotonicity.
\end{proof}

\begin{theorem}[Universal non-kernel log-optimality]\label{thm:partition-optimal}
Run one Laplace expert per partition $\mathcal P_j$, with any deterministic finite activation times satisfying \cref{sec:activation}, and combine them by \eqref{eq:outer}. Use either no stake wrapper, the exact continuous wrapper, or a fixed discrete wrapper with positive mass at full betting. Then for every pair of probability laws on the standard Borel space,
\begin{equation}\label{eq:partition-limit}
\frac{\log W_n}{n}\longrightarrow G(P,Q)\as
\end{equation}
The process is anytime valid under the conditional independent-pair null. With $\pi_0>0$, convergence also holds in $L^1$ and in expected normalized log wealth.
\end{theorem}
\begin{proof}
Validity has already been proved. Fix a partition and a cellwise-constant predictor $u$ with all coordinates in $(0,1)$. Its log loss is bounded. Combine \cref{thm:partition-regret,prop:aggregation}, divide by $n$, and apply the strong law to the tail beginning at $s_j$. All regret, delay, and prior terms are $o(n)$. Therefore
\[
\liminf_{n\to\infty}\frac{\log W_n}{n}
\ge\log4-\mathcal R(u)\as
\]
Taking a countable sequence of interior predictors tending to the cellwise Bayes predictor gives a lower bound $G_j$. Intersect these probability-one events over $j$ and use \cref{lem:partition-info}. The resulting lower bound is $G(P,Q)$, and \cref{thm:optimal} supplies the matching upper bound.

Every base multiplier is at most four; fractional multipliers are also at most four. Hence $\pi_0\le W_n\le4^n$. Consequently $n^{-1}\log W_n$ lies between $\log\pi_0$ and $\log4$ for $n\ge1$. Dominated convergence gives the asserted $L^1$ and expectation convergence.
\end{proof}

\paragraph{Implementation on non-Euclidean spaces.}
A countable measurable generator specifies the partition construction. On a countable domain, we form the $j$th partition by separating its first $j$ elements and collecting the remainder in one cell. On $\R^d$, a coordinatewise arctangent map into $(0,1)^d$ followed by dyadic partitions gives an explicit generator. The implementation accesses each observation through its cell map. A fixed partition with $d$ cells uses $O(\min\{d,m\})$ counter storage and constant counter-update time per pair, excluding cell evaluation. The schedule \eqref{eq:schedule} has $O(\log n)$ active partitions. Universal attainment uses the generating sequence; a finite collection has the corresponding fixed-partition regret guarantees.

\section{A log-optimal kernel construction on general spaces}\label{sec:kernel}
\subsection{Assumptions and the role of a constant feature}
Let $k:\X\times\X\to\R$ be a measurable positive semidefinite kernel with separable RKHS $\HH_k$ and $\sup_x k(x,x)\le1$. Its mean embedding is $\mu_k(R)=\int k(x,\cdot)\,dR(x)$. The integral exists as a Bochner integral: the feature map is measurable into a separable Hilbert space and is bounded. Call $k$ characteristic if $\mu_k(R)=\mu_k(S)$ implies $R=S$ for probability laws on $(\X,\B)$.

To include an intercept while preserving normalization, use
\begin{equation}\label{eq:augmented-kernel}
K(x,y)=\frac{1+k(x,y)}{2}.
\end{equation}
Its feature map can be realized as $(1,k(x,\cdot))/\sqrt2$ in $\R\oplus\HH_k$. In particular, $K(x,x)\le1$, its RKHS is separable, it contains the constant functions, and it contains $\HH_k$ as a set of functions. Characteristicness of $k$ is inherited by $K$.

\begin{lemma}[Density of the augmented RKHS]\label{lem:characteristic-density}
Under the preceding assumptions, if $k$ is characteristic, then $\HH_K$ is dense in $L^2(\mu)$ for every probability law $\mu$ on $(\X,\B)$.
\end{lemma}
\begin{proof}
All RKHS functions belong to $L^2(\mu)$ by bounded evaluation. Suppose $f\in L^2(\mu)$ is orthogonal to $\HH_K$. Since $1\in\HH_K$, $\int f\,d\mu=0$. Since $\HH_k\subseteq\HH_K$,
\[
\int f(x)h(x)\,d\mu(x)=0\quad\text{for every }h\in\HH_k.
\]
The finite signed measure $\nu=f\mu$ therefore has total mass zero and mean embedding zero in $\HH_k$. If $\nu\ne0$, its positive and negative parts have the same positive finite mass $c$. Their normalized versions $R=\nu^+/c$ and $S=\nu^-/c$ are distinct probability measures with $\mu_k(R)=\mu_k(S)$, contradicting characteristicness. Hence $f=0$ $\mu$-almost surely. A proper closed subspace of a Hilbert space has a nonzero orthogonal complement; consequently the closure of $\HH_K$ in $L^2(\mu)$ is all of $L^2(\mu)$.
\end{proof}

\begin{remark}
Characteristicness gives injectivity on probability measures. The constant feature supplies the total-mass constraint needed to extend the argument to signed measures and obtain function approximation. Alternatively, any bounded measurable separable RKHS already dense in $L^2(M)$ can be used directly for the fixed alternative, without a characteristicness assumption.
\end{remark}

\subsection{Approximation of the Bayes risk by logistic scores}
For $h\in\HH_K$, put $\eta_h=\sig\circ h$, where $\sig(v)=1/(1+e^{-v})$, and define $\spplus(v)=\log(1+e^v)$. The loss and multiplier are
\begin{align}
\ell_t(h)&=\spplus(-h(X_t))+\spplus(h(Y_t)),\label{eq:kernel-loss}\\
E_t(h)&=4\sig(h(X_t))\sig(-h(Y_t)),\label{eq:kernel-factor}\\
\mathcal R(h)&=\E_A\ell_t(h).
\end{align}
The extended-real Bayes score is $h^*=\log\{\eta^*/(1-\eta^*)\}$; it can be infinite on sets where the sources are perfectly distinguishable.

\begin{lemma}[Bayes-risk approximation in the RKHS]\label{lem:logistic-approx}
If $\HH_K$ is dense in $L^2(M)$, then
\begin{equation}\label{eq:rkhs-inf}
\inf_{h\in\HH_K}\mathcal R(h)=\mathcal R^*.
\end{equation}
More explicitly, for $B>0$, put $h_B^*=\clip(h^*,-B,B)$. Whenever $\norm{h-h_B^*}_{L^2(M)}\le\delta$,
\begin{equation}\label{eq:approx-rate}
\mathcal R(h)-\mathcal R^*\le2\log(1+e^{-B})+2\delta.
\end{equation}
\end{lemma}
\begin{proof}
The conditional logistic risk at class probability $u$ is
$r_u(v)=u\spplus(-v)+(1-u)\spplus(v)$, with derivative $\sig(v)-u$. Hence $r_u$ is 1-Lipschitz, and
\[
|\mathcal R(h)-\mathcal R(g)|\le2\int|h-g|\,dM\le2\norm{h-g}_{L^2(M)}.
\]
Let $a=\sig(-B)$. Clipping the optimal probability to $[a,1-a]$ incurs Bernoulli KL at most $-\log(1-a)=\log(1+e^{-B})$. To see this, the loss is zero when $u\in[a,1-a]$; on $[0,a]$, convexity of $u\mapsto\kl(u\|a)$ places its maximum at an endpoint, with values $-\log(1-a)$ and zero. The other interval is symmetric. The calibration identity thus gives $\mathcal R(h_B^*)-\mathcal R^*\le2\log(1+e^{-B})$. Add the Lipschitz bound to obtain \eqref{eq:approx-rate}. Since $h_B^*$ is bounded, it lies in $L^2(M)$ and can be approximated there by RKHS functions. First take $\delta\downarrow0$, then $B\to\infty$.
\end{proof}

\subsection{Projected online learning with an explicit regret bound}
Fix a radius $R>0$. At local time $t$, use $h_t\in\HH_K$ with $\norm{h_t}_K\le R$ to score the next pair. Start from $h_1=0$. After scoring, compute
\begin{equation}\label{eq:gradient}
g_t=\{\sig(h_t(X_t))-1\}K(X_t,\cdot)
       +\sig(h_t(Y_t))K(Y_t,\cdot)
\end{equation}
and update
\begin{equation}\label{eq:OGD}
h_{t+1}=\Pi_R(h_t-\gamma_t g_t),\qquad
\gamma_t=\frac{R}{2\sqrt t},\qquad
\Pi_R(f)=f\min\{1,R/\norm f_K\}.
\end{equation}
The expression for $\Pi_R(0)$ is interpreted as zero.

\begin{theorem}[Kernel learning regret]\label{thm:kernel-regret}
For every data sequence and every $h\in\HH_K$ with $\norm h_K\le R$, the iterates \eqref{eq:OGD} satisfy
\begin{equation}\label{eq:kernel-regret}
\sum_{t=1}^m\{\ell_t(h_t)-\ell_t(h)\}\le6R\sqrt m.
\end{equation}
Consequently their base wealth satisfies
\begin{equation}\label{eq:kernel-wealth-regret}
\log B_{K,R,m}\ge m\log4-\sum_{t=1}^m\ell_t(h)-6R\sqrt m.
\end{equation}
\end{theorem}
\begin{proof}
The loss is convex and Fr\'echet differentiable, with gradient \eqref{eq:gradient}. Each coefficient has magnitude at most one and $\norm{K(x,\cdot)}_K\le1$, so $\norm{g_t}_K\le2$. Nonexpansiveness of projection gives
\[
\norm{h_{t+1}-h}_K^2
\le\norm{h_t-h}_K^2-2\gamma_t\inner{g_t}{h_t-h}_K
 +\gamma_t^2\norm{g_t}_K^2.
\]
Convexity implies $\ell_t(h_t)-\ell_t(h)\le\inner{g_t}{h_t-h}_K$. Writing $D_t=\norm{h_t-h}_K^2\le4R^2$ and summing,
\begin{align*}
\sum_{t=1}^m\{\ell_t(h_t)-\ell_t(h)\}
&\le\frac{D_1}{2\gamma_1}
 +\sum_{t=2}^mD_t\left(\frac1{2\gamma_t}-\frac1{2\gamma_{t-1}}\right)
 +\frac12\sum_{t=1}^m\gamma_t\norm{g_t}_K^2\\
&\le\frac{4R^2}{2\gamma_m}+2\sum_{t=1}^m\gamma_t
\le4R\sqrt m+2R\sqrt m.
\end{align*}
We used decreasing $\gamma_t$ and $\sum_{t=1}^m t^{-1/2}\le2\sqrt m$. This proves \eqref{eq:kernel-regret}; \eqref{eq:kernel-wealth-regret} follows from \eqref{eq:wealth-loss}.
\end{proof}

\subsection{Mixing norm budgets and kernels}
Take a finite or countable bank $(K_b)$ of bounded measurable kernels with separable RKHSs, all normalized by $K_b(x,x)\le1$. For each kernel include unbounded radii, for example $R_r=2^{r-1}$, $r\ge1$. Enumerate the pairs $(b,r)$ so every pair receives a finite index, activation time, and positive prior weight. For a finite bank of size $B_0$, one may take $j=B_0(r-1)+b$. For a countably infinite bank, use diagonal enumeration. Run \eqref{eq:OGD} independently for each expert and aggregate as in \eqref{eq:outer}.

\begin{theorem}[Universal kernel log-optimality]\label{thm:kernel-optimal}
Suppose at least one kernel $K_b$ in the bank has RKHS dense in $L^2(M)$. For the resulting kernel wealth, with any of the three stake options of \cref{thm:partition-optimal},
\begin{equation}\label{eq:kernel-limit}
\frac{\log W_n}{n}\longrightarrow G(P,Q)\as
\end{equation}
If the bank contains the augmentation \eqref{eq:augmented-kernel} of a bounded measurable characteristic kernel with separable RKHS, the conclusion holds for every fixed $(P,Q)$ on the space. Anytime validity holds without any approximation assumption. With positive cash weight, convergence is also in $L^1$ and expected normalized log wealth.

In particular, for an active expert $j=(b,r)$, $m=n-s_j+1$, and $\norm h_{K_b}\le R_r$, the continuous-portfolio version obeys
\begin{equation}\label{eq:master-kernel}
\log W_n\ge m\log4-\sum_{t=s_j}^n\ell_t(h)
-6R_r\sqrt m-\log(m+1)-\log(1/\pi_j).
\end{equation}
\end{theorem}
\begin{proof}
The finite-sample inequality is \cref{thm:kernel-regret,prop:aggregation}. Each fixed comparator $h$ is bounded pointwise by its RKHS norm, so $\ell_t(h)$ is bounded. Divide \eqref{eq:master-kernel} by $n$ and apply the strong law to the fixed comparator. For every such $h$,
\[
\liminf_n n^{-1}\log W_n\ge\log4-\mathcal R(h)\as.
\]
Choose a countable sequence of comparator functions in the dense RKHS with risks tending to $\mathcal R^*$, as supplied by \cref{lem:logistic-approx}. Each belongs to some finite radius and therefore has an expert assigned positive prior mass and finite activation time. Intersecting the corresponding probability-one events gives the lower bound $G(P,Q)$. The upper bound is \cref{thm:optimal}. \Cref{lem:characteristic-density} gives the distribution-universal sufficient condition. The other stake options merely change an $o(n)$ term. The cash lower bound and multiplier upper bound give $L^1$ convergence exactly as in \cref{thm:partition-optimal}.
\end{proof}

\paragraph{Adaptation over norms, kernels, and betting fractions.}
The RKHS radius controls the expressiveness and magnitude of the score, the kernel bandwidth determines its geometry, and the betting fraction controls capital allocation. One augmented characteristic kernel with unbounded norm budgets suffices for first-order optimality. Mixtures over bandwidths provide additional adaptation. A finite set of betting fractions containing $1$ preserves the optimal rate.

\subsection{Concrete general-space kernels}\label{sec:general-spaces}
\begin{example}[Countable structured observations]
For a countable set of strings, sequences, or canonically represented finite graphs, the equality kernel $k_\delta(x,y)=\one\{x=y\}$ has separable RKHS $\ell^2(\X)$ and is characteristic: its mean embedding lists the point probabilities. Thus $(1+k_\delta)/2$ satisfies the theorem. For a bounded normalized domain-specific positive semidefinite kernel $k_0$, the mixture
\[
k_\rho=(1-\rho)k_0+\rho k_\delta,\qquad 0<\rho\le1,
\]
is characteristic as well, because a zero squared mean-embedding distance forces the equality-kernel distance to be zero. The domain kernel shares information between related objects, while the equality component ensures characteristicness on the countable domain.
\end{example}

\begin{example}[An existence construction on every standard Borel space]
Let $\mathcal P_0=\{\X\}$ and let $(\mathcal P_j)_{j\ge1}$ be increasing finite generating partitions. Define
\begin{equation}\label{eq:partition-kernel}
K_{\mathrm{part}}(x,y)
=\sum_{j=0}^\infty2^{-(j+1)}
 \sum_{C\in\mathcal P_j}\one_C(x)\one_C(y).
\end{equation}
This is a measurable positive semidefinite kernel with diagonal one and a countable-dimensional feature representation in $\ell^2$, hence a separable RKHS. Its mean embedding contains a positive multiple of every partition-cell probability. Equality of embeddings implies equality on the algebra generated by these cells and therefore equality of probability laws on $\B$. It is characteristic and already contains the constant feature from $j=0$. This construction supplies a kernel satisfying the theorem on every standard Borel space. Its computational cost depends on the representation of the partitions and evaluation of the series.
\end{example}

\paragraph{Geometry and representation.}
Both constructions operate on general objects through cell maps or kernel evaluations. A kernel's RKHS determines the available score functions. Under the approximation assumptions of \cref{thm:kernel-optimal}, this class attains the full information bound; for other kernels, \eqref{eq:master-kernel} gives a comparator guarantee within the chosen RKHS.

\section{MMD witnesses and the exchangeable-pair benchmark}\label{sec:mmd}
\subsection{Information bound for additive witnesses}
The betting framework of \citet{shekhar2024betting} includes predictable multipliers
\begin{equation}\label{eq:additive}
1+\lambda_t\{g_t(X_t)-g_t(Y_t)\},\qquad
\norm{g_t}_\infty\le\frac12,\quad |\lambda_t|\le1.
\end{equation}
These multipliers satisfy $e(x,y)+e(y,x)=2$, are at most two, and have mean exactly one under any exchangeable pair law. Their validity under within-pair dependence leads to the exchangeable-pair information bound.

\begin{theorem}[Exchangeable-pair optimum]\label{thm:swap}
Let $A=P\otimes Q$, $A^{\mathrm{sw}}=Q\otimes P$, and $S=(A+A^{\mathrm{sw}})/2$. Among processes valid under every i.i.d. exchangeable pair law, the maximal expected log-growth per pair and attainable asymptotic log-growth rate are
\begin{equation}\label{eq:swap-rate}
G_{\mathrm{sw}}(P,Q)=\KL(A\|S)=\JS(A,A^{\mathrm{sw}})\le\log2.
\end{equation}
An oracle multiplier is
\begin{equation}\label{eq:swap-oracle}
e^*_{\mathrm{sw}}(x,y)=
\frac{2p(x)q(y)}{p(x)q(y)+q(x)p(y)},
\end{equation}
where $p=dP/d(P+Q)$, $q=dQ/d(P+Q)$, and the value is set to one when the denominator is zero. Moreover,
\begin{equation}\label{eq:rate-gap}
G(P,Q)=G_{\mathrm{sw}}(P,Q)+\KL(S\|M\otimes M).
\end{equation}
\end{theorem}
\begin{proof}
The multiplier in \eqref{eq:swap-oracle} is $dA/dS$ and satisfies $e(x,y)+e(y,x)=2$ everywhere. Its expectation under every exchangeable law is consequently one. Its logarithm is integrable under $A$, by the same bounded-density-ratio argument as in \cref{thm:optimal}, now using the bound two. The strong law gives its rate $\KL(A\|S)$. Since $S$ itself is an exchangeable null, \cref{lem:comparison} and Borel--Cantelli give the matching upper bound for every valid process. Swapping coordinates shows $\KL(A\|S)=\KL(A^{\mathrm{sw}}\|S)$, establishing the JS identity.

Finally $dS/d(M\otimes M)$ is symmetric in the coordinates. Decompose $\log(dA/d(M\otimes M))$ through $S$. The expectation under $A$ of the symmetric remaining log ratio equals its expectation under $S$, yielding \eqref{eq:rate-gap}. All divergences in this identity are finite; the first two are bounded by $\log4$ and $\log2$, respectively, and $dS/d(M\otimes M)$ is bounded.
\end{proof}
For mutually singular sources, $G=\log4$ while $G_{\mathrm{sw}}=\log2$. Every additive wealth in \eqref{eq:additive} is at most $2^n$, and every fixed wealth mixture of them is also at most $2^n$. The same $2^n$ bound therefore applies to mixtures over kernels and betting fractions within this payoff family.

The role of independence in the product multiplier is illustrated by an exchangeable pair equal to $(0,1)$ or $(1,0)$ with probability $1/2$ each. For $\eta(0)=1$ and $\eta(1)=0$, its mean product multiplier is two. Open-interval approximations to this predictor also give a mean greater than one.

\subsection{Full-range fraction adaptation for MMD witnesses}
For any predictable $V_{k,t}=g_{k,t}(X_t)-g_{k,t}(Y_t)\in[-1,1]$, use
\begin{equation}\label{eq:mmd-UP}
U_{k,n}=\frac12\int_{-1}^1\prod_{t=1}^n(1+\lambda V_{k,t})\,d\lambda.
\end{equation}
Taking $a_t=1-V_{k,t}$, $b_t=1+V_{k,t}$ in \cref{lem:UP} proves
\begin{equation}\label{eq:mmd-regret}
\log U_{k,n}\ge\max_{|\lambda|\le1}\sum_{t=1}^n\log(1+\lambda V_{k,t})-\log(n+1).
\end{equation}
A prior-weighted mixture over kernels adds only $\log(1/\pi_k)$ to the comparator penalty. Optimizing over the full fraction interval recovers the fixed-witness oracle comparator; \citet{waudbysmith2026} analyze the losses caused by restricting this interval.

Suppose $k(x,x)\le1$, its RKHS is separable, and $\Delta=\mu_k(P)-\mu_k(Q)\ne0$. The empirical normalized witness is
\[
g_{k,t}=\frac{\widehat\Delta_{t-1}}{2\norm{\widehat\Delta_{t-1}}_k},
\qquad g_k^*=\frac{\Delta}{2\norm\Delta_k},
\]
with $g_{k,t}=0$ when its denominator vanishes. The Hilbert-space strong law gives RKHS convergence and therefore uniform convergence $g_{k,t}\to g_k^*$. For every fixed $|\lambda|<1$, logarithms in \eqref{eq:mmd-regret} have denominators bounded below by $1-|\lambda|$. Their empirical averages consequently converge to
\[
F_k(\lambda)=\E\log\{1+\lambda(g_k^*(X)-g_k^*(Y))\}.
\]
Using a countable dense set of interior fractions in \eqref{eq:mmd-regret} gives
\begin{equation}\label{eq:mmd-oracle-lower}
\liminf_n n^{-1}\log U_{k,n}\ge\sup_{|\lambda|\le1}F_k(\lambda)\as.
\end{equation}
To include endpoints, use concavity: $F_k(r\lambda)\ge rF_k(\lambda)$ for $0<r<1$, since $F_k(0)=0$. An endpoint with value $-\infty$ cannot improve the supremum. Thus interior fractions approximate the full fixed-witness optimum.

The MMD witness maximizes a mean difference. Expected log wealth can instead be optimized within the additive family by learning $f$ with the convex loss $-\log\{1+f(X)-f(Y)\}$. On the RKHS ball of radius $(1-\varepsilon)/2$, the multiplier is at least $\varepsilon$ and the gradient norm is at most $2/\varepsilon$. The projected-gradient proof in \cref{thm:kernel-regret}, with this gradient bound, gives sublinear regret for every fixed $\varepsilon>0$. Mixing over radii, kernels, and $\varepsilon\downarrow0$ then competes with that additive function class. The resulting strategy is subject to the exchangeable-pair bound. Appendix~\ref{app:swap} gives a kernel method attaining that bound.

\section{Practical implementation}\label{sec:implementation}
\subsection{An exact kernel update without a Gram matrix}
A local kernel expert stores $h_t=\sum_{i=1}^{N_t}c_iK(z_i,\cdot)$ and $q_t=\norm{h_t}_K^2$. Evaluate $u=h_t(x)$ and $v=h_t(y)$ on the new pair before changing any coefficient. With $a=\sig(u)-1$, $b=\sig(v)$, and $\gamma=R/(2\sqrt t)$, the unprojected update appends coefficients $-\gamma a$ at $x$ and $-\gamma b$ at $y$. Its squared norm is
\begin{equation}\label{eq:norm-update}
\widetilde q=q_t-2\gamma(au+bv)
+\gamma^2\{a^2K(x,x)+2abK(x,y)+b^2K(y,y)\}.
\end{equation}
Multiply all resulting coefficients by $\rho=\min\{1,R/\sqrt{\widetilde q}\}$ and store $q_{t+1}=\rho^2\widetilde q$. This is exactly \eqref{eq:OGD} in exact arithmetic. Repeated support points may be merged to reduce storage. The code uses a generic scalar kernel callback, so observations may be arbitrary Python objects supported by that callback.

\begin{algorithm}[t]
\caption{One projected kernel logistic expert (local time)}\label{alg:kernel}
\begin{algorithmic}[1]
\Require Normalized kernel $K$, radius $R>0$, optional stake wrapper
\State Initialize empty support list, empty coefficients, $q=0$, $\log B=0$
\For{$t=1,2,\ldots$}
\State Observe pair $(x,y)$; compute $u=\sum_i c_iK(z_i,x)$, $v=\sum_i c_iK(z_i,y)$
\State $z\gets\log4-\spplus(-u)-\spplus(v)$
\State Add $z$ to base log wealth, or pass it to the stake wrapper
\State $a\gets\sig(u)-1$; $b\gets\sig(v)$; $\gamma\gets R/(2\sqrt t)$
\State Compute $\widetilde q$ by \eqref{eq:norm-update}
\State Append $(x,-\gamma a)$ and $(y,-\gamma b)$ to the expansion
\State $\rho\gets\min\{1,R/\sqrt{\widetilde q}\}$, with $\rho=1$ if $\widetilde q=0$
\State Multiply all coefficients by $\rho$; $q\gets\rho^2\widetilde q$
\EndFor
\end{algorithmic}
\end{algorithm}

With $m$ previous local pairs, the expansion has at most $2m$ terms. Each next pair uses $O(m)$ kernel evaluations and arithmetic, and $O(m)$ storage, assuming a kernel call has constant cost. The update maintains the squared norm directly, with total cost $O(n^2)$ per expert over $n$ pairs. Sparse dictionaries and random-feature approximations can reduce this cost. The growth theorem extends to an approximate learner whenever its cumulative excess predictive log loss relative to the exact learner is $o(n)$.

\subsection{Stake implementation and numerical stability}
The exact uniform portfolio uses \eqref{eq:bernstein-update} in logarithmic coordinates. If $D_{m,j}=\log d_{m,j}$, each update is a log-sum-exp of its two nonnegative summands, with missing terms represented by $-\infty$. This requires $O(m)$ time and memory at local pair $m$, hence $O(n^2)$ total time. The base product case uses $(a,b)=(1,E)$; the additive case uses $(a,b)=(1-V,1+V)$. The scalar implementation supports both.

A lower-cost alternative uses a fixed finite grid $\Lambda\subseteq[0,1]$ containing $1$, for example $\{0,1/2,1\}$, with fixed positive weights. We maintain one log wealth per fraction and combine them by log-sum-exp. This costs $O(|\Lambda|)$ per round and retains \cref{thm:partition-optimal,thm:kernel-optimal}, since the full-bet component receives positive prior mass. Its comparator penalty is the negative log prior weight of the chosen grid point. Fixed quadrature rules with positive normalized weights give the same finite-mixture guarantee.

We evaluate $\log\sig(u)=-\spplus(-u)$ using the stable identity $\spplus(u)=\logaddexp(0,u)$. Laplace log probabilities are computed as differences of log counts. Wealth, prior weights, inactive reserves, and stopping thresholds are stored in log space. The implementation assumes a positive semidefinite, normalized kernel and checks the norm update \eqref{eq:norm-update} for numerical inconsistencies.

\subsection{Growing-bank bookkeeping}
At time $t$, activate every expert with $s_j=t$ \emph{before} observing the pair, initializing its predictor and stake state without using that pair. Subtract its predetermined capital from the inactive reserve, leaving total wealth unchanged. Score every active expert using its old predictor; then update its learner for the next round. The reported wealth is
\begin{equation}\label{eq:reserve-wealth}
W_t=\pi_0+\sum_{j:s_j>t}\pi_j+\sum_{j:s_j\le t}\pi_jU_{j,t}.
\end{equation}
The reported wealth retains the original prior weights. With \eqref{eq:schedule}, the number of active experts is $1+\lfloor\log_2t\rfloor$. A slower-growing activation time or an application-specific finite leading bank can reduce delays for useful models, provided each comparator is eventually activated with positive capital.



\subsection{Adaptive representations and neural fine-tuning}\label{sec:finetuning}
The predictor in \eqref{eq:predictive} may be a neural network whose representation and classification head both learn from the streams. Write
\begin{equation}\label{eq:neural-predictor}
 \eta_t(z)=\sig\!\left(g_{\omega_{t-1}}\bigl(\phi_{\theta_{t-1}}(z)\bigr)\right),
 \qquad E_t=4\eta_t(X_t)\{1-\eta_t(Y_t)\},
\end{equation}
where the encoder parameters $\theta_{t-1}$, head parameters $\omega_{t-1}$, and all preprocessing are $\F_{t-1}$-measurable. A pretrained initialization may be obtained from data independent of the testing streams. Observations from $P$ receive source label one and those from $Q$ source label zero. These labels are available from the sampling protocol, even when semantic image labels are unknown.

The model scores both members of a pair before either enters its training buffer. After recording $E_t$, any update based on the observations through time $t$ may change the head, selected encoder layers, or the entire network. Earlier observations may be reused for subsequent gradient steps while their recorded wealth increments remain fixed. Updates can occur after every pair or at less frequent checkpoints. This predict-then-train protocol is also used in deep anytime-valid testing \citep{pandeva2024deep}. For the product bet, ordinary source cross-entropy gives exactly the wealth identity \eqref{eq:wealth-loss}.

Conditionally on the past, put $a_t=\int\eta_t\,dR_t$ under the null. Then
\[
 \E[E_t\mid\F_{t-1}]=4a_t(1-a_t)\le1.
\]
Thus fine-tuning preserves the test-supermartingale property for every predictable training rule. Scoring uses the preceding normalization statistics and a pointwise predictor; evaluation-mode network layers implement this requirement. Training quality enters through \cref{prop:calibration}. For two predictable learning rules, their unfractionated product wealths satisfy
\begin{equation}\label{eq:neural-loss-comparison}
 \log W_n^{\mathrm{tuned}}-\log W_n^{\mathrm{frozen}}
 =\sum_{t=1}^n\{\ell_t^{\mathrm{frozen}}-\ell_t^{\mathrm{tuned}}\}.
\end{equation}
Every reduction in prequential logarithmic loss increases log wealth by the same amount. \Cref{prop:neural-growth} gives a sufficient accumulated-loss condition for almost-sure attainment of $G(P,Q)$ by an arbitrary learner.

\paragraph{Learning the representation.}
A fixed feature map $\phi$ restricts the available information to $2\JS(P_\phi,Q_\phi)$, where $P_\phi,Q_\phi$ are the induced laws. The gap to $2\JS(P,Q)$ is zero precisely when the source posterior is measurable with respect to the features. Fine-tuning can recover distinctions suppressed by a generic representation and can make them easier for a simple head to express. \Cref{prop:representation-gap} decomposes the total growth shortfall into information lost by the representation and predictive error within that representation. Estimation and calibration costs also change during adaptation, so both prequential loss and early rejection times are relevant performance measures.

\paragraph{Retaining a frozen-model component.}
For a finite set of temperatures $T>0$ and fractions $\lambda\in[0,1]$, each arm can maintain a fixed mixture of the valid multipliers
\[
 1-\lambda+4\lambda\sig(h_t(X_t)/T)\{1-\sig(h_t(Y_t)/T)\}.
\]
Here $h_t$ is the arm's source logit, fixed before pair $t$. Let $W^F$ and $W^T$ be the resulting frozen-encoder and fine-tuned wealth processes. The equal mixture satisfies
\begin{equation}\label{eq:frozen-tuned-mixture}
 W_t^{\mathrm{mix}}=\tfrac12W_t^F+\tfrac12W_t^T,
 \qquad
 \log W_t^{\mathrm{mix}}\ge
 \max\{\log W_t^F,\log W_t^T\}-\log2.
\end{equation}
The two arms may share their training observations. Their mixture remains valid and its wealth stays within a factor two of either constituent at every time. A positive-weight universal component from \cref{sec:partition,sec:kernel} additionally supplies the corresponding universal growth guarantee. The all-allocation version of this combination is given after \cref{cor:all-ratios}.

\begin{algorithm}[tbp]
\caption{Predict-then-train image betting}\label{alg:finetuning}
\begin{algorithmic}[1]
\Require Independent pretrained encoder; update schedule; temperatures and fractions; fixed mixture weights
\State Initialize frozen and trainable encoders from the same pretrained checkpoint
\State Set both source predictions to $1/2$ and all wealth accounts to one
\State Initialize an empty source-labeled replay buffer
\For{each newly observed pair $(x,y)$}
 \State Evaluate both arms' existing source predictors on $x$ and $y$
 \State Update each temperature/fraction account using \eqref{eq:neural-predictor}
 \State Report the fixed mixture of the account wealths; compare it with $1/\alpha$
 \State Append $(x,1)$ and $(y,0)$ to the replay buffer
 \If{an update is scheduled}
  \State Refit the frozen arm's source head on preceding source-labeled observations
  \State Fine-tune the trainable arm's selected layers using only the replay buffer
  \State Refit its source head if specified; fix both predictors for subsequent scoring
 \EndIf
\EndFor
\end{algorithmic}
\end{algorithm}
A warm-up period is implemented by neutral predictions $\eta_t\equiv1/2$. For unequal prospective blocks, the same learned classifier supplies the weighted product bet of \cref{thm:imb-oracle}, with updates after each scored block. The sample-size and reporting protocols are described in \cref{sec:imbalance}.

\section{Consequences for rejection times}\label{sec:stopping}
The comparator bounds yield power one and first-order optimality of both almost-sure and expected rejection times in the high-confidence limit.

\begin{theorem}[First-order rejection-time efficiency]\label{thm:stopping}
Fix $P\ne Q$, and use either universal construction with a fixed prior, fixed activation schedule, and one of the stated stake wrappers. Then $\PP_{P,Q}(\tau_\alpha<\infty)=1$, $\E_{P,Q}\tau_\alpha<\infty$ for every $0<\alpha<1$, and
\begin{equation}\label{eq:stopping-limit}
\frac{\tau_\alpha}{\log(1/\alpha)}\longrightarrow\frac1{G(P,Q)}\as,
\qquad
\frac{\E_{P,Q}\tau_\alpha}{\log(1/\alpha)}\longrightarrow\frac1{G(P,Q)}
\quad(\alpha\downarrow0).
\end{equation}
For any fixed process satisfying the deterministic-time validity conditions of \cref{thm:optimal}, the corresponding $\liminf$ of each ratio is at least $1/G(P,Q)$, allowing infinite rejection times and expectations.
\end{theorem}
\begin{proof}
Write $b=\log(1/\alpha)$ and $\tau_b=\inf\{n:\log W_n\ge b\}$. Since $G>0$, almost-sure log-growth gives power one. For any $\varepsilon\in(0,G)$, eventually
\[
(G-\varepsilon)n\le\log W_n\le(G+\varepsilon)n
\]
on the probability-one convergence event. The maximum over any finite prefix is finite, so it cannot cross arbitrarily large thresholds. These bounds imply $\tau_b/b\to1/G$. For an arbitrary valid process, only the upper bound on eventual log growth is needed to obtain $\liminf_b\tau_b/b\ge1/G$. Fatou's lemma gives the corresponding expected $\liminf$ along every sequence $b\to\infty$.

It remains to prove the expected upper bound for our constructions. Choose a fixed comparator with bounded log multiplier $Z_t$, mean $\gamma>0$, and range width $L<\infty$. In the partition construction use any clipped interior cellwise predictor; in the kernel construction use a fixed finite-norm score. Their approximation lemmas allow $\gamma$ to be arbitrarily close to $G$. For some fixed expert activation time $s$ and a deterministic function $r(m)=o(m)$, the finite-sample regret bounds imply
\[
\log W_{s-1+m}\ge\sum_{t=s}^{s-1+m}Z_t-r(m).
\]
The prior and stake penalties are included in $r$. Fix $0<\varepsilon<\gamma/2$ and take $m_0$ so $r(m)\le\varepsilon m$ for $m\ge m_0$. If also $m\ge b/(\gamma-2\varepsilon)$, then
\begin{align*}
\PP(\tau_b>s-1+m)
&\le\PP(\log W_{s-1+m}<b)\\
&\le\PP\left(\sum_{t=s}^{s-1+m}(Z_t-\gamma)<-\varepsilon m\right)
\le\exp(-2\varepsilon^2m/L^2),
\end{align*}
by Hoeffding's inequality. If $L=0$, the last event is impossible and the same conclusion follows directly. Summing the tail probabilities yields a finite expectation and
\[
\limsup_{b\to\infty}\frac{\E\tau_b}{b}\le\frac1{\gamma-2\varepsilon}.
\]
Let $\varepsilon\downarrow0$ and choose comparators with $\gamma\uparrow G$. This proves the matching expected upper bound for each fixed alternative.
\end{proof}

\section{Unequal and strongly imbalanced sample sizes}\label{sec:imbalance}
We now observe independent streams $X_1,X_2,\ldots\sim P$ and $Y_1,Y_2,\ldots\sim Q$, with respective sample sizes $a$ and $b$. Under the null, both streams have a fixed common law $R$. The stream sampled next may be selected using the observations collected so far. Put
\[
N=a+b,\quad m=\min(a,b),\quad w=a/N,\quad v=b/N,\quad M_w=wP+vQ.
\]
Strong imbalance calls for learning and validity costs of order $o(m)$. We obtain such bounds by universal coding and by training a kernel score on a vanishing fraction of the smaller sample.

\subsection{The exact information bound and its one-sample limit}
Define the weighted information
\begin{equation}\label{eq:imb-information}
I_{a,b}(P,Q)=a\KL(P\|M_w)+b\KL(Q\|M_w).
\end{equation}
For $0<\rho<1$, also write
\[
G_\rho(P,Q)=\rho\KL(P\|\rho P+(1-\rho)Q)
 +(1-\rho)\KL(Q\|\rho P+(1-\rho)Q).
\]
The following result gives a Radon--Nikodym formulation of the unequal-block oracle of \citet{turner2024generic}. The subsequent constructions learn its optimal rate simultaneously for all unknown laws in the stated general-space models.

\begin{theorem}[Finite-sample optimum]\label{thm:imb-oracle}
For every nonnegative $E_{a,b}$ satisfying $\E_{R^{a+b}}E_{a,b}\le1$ for all common laws $R$,
\begin{equation}\label{eq:imb-log-upper}
\E_{P^a\otimes Q^b}\log E_{a,b}\le I_{a,b}(P,Q)
 =\inf_R\{a\KL(P\|R)+b\KL(Q\|R)\}.
\end{equation}
Equality is attained by
\begin{equation}\label{eq:imb-oracle-e}
E^*_{a,b}=\prod_{i=1}^a\frac{dP}{dM_w}(X_i)
                 \prod_{j=1}^b\frac{dQ}{dM_w}(Y_j).
\end{equation}
More generally, every fixed bounded measurable score $h$ gives an e-variable
\begin{align}
E_{a,b}(h)&=\exp F_{a,b}(h),\label{eq:imb-score-e}\\
F_{a,b}(h)&=\sum_{j=1}^b h(Y_j)
 -\sum_{z\in(X_1,\ldots,X_a,Y_1,\ldots,Y_b)}
       \log\{w+v e^{h(z)}\}.\label{eq:imb-score-F}
\end{align}
The sum is over observations with multiplicity. A score trained on an independent pilot also gives a conditionally valid e-variable on a fresh testing block, with $w,v$ determined by that block.
\end{theorem}
\begin{proof}
For any $R$, the entropy variational inequality applied to $P^a\otimes Q^b$ and $R^{a+b}$ gives the upper bound by $a\KL(P\|R)+b\KL(Q\|R)$. Decomposition of the log densities gives
\[
a\KL(P\|R)+b\KL(Q\|R)
 =I_{a,b}(P,Q)+N\KL(M_w\|R),
\]
including infinite divergences, so the minimum occurs at $R=M_w$.
Let $f_0=dP/dM_w$ and $f_1=dQ/dM_w$, choosing versions such that $wf_0+vf_1=1$ everywhere (set both to one outside a common supporting set). Under a common law $R$, put $u=\E_R f_0$ and $z=\E_R f_1$. Independence and weighted AM--GM give
\[
\E_{R^{a+b}} E^*_{a,b}=u^a z^b
 =(u^w z^v)^N\le(wu+vz)^N=1.
\]
The oracle's expected logarithm is exactly \eqref{eq:imb-information}. For \eqref{eq:imb-score-e}, use $f_0=(w+ve^h)^{-1}$ and $f_1=e^h/(w+ve^h)$; their weighted sum is one pointwise. Conditional on an independent pilot, the identical calculation applies.
\end{proof}

The score targets $h^*=\log(q/p)$, with densities taken relative to $P+Q$ and with extended-real values permitted. For a fixed bounded $h$, as $a/b\to\infty$,
\begin{equation}\label{eq:imb-variational-limit}
\frac{\E F_{a,b}(h)}b\longrightarrow
J(h):=\E_Qh-\E_Pe^h+1.
\end{equation}
Indeed, $\log\{w+ve^h\}=\log\{1+v(e^h-1)\}
=v(e^h-1)+O_h(v^2)$ uniformly in the observation. Consequently the $X$-sum, divided by $b$, has expectation tending to $\E_P(e^h-1)$, whereas the corresponding $Y$-term tends to zero.

\begin{lemma}[Variational limit]\label{lem:imb-variational}
For arbitrary $P,Q$,
\begin{equation}\label{eq:imb-J-sup}
\sup_{h\ \mathrm{bounded,\ measurable}}J(h)=\KL(Q\|P).
\end{equation}
Moreover, if $h_B^*=\clip(\log(q/p),-B,B)$, then $J(h_B^*)\to\KL(Q\|P)$ as $B\to\infty$, including the value $+\infty$. Therefore
\begin{equation}\label{eq:imb-info-limit}
\frac{I_{a,b}(P,Q)}b\longrightarrow\KL(Q\|P)
\quad\text{whenever }a/b\to\infty.
\end{equation}
\end{lemma}
\begin{proof}
If $Q\ll P$, write $r=dQ/dP$. For finite $\KL(Q\|P)$,
\[
\KL(Q\|P)-J(h)
 =\int\{r\log(r/e^h)-r+e^h\}\,dP\ge0.
\]
For $h_B^*$, the clipped density ratio converges to $r$, is bounded above by $r+1$, and hence its $P$-integral tends to one. The negative part of $\log r$ under $Q$ is integrable because $r|\log r|\one\{r<1\}\le1/e$ pointwise. Truncation therefore gives $\E_Qh_B^*\to\KL(Q\|P)$, also when the positive part is infinite. If $Q$ has singular mass $\beta>0$ relative to $P$, then $h_B^*=B$ on that singular part, its contribution to $\E_Qh_B^*$ is $\beta B$, and its $P$-integral is zero. The same argument on the absolutely continuous part shows divergence to infinity. This proves \eqref{eq:imb-J-sup}. By \cref{thm:imb-oracle}, $I_{a,b}\ge\E F_{a,b}(h)$ for every bounded $h$, while $I_{a,b}\le b\KL(Q\|P)$ when the latter is finite. Apply \eqref{eq:imb-variational-limit} and then take the supremum over $h$.
\end{proof}

\subsection{Almost-sure information upper bounds}
We evaluate a field indexed by two sample sizes along the actual sampling path. Optional-stopping validity is imposed on the resulting adapted process.

\begin{lemma}[Change of measure for e-processes]\label{lem:imb-change}
Let $W_t$ be an e-process under a law $\mathbb P_0$, and suppose an alternative $\mathbb P_1$ is locally absolutely continuous with likelihood-ratio process $L_t$. Then $W_t/L_t$, defined $\mathbb P_1$-almost surely, is an e-process under $\mathbb P_1$. In particular,
\[
\sup_{t\ge0} W_t/L_t<\infty\qquad\mathbb P_1\text{-a.s.}
\]
\end{lemma}
\begin{proof}
For every bounded stopping time $\tau$,
$\E_1[W_\tau/L_\tau]=\E_0[W_\tau\one\{L_\tau>0\}]\le1$.
Applying this inequality at the first crossing of a level $c$, truncated at a deterministic horizon, gives
$\mathbb P_1(\sup_t W_t/L_t\ge c)\le1/c$.
Let $c\to\infty$.
\end{proof}

\begin{corollary}[Pathwise information upper bounds]\label{cor:imb-upper}
Under independent i.i.d. streams and any predictable allocation rule, every e-process for the common-law null satisfies the following. If the sample counts $a_t,b_t$ obey $b_t\to\infty$ and $\KL(Q\|P)<\infty$, then
\begin{equation}\label{eq:imb-as-upper}
\limsup_t\frac{\log W_t}{b_t}\le\KL(Q\|P)\qquad\text{a.s.}
\end{equation}
If instead $a_t/(a_t+b_t)\to\rho\in(0,1)$ and $a_t+b_t\to\infty$, then
$\limsup_t\log W_t/(a_t+b_t)\le G_\rho(P,Q)$ almost surely.
\end{corollary}
\begin{proof}
For the first claim choose the common null $R=P$. The likelihood ratio is
$L_t=\prod_{j=1}^{b_t}(dQ/dP)(Y_j)$.
Its normalized logarithm converges almost surely to $\KL(Q\|P)$ by the strong law. \Cref{lem:imb-change} gives the result. For the second claim choose the fixed null $R=\rho P+(1-\rho)Q$ and apply the strong law separately to both streams. The log density ratios are integrable, since the densities relative to $R$ are bounded above and the negative log parts are integrable. These arguments also cover random subsequences of the stream counts.
\end{proof}

\subsection{Non-kernel construction: a universal code over the pooled null MLE}
Let $\mathcal A=\{C_1,\ldots,C_K\}$ be a fixed measurable partition. Write $A_c$ and $B_c$ for cell counts, $C_c=A_c+B_c$, and define the Laplace probability of a categorical sequence with count vector $u$ and length $n$ by
\begin{equation}\label{eq:imb-code}
q_K(u)=\frac{\Gamma(K)\prod_{c=1}^K\Gamma(u_c+1)}{\Gamma(n+K)}.
\end{equation}
Here $q_K(u)$ assigns a probability to each ordered sequence with count vector $u$. Set
\begin{equation}\label{eq:imb-code-e}
V^{\mathcal A}_{a,b}
 =\frac{q_K(A)q_K(B)}{\prod_{c=1}^K(C_c/N)^{C_c}},
\end{equation}
with $0^0=1$ and value one at $(0,0)$. This is a prequential alternative likelihood divided by the pooled null maximum likelihood. The general domination principle is part of universal inference \citep{wasserman2020universal}; universal-predictor/null-MLE processes are also studied by \citet{shekhar2025optimal}.

\begin{theorem}[Validity and finite-partition regret]\label{thm:imb-code}
Along every predictable allocation path, \eqref{eq:imb-code-e} is an e-process for the fixed common-law i.i.d. null. For $a,b\ge1$, write $\widehat p=A/a$, $\widehat q=B/b$, and $\widehat m=C/N$. Then
\begin{align}\label{eq:imb-code-regret}
\log V^{\mathcal A}_{a,b}\ge{}&
 a\KL(\widehat p\|\widehat m)+b\KL(\widehat q\|\widehat m)\\[-3pt]
&-\log\binom{a+K-1}{K-1}-\log\binom{b+K-1}{K-1}.\nonumber
\end{align}
Each binomial penalty is at most $(K-1)\log(n+1)$ for the corresponding sample size $n$.
\end{theorem}
\begin{proof}
For a true null cell vector $r$, the denominator in \eqref{eq:imb-code-e} is at least $\prod_c r_c^{C_c}$. The resulting upper bound
$q_K(A)q_K(B)/\prod_c r_c^{C_c}$ is the product of separate Laplace predictive probabilities divided by the true probability of each next observation. It is a nonnegative supermartingale under the null, including when some $r_c=0$, along every predictable allocation path. Domination gives the e-process property.
For a sequence with counts $u$, the maximum likelihood is $\prod_c(u_c/n)^{u_c}$. Its multinomial count probability is at most one, so
\[
\prod_c(u_c/n)^{u_c}\le \frac{\prod_c u_c!}{n!}
 =\binom{n+K-1}{K-1}q_K(u).
\]
Apply this to each stream and subtract the pooled maximum log likelihood. Finally, weak compositions of $n$ into $K$ cells inject into $\{0,\ldots,n\}^{K-1}$, giving the claimed penalty bound.
\end{proof}

On a standard Borel space choose nested finite partitions $\mathcal A_j$ generating the sigma-field, include the one-cell partition, and assign every partition a fixed positive weight $\pi_j$ with total one. Define
\begin{equation}\label{eq:imb-partition-mix}
V^{\mathrm{part}}_{a,b}=\sum_{j\ge0}\pi_jV^{\mathcal A_j}_{a,b}.
\end{equation}
The one-cell component is identically one and supplies a cash floor.

\begin{theorem}[Universal asymmetric attainment by partition coding]\label{thm:imb-partition-opt}
The field \eqref{eq:imb-partition-mix}, evaluated along any predictable allocation path, is an e-process. For every fixed $(P,Q)$:
\begin{align}
\frac{\log V^{\mathrm{part}}_{a,b}}N&\longrightarrow G_\rho(P,Q)
&&\text{if }a/N\to\rho\in(0,1),\ N\to\infty;\label{eq:imb-partition-interior}\\
\frac{\log V^{\mathrm{part}}_{a,b}}b&\longrightarrow\KL(Q\|P)
&&\text{if }a/b\to\infty,\ b\to\infty,\ \log a=o(b).
\label{eq:imb-partition-edge}
\end{align}
Convergence is almost sure; the second limit is in the extended sense if the divergence is infinite. The reverse-imbalance statement follows by interchanging the streams. Along deterministic size sequences and for finite limiting information, convergence also holds in $L^1$ and in expectation.
\end{theorem}
\begin{proof}
A fixed weighted sum of the null-dominating supermartingales proves validity. Fix a partition. For an interior limiting allocation, continuity of finite-alphabet entropy and the strong law give the limiting empirical weighted information, while the logarithmic code penalties are $o(N)$.
For extreme imbalance, first suppose $Q(C_c)>0$ only when $P(C_c)>0$. Eventually all positive $P$ cells have empirical probabilities bounded away from zero. Expanding $\log(1+v(\widehat q_c/\widehat p_c-1))$ on this finite support gives
\[
(a/b)\KL(\widehat p\|\widehat m)\to0,
\qquad \KL(\widehat q\|\widehat m)\to\KL(Q_{\mathcal A}\|P_{\mathcal A}).
\]
For example, the first convergence follows from a uniform $O(v^2)$ bound on the divergence, multiplied by $a/b=O(1/v)$. If a cell has $P(C_c)=0<Q(C_c)$, then $A_c=0$ identically and its contribution to $\KL(\widehat q\|\widehat m)$ is $\widehat q_c\log(1/v)\to\infty$; all other negative contributions together are bounded below by $-1/e$. Thus the same conclusion holds with value infinity. The code penalties divided by $b$ vanish under $\log a=o(b)$.
Since $\log V^{\mathrm{part}}\ge\log\pi_j+\log V^{\mathcal A_j}$, the limit inferior is at least every finite-partition information. These partition informations increase to the full information. To verify this without a topological assumption, approximate any bounded score $h$ by $h_j=\E_{(P+Q)/2}[h\mid\sigma(\mathcal A_j)]$. The bounded martingale convergence theorem gives convergence in $L^1(P+Q)$; consequently $J(h_j)\to J(h)$ and the analogous weighted logistic objectives converge. Taking suprema and applying \cref{lem:imb-variational}, or the truncated oracle for the interior case, proves the assertion. \Cref{cor:imb-upper} gives the matching upper bounds.
For deterministic sequences, the cash floor bounds normalized log wealth below by a deterministic quantity tending to zero. Fatou's lemma and \eqref{eq:imb-log-upper} sandwich its expectations at the almost-sure limit. Convergence in $L^1$ follows by applying the standard nonnegative convergence lemma to normalized log wealth shifted by its vanishing deterministic lower bound.
\end{proof}

\paragraph{Finite implementation.}
A deterministic increasing finite subset of the partition sum, containing the one-cell component and eventually every fixed partition, remains valid: it is pointwise bounded by the full mixture. Its logarithmic rate is unchanged. The reported statistic is the unnormalized partial sum with the original prior weights. Each active partition needs two count vectors. A new observation in cell $c$ of a stream of current size $n$ adds $\log(u_c+1)-\log(n+K)$ to that stream's log code. Maintaining $\sum_c C_c\log C_c-N\log N$ updates the pooled maximum log likelihood in constant time per active partition. Newly activated partitions can be evaluated retrospectively from stored data. The partition sequence and its prior weights are fixed before observing the data.

\subsection{A kernel learner whose training cost is negligible for the smaller sample}
The kernel method trains on a vanishing fraction of balanced pairs and uses the remaining observations in an unequal-block bet. Its training cost is therefore determined by the smaller sample.

Let $k$ be bounded, measurable and characteristic, with separable RKHS and $k(x,x)\le1$, and use $K=(1+k)/2$. These are precisely the general-space conditions from the balanced construction. The role of characteristicness versus function approximation is discussed by \citet{sriperumbudur2011}. By \cref{lem:characteristic-density}, the RKHS of $K$ is dense in $L^2(\mu)$ for every probability law $\mu$.

We reverse the label convention for this section: $Y$ has label one and $X$ has label zero. The balanced training loss is
\begin{equation}\label{eq:imb-train-loss}
\ell_i(h)=\log(1+e^{h(X_i)})+\log(1+e^{-h(Y_i)}),
\qquad \mathcal R(h)=\E\ell_i(h).
\end{equation}
Its Bayes probability is $\eta^*=q/(p+q)$ and its Bayes score is $h^*=\log(q/p)$.

For a training size $s$, set $R_s=s^{1/4}$ and $\gamma_s=R_s/(2\sqrt s)$. Starting at $h_{s,1}=0$, run
\begin{align}
g_{s,i}&=\sig(h_{s,i}(X_i))K(X_i,\cdot)
       +\{\sig(h_{s,i}(Y_i))-1\}K(Y_i,\cdot),\nonumber\\
h_{s,i+1}&=\operatorname{Proj}_{\{\|h\|_{\HH_K}\le R_s\}}
               (h_{s,i}-\gamma_s g_{s,i}),\qquad 1\le i\le s,\label{eq:imb-ogd}\\
\overline h_s&=s^{-1}\sum_{i=1}^s h_{s,i}.\nonumber
\end{align}
The horizon-dependent training run is computed separately when $s$ changes. Averaging is over the pre-update scores.

\begin{lemma}[Strong universal calibration of the pilot learner]\label{lem:imb-training}
For every fixed $P,Q$, the learner \eqref{eq:imb-ogd} satisfies
\[
\mathcal R(\overline h_s)\to\mathcal R^*,\qquad
\int|\sig(\overline h_s)-\eta^*|\,d\mu\to0
\quad\text{a.s.},\qquad \mu=(P+Q)/2.
\]
For every fixed $B<\infty$,
\begin{equation}\label{eq:imb-clip-consistency}
\int|\clip(\overline h_s,-B,B)-h_B^*|\,d(P+Q)\to0
\qquad\text{a.s.}
\end{equation}
\end{lemma}
\begin{proof}
The gradient norm is at most two. Projection and convexity give, for every $\|h\|\le R_s$,
\[
\sum_{i=1}^s\ell_i(h_{s,i})-\sum_{i=1}^s\ell_i(h)
\le\frac{R_s^2}{2\gamma_s}+2\gamma_s s=2R_s\sqrt s.
\]
The loss lies in $[0,L_s]$, where $L_s=2\log(1+e^{R_s})$. Within each training run, $\ell_i(h_{s,i})-\mathcal R(h_{s,i})$ is a martingale difference. Conditional Hoeffding bounds show that, with probability one,
\[
\left|s^{-1}\sum_{i=1}^s
 [\ell_i(h_{s,i})-\mathcal R(h_{s,i})]\right|
 =O\left(L_s\sqrt{\frac{\log(s+1)}s}\right)=o(1).
\]
Choosing the Hoeffding threshold so that its failure probability is at most $2(s+1)^{-3}$ allows Borel--Cantelli to be applied across training horizons. Convexity gives $\mathcal R(\overline h_s)\le s^{-1}\sum_i\mathcal R(h_{s,i})$. For any fixed RKHS comparator, the regret bound, the preceding martingale bound, and the ordinary strong law imply $\limsup_s\mathcal R(\overline h_s)\le\mathcal R(h)$. A countable sequence of comparators approximates the Bayes risk. Indeed, clipping an optimal Bernoulli probability to $[\sig(-C),\sig(C)]$ incurs Bernoulli KL at most $\log(1+e^{-C})$: on each clipped interval convexity places the maximum at its endpoint. Thus $\mathcal R(h_C^*)-\mathcal R^*\le2\log(1+e^{-C})$. Truncated Bayes scores belong to $L^2(\mu)$, logistic loss is Lipschitz, and the RKHS is dense, proving the approximation assertion. This proves risk convergence.
The calibration identity is
\[
\mathcal R(h)-\mathcal R^*
 =2\int\kl(\eta^*\|\sig(h))\,d\mu
 \ge4\int|\eta^*-\sig(h)|^2\,d\mu.
\]
Finally, the map $u\mapsto\clip(\log(u/(1-u)),-B,B)$ is a bounded Lipschitz function on $[0,1]$, with the endpoint values defined by continuity. Applying it proves \eqref{eq:imb-clip-consistency}.
\end{proof}

\subsection{Unequal-block betting and protection over sample sizes}
For $a,b\ge2$, let
\[
s=\lfloor\sqrt m\rfloor,\qquad A=a-s,\quad B_0=b-s,\quad
w'=A/(A+B_0),\quad v'=B_0/(A+B_0).
\]
Only the first $s$ observations of each stream train $\overline h_s$. For each clipping level $L\in\N$, evaluate \eqref{eq:imb-score-F} on the remaining observations, with score
$h_{s,L}=\clip(\overline h_s,-L,L)$ and weights $w',v'$. Denote this log statistic by $F^{\rm test}_{a,b,L}$ and its exponential by $E^{\rm test}_{a,b,L}$. Let
\begin{equation}\label{eq:imb-kernel-base}
H_{a,b}=\frac12+\frac12\sum_{L\ge1}\frac{E^{\rm test}_{a,b,L}}{L(L+1)}.
\end{equation}
For $a<2$ or $b<2$, set $H_{a,b}=1$. Each $H_{a,b}$ is a fixed-size e-variable. We combine the endpoints through the following sample-size correction.

Put
\begin{equation}\label{eq:imb-size-weights}
\kappa_{a,b}=\frac{1}{a(a+1)b(b+1)},\qquad a,b\ge1,
\qquad W^{\rm ker}_{a,b}=\kappa_{a,b}H_{a,b}.
\end{equation}
The double sum of $\kappa_{a,b}$ equals one. The correction is independent of the significance level.

\begin{lemma}[Simultaneous sample-size protection]\label{lem:imb-stitch}
If $H_{a,b}$ is a fixed-size e-variable for each $a,b\ge1$, then, under every common-law null,
\[
\E\!\left[\kappa_{A,B}H_{A,B}\right]\le1
\]
for every random finite pair $(A,B)$, even if selected after inspecting additional observations. In particular, the corrected field is an e-process along any sampling path. It also obeys
\[
\PP_0\left(\sup_{a,b\ge1}\kappa_{a,b}H_{a,b}\ge1/\alpha\right)\le\alpha.
\]
At calendar time zero its value may be set to one; until both samples are present it may be set to zero.
\end{lemma}
\begin{proof}
Pointwise, $\kappa_{A,B}H_{A,B}\le Z:=\sum_{a,b}\kappa_{a,b}H_{a,b}$, and Tonelli gives $\E_0Z\le1$. The stopping and crossing claims follow. The initialization at one applies to a trivial initial sigma-field; the same conclusion holds conditionally on independent initial randomization.
\end{proof}

\begin{lemma}[Testing fluctuations on the minority-sample scale]\label{lem:imb-concentration}
For every fixed clipping level $L$, almost surely,
\begin{equation}\label{eq:imb-testing-fluctuation}
\frac{F^{\rm test}_{a,b,L}
 -\E[F^{\rm test}_{a,b,L}\mid X_{1:s},Y_{1:s}]}{m}\to0
\end{equation}
along every sequence with $m\to\infty$ and $\log\max(a,b)=o(m)$. The probability-one event can be chosen simultaneously for all such sequences and all integer $L$.
\end{lemma}
\begin{proof}
Condition on the pilot. The test observations are independent. On $[-L,L]$, the $X$ log factor $-\log(w'+v'e^h)$ has range at most $c_Lv'$, and the $Y$ log factor $h-\log(w'+v'e^h)$ has range at most $c_Lw'$, where $c_L=2Le^{2L}$. These follow by differentiating: the derivative magnitudes are respectively at most $v'e^{2L}$ and $w'e^{2L}$. Consequently the sum of squared ranges is at most
\[
c_L^2\{A(v')^2+B_0(w')^2\}
 =c_L^2\frac{AB_0}{A+B_0}\le c_L^2m.
\]
Hoeffding's inequality therefore bounds the probability that the absolute numerator in \eqref{eq:imb-testing-fluctuation} exceeds $\varepsilon m$ by $2\exp(-2\varepsilon^2m/c_L^2)$. This bound also holds unconditionally. For a fixed $m$, there are at most $2e^{\delta m}$ pairs with that minimum and maximum at most $e^{\delta m}$. Choose $0<\delta<2\varepsilon^2/c_L^2$ and sum the resulting probabilities over $m$. Borel--Cantelli applies. Every sequence in the statement eventually lies in this set of pairs. Intersect over positive rational $\varepsilon$ and integer $L$.
\end{proof}

\begin{theorem}[Universal asymmetric attainment by kernels]\label{thm:imb-kernel-opt}
Under the stated measurable characteristic-kernel and separability assumptions, \eqref{eq:imb-size-weights} is an e-process for equality of two independent i.i.d. laws. For every fixed $(P,Q)$ it satisfies, almost surely,
\begin{align}
\frac{\log W^{\rm ker}_{a,b}}N&\to G_\rho(P,Q)
&&\text{if }a/N\to\rho\in(0,1),\ N\to\infty;\label{eq:imb-kernel-interior}\\
\frac{\log W^{\rm ker}_{a,b}}b&\to\KL(Q\|P)
&&\text{if }a/b\to\infty,\ b\to\infty,\ \log a=o(b).
\label{eq:imb-kernel-edge}
\end{align}
The second convergence includes infinite divergence, and the reverse statement holds with the streams exchanged. Along deterministic size sequences, finite limits also hold in $L^1$ and in expectation.
\end{theorem}
\begin{proof}
Conditional validity follows from \cref{thm:imb-oracle}, mixtures preserve it, and \cref{lem:imb-stitch} proves optional-stopping validity. Fix $L$. In the extreme-imbalance regime, the uniform Taylor expansion used in \eqref{eq:imb-variational-limit}, now for the random but uniformly bounded score $h_{s,L}$, gives
\[
\frac{\E[F^{\rm test}_{a,b,L}\mid\text{pilot}]}b
 -J(h_{s,L})\to0.
\]
Here $s=o(b)$ and the testing block has the same limiting imbalance. By \cref{lem:imb-training}, $J(h_{s,L})\to J(h_L^*)$: both $h$ and $e^h$ are Lipschitz on the clipping interval. \Cref{lem:imb-concentration} proves that the realized normalized log statistic has the same limit. Further,
\[
\log W^{\rm ker}_{a,b}
 \ge F^{\rm test}_{a,b,L}+\log\kappa_{a,b}-\log\{2L(L+1)\},
\]
and $-\log\kappa_{a,b}=2\log a+2\log b+O(1)=o(b)$. Let $L\to\infty$ and apply \cref{lem:imb-variational} to prove the lower bound. \Cref{cor:imb-upper} proves the finite-divergence upper bound.
For an interior allocation, the same steps yield the weighted population log-growth at score $h_L^*$. As $L\to\infty$ it converges to $G_\rho$. To justify the last passage, each individual log factor at the clipped oracle lies between its value zero at score zero and its value at the full oracle, because the $X$ factor decreases and the $Y$ factor increases with the score. The full-oracle log factors are integrable. Dominated convergence applies. The correction is $o(N)$ and the same information upper bound completes the proof.
Finally, $W^{\rm ker}_{a,b}\ge\kappa_{a,b}/2$. Its normalized log lower bound tends to zero under the displayed size conditions. Fatou's lemma, the deterministic-size information upper bound, and the nonnegative convergence argument used for partitions prove the expectation and $L^1$ assertions.
\end{proof}

\paragraph{Clipping and calibration.}
The imbalanced objective involves both $\E_Qh$ and $\E_Pe^h$. Fixed clipping makes these functionals continuous under the probability-calibration convergence supplied by the training lemma. The mixture over fixed clipping levels then recovers the full variational optimum, using consistency at each level.

\paragraph{Finite kernel implementation.}
At minimum sample size $m$, one may retain only $L=1,\ldots,J(m)$, where $J(m)\uparrow\infty$ deterministically, and replace the omitted prior mass $1/(J(m)+1)$ by cash in the fixed-size statistic $H_{a,b}$. This is valid before the size correction and leaves every fixed-level lower bound intact. A finite mixture over kernels is also valid; a countable mixture with finite active subsets remains universal if every fixed component is eventually retained and at least one component satisfies the approximation conditions. Pilot-only kernel tuning preserves validity. Retaining a fixed universal component also preserves the stated consistency guarantee.

The averaged score is a linear combination of at most $2s$ kernel sections. Its projected-gradient training can be performed in $O(s^2)$ kernel evaluations and arithmetic, with $O(s^2)$ Gram-matrix storage. In particular, a proposed update $h-\gamma g$ has squared norm
\[
\|h\|^2-2\gamma\langle h,g\rangle+\gamma^2\|g\|^2;
\]
$g$ contains two sections, allowing the squared norm to be maintained from their kernel evaluations and current scores. Testing requires $O((a+b)s)$ kernel evaluations, after which all clipping-level scores reuse the same predictions. With $s\asymp\sqrt m$, the training cost depends on the smaller sample. The procedure remains valid with any pilot-trained predictor; the particular averaged gradient algorithm supplies the universal-rate proof. Evaluating at selected checkpoints saves repeated testing work. The same sample-size correction applies when reporting statistics at a subset of checkpoints.

\subsection{Removing the condition on the logarithm of the larger sample}
The preceding all-data results require $\log\max(a,b)=o(m)$. Retaining a slowly diverging ratio of majority to minority observations removes this condition while preserving first-order information.

Let $L(m)\ge1$ be deterministic with $L(m)\to\infty$ and $\log L(m)=o(m)$; for example $L(m)=\lceil\log(m+1)\rceil$. Define retained sizes
\begin{equation}\label{eq:imb-capped-counts}
u=\min\{a,\lfloor mL(m)\rfloor\},\qquad
z=\min\{b,\lfloor mL(m)\rfloor\}.
\end{equation}
We use the first $u$ and $z$ observations of the respective streams.

\begin{corollary}[Optimality under arbitrary diverging imbalance]\label{cor:imb-cap}
For the kernel method, report
$\widetilde W^{\rm ker}_{a,b}=\kappa_{u,z}H_{u,z}$.
For the partition method, report
$\widetilde W^{\rm part}_{a,b}=\kappa_{u,z}V^{\rm part}_{u,z}$.
Both are anytime-valid, even in the filtration containing all observations, including those not retained for the displayed statistic. If $a/b\to\infty$ and $b\to\infty$, then, without any condition on $\log a$,
\begin{equation}\label{eq:imb-cap-limit}
\frac{\log\widetilde W^{\rm ker}_{a,b}}b\to\KL(Q\|P),
\qquad
\frac{\log\widetilde W^{\rm part}_{a,b}}b\to\KL(Q\|P)
\qquad\text{a.s.}
\end{equation}
The reverse-imbalance and interior-allocation conclusions also hold.
\end{corollary}
\begin{proof}
Both uncorrected fields are fixed-size e-variables, so \cref{lem:imb-stitch} proves validity even when the retained endpoint is selected after seeing additional data. Under $a/b\to\infty$, eventually $m=b$ and $z=b$, while
\[
\frac ub=\min\{a/b,\lfloor bL(b)\rfloor/b\}\to\infty,
\qquad \log u\le\log b+\log L(b)=o(b).
\]
Thus the lower-bound proofs of \cref{thm:imb-partition-opt,thm:imb-kernel-opt} apply to retained prefixes. Those proofs hold on a single probability-one event for all size sequences satisfying their conditions: partitions use ordinary streamwise strong laws, and the kernel testing bound is simultaneous by \cref{lem:imb-concentration}. The additional partition size penalty is $o(b)$. The finite-divergence upper bound follows from \cref{cor:imb-upper} applied in the actual full filtration, since the original minority sample count is still $b$. Infinite divergence needs no upper bound. For interior limiting proportions, neither sample is eventually capped. Interchanging the streams proves the remaining case.
\end{proof}

\subsection{One procedure for every limiting sample-size ratio}
Fix the cap $L(m)=\lceil\log(m+1)\rceil\vee1$, the partitions or kernel, and all prior weights before observing data. For either construction in \cref{cor:imb-cap}, let $\mathcal W_{a,b}=\widetilde W^{\rm part}_{a,b}$ or $\mathcal W_{a,b}=\widetilde W^{\rm ker}_{a,b}$, respectively. This is one rule for the entire sample-size field; its definition contains no limiting allocation parameter.

\begin{corollary}[Simultaneous optimality across allocation regimes]\label{cor:all-ratios}
For every fixed pair $(P,Q)$, under the corresponding observation-space assumptions, the process $\mathcal W$ is anytime-valid for the common-law null. Along any predictable allocation with $a,b\to\infty$ and $a/(a+b)\to\rho$, for a fixed $\rho\in[0,1]$, it satisfies
\begin{equation}\label{eq:all-ratios}
\frac{\log\mathcal W_{a,b}}{\min(a,b)}\longrightarrow
\begin{cases}
\KL(P\|Q),&\rho=0,\\[2pt]
G_\rho(P,Q)/\min\{\rho,1-\rho\},&0<\rho<1,\\[2pt]
\KL(Q\|P),&\rho=1,
\end{cases}
\qquad\text{a.s.}
\end{equation}
The endpoint limits hold in the extended-real sense. The same specified procedure attains these rates for every $\rho$ and every pair of laws in its model. Whenever $P\ne Q$ and the displayed limit is finite, it is strictly positive and
\[
\frac{\log\mathcal W_{a,b}}{I_{a,b}(P,Q)}\longrightarrow1
\qquad\text{a.s.}
\]
\end{corollary}
\begin{proof}
If $0<\rho<1$, the sample-size ratio is eventually bounded above and below, so the cap is eventually inactive. The interior conclusions of \cref{thm:imb-partition-opt,thm:imb-kernel-opt}, together with $\min(a,b)/(a+b)\to\min\{\rho,1-\rho\}$, give the middle case; the additional partition sample-size correction is $o(\min(a,b))$. At $\rho=1$, the majority sample is $a$ and \cref{cor:imb-cap} gives the third case. Interchanging the streams gives the first case. Validity is already simultaneous over all reported endpoints by \cref{lem:imb-stitch}. The limiting normalized oracle information is the same in the interior by its definition and continuity, and at the endpoints by \cref{lem:imb-variational}. Division proves the final assertion when this limit is finite and positive.
\end{proof}

\paragraph{Combining a practical learner with the universal construction.}
Let $V$ be any e-process in the same sampling filtration, and let $0<\varepsilon<1$ be fixed. The mixture
\[
W=\varepsilon\mathcal W+(1-\varepsilon)V
\]
is an e-process and retains every limit in \cref{cor:all-ratios}. Indeed, $\log W\ge\log\mathcal W+\log\varepsilon$, and the information upper bounds apply to the mixture. Thus a universal component can be combined with a specialized learner at a constant log-wealth cost. In particular, the allocation-adaptive coding mixture
\begin{equation}\label{eq:adaptive-code-mixture}
\widehat W^{\rm part}_{a,b}
 =\tfrac12 V^{\rm part}_{a,b}+\tfrac12\widetilde W^{\rm part}_{a,b}
\end{equation}
retains all observations through its first component and the unrestricted-imbalance guarantee through its second. Its log wealth is at least that of the all-data code minus $\log2$, while it satisfies \cref{cor:all-ratios} at every limiting ratio. Finite model banks used on their own retain validity and their comparator guarantees; expanding the bank supplies general-space universal attainment.

\begin{remark}[Sample allocation and first-order information]
Capping preserves the first-order information on the minority-sample scale while changing the finite-sample oracle information. The capped statistics and corrected kernel field are e-processes in the full filtration. Another valid construction multiplies bets over disjoint blocks whose sizes are fixed prospectively; its oracle information is the sum of the block informations. In particular, one-pair blocks yield the balanced rate $2\JS(P,Q)$ per pair. Retaining a diverging reference-to-test ratio allows the rate per minority observation to approach $\KL(Q\|P)$.
\end{remark}

\paragraph{Implementation choices.}
Partition coding gives an online all-data e-process with logarithmic redundancy. The kernel construction trains a calibrated score on $\lfloor\sqrt m\rfloor$ pairs, freezes and clips it, evaluates the weighted product on the remaining data, and applies a logarithmic sample-size correction. Both achieve the exact weighted information for fixed limiting proportions and the one-sample KL rate under strong imbalance. Capping the overwhelming sample at $m\lceil\log(m+1)\rceil$ removes the logarithmic-size restriction without changing either first-order limit.

\section{Fixed-sample tests and classical efficiency}\label{sec:batch-efficiency}
The information rate determines the significance attained under a fixed alternative. Error probabilities and local alternatives require additional analysis. This section gives the general-space Bahadur consequence and a complete finite-alphabet analysis of the coding construction. The latter identifies a single statistic with sequential, large-deviation, and classical conditional calibrations. All sample sizes in this section are deterministic; sequential validity of the underlying e-process is retained.

\subsection{Optimal attained-significance rates}
Write $N=a+b$ and $\rho_N=a/N$. For an e-process evaluated after these observations, let
\[
 p^e_{a,b}=\min\{1,W_{a,b}^{-1}\}.
\]
Its running-maximum version along a prescribed sampling path is also a valid p-value. We use the classical factor-two convention for the Bahadur slope of the reported p-value \citep{bahadur1960}: the slope at normalization $s$ is the limit of $-2\log p/s$.

\begin{theorem}[Optimal Bahadur slopes]\label{thm:ts-bahadur}
Suppose $P\ne Q$, and use either universal construction in \cref{cor:all-ratios}, or a fixed positive-weight mixture containing it. If $\rho_N\to\rho\in(0,1)$, then
\begin{equation}\label{eq:ts-bahadur-interior}
 -N^{-1}\log p^e_{a,b}\longrightarrow
 G_\rho(P,Q)=\rho\KL(P\|M_\rho)+(1-\rho)\KL(Q\|M_\rho)
 \qquad\text{a.s.}
\end{equation}
The same limit holds for the running-maximum p-value. For every sequence $p_{a,b}$ of p-values valid under all common laws,
\[
 \limsup_N -N^{-1}\log p_{a,b}\le G_\rho(P,Q)\qquad\text{a.s.}
\]
Thus the optimal slope per individual observation is $2G_\rho(P,Q)$. In the balanced experiment it is $2\JS(P,Q)$ per observation, or $4\JS(P,Q)$ per pair.

If $a_m/m\to\infty$ and the other sample has size $b_m=m\to\infty$, then
\begin{equation}\label{eq:ts-bahadur-edge}
 -m^{-1}\log p^e_{a_m,m}\longrightarrow\KL(Q\|P)\qquad\text{a.s.},
\end{equation}
including infinite divergence. For finite divergence this is an upper bound on the limit superior for every valid sequence of p-values at these endpoints. The reverse allocation exchanges $P$ and $Q$. The same endpoint limits hold for the running maximum of the universal e-process.
\end{theorem}
\begin{proof}
For the interior converse use the common null $R=M_\rho$ and its likelihood ratio $L_{a,b}$ against $P^a\otimes Q^b$. Streamwise strong laws give $N^{-1}\log L_{a,b}\to G_\rho$. Superuniformity gives, for every $\varepsilon>0$,
\[
 \PP_{P,Q}\{p_{a,b}\le e^{-N(G_\rho+\varepsilon)},\
 L_{a,b}\le e^{N(G_\rho+\varepsilon/2)}\}
 \le e^{-N\varepsilon/2}.
\]
Summation over integer $N$ and Borel--Cantelli prove the upper bound. For the boundary converse take $R=P$: the likelihood ratio is $\prod_{j=1}^m(dQ/dP)(Y_j)$, so the same proof sums $e^{-m\varepsilon/2}$. The endpoint indexing places one comparison at each minority sample size; it avoids making a multiple-look claim for unrelated fixed-endpoint p-values. Attainment follows from \cref{cor:all-ratios}. For running maxima, eventual upper and lower growth bounds and monotonicity of the sample counts give the same limits. At the boundary one can alternatively apply \cref{lem:imb-change} to bound wealth at every earlier time by a finite random multiple of the likelihood ratio involving the observations from the smaller stream.
\end{proof}
In particular, for a deterministic level $e^{-Nc}$ in the interior, terminal power tends to one if $c<G_\rho(P,Q)$ and to zero if $c>G_\rho(P,Q)$. This statement does not quantify the rate of the miss probability.

\subsection{An exact likelihood-ratio comparison on a finite alphabet}
For the remainder of this subsection and the next two, $P=(p_i)$ and $Q=(q_i)$ are distinct and strictly positive on a fixed alphabet of size $K\ge2$, and $\rho_N\to\rho\in(0,1)$. Let $V_{a,b}$ be the identity-partition code in \eqref{eq:imb-code-e}, and write
\[
 T_{a,b}=a\KL(\widehat P_a\|\widehat M_N)
       +b\KL(\widehat Q_b\|\widehat M_N)
       =NG_{\rho_N}(\widehat P_a,\widehat Q_b).
\]
This is the log generalized likelihood ratio for homogeneity. The coding bounds imply the deterministic sandwich
\begin{equation}\label{eq:ts-code-sandwich}
 T_{a,b}-d_{a,b}\le\log V_{a,b}\le T_{a,b},\qquad
 d_{a,b}=\log\binom{a+K-1}{K-1}+\log\binom{b+K-1}{K-1}.
\end{equation}
The lower bound is \cref{thm:imb-code}. The upper bound follows because each Laplace predictive probability is a mixture of iid likelihoods and is at most their maximum. Thus $d_{a,b}=O_K(\log N)$ uniformly over all data, including empty cells.

For candidate distributions $S,T$ on the same alphabet define
\begin{align}
 \mathcal D_\rho(S,T;P,Q)&=\rho\KL(S\|P)+(1-\rho)\KL(T\|Q),\nonumber\\
 B_\rho(c;P,Q)&=\inf_{G_\rho(S,T)\le c}\mathcal D_\rho(S,T;P,Q),\label{eq:ts-B}\\
 J_\rho(P,Q)&=B_\rho(0;P,Q)
 =\inf_R\{\rho\KL(R\|P)+(1-\rho)\KL(R\|Q)\}\nonumber\\
 &=-\log\sum_{i=1}^Kp_i^\rho q_i^{1-\rho}.\label{eq:ts-J}
\end{align}
For a test $\varphi_{a,b}$, let
$\alpha_{a,b}=\sup_R\E_{R^N}\varphi_{a,b}$ and
$\beta_{a,b}(P,Q)=\E_{P^a\otimes Q^b}(1-\varphi_{a,b})$.

\begin{theorem}[Optimal Stein and Hoeffding exponents]\label{thm:ts-errors}
At any fixed $0<\alpha<1$, the terminal test $\one\{V_{a,b}\ge1/\alpha\}$ satisfies
\begin{equation}\label{eq:ts-stein}
 \lim_N -N^{-1}\log\beta_{a,b}(P,Q)=J_\rho(P,Q).
\end{equation}
No sequence of tests with $\alpha_{a,b}\le\alpha$ has a larger limit superior exponent. For every $0<c<G_\rho(P,Q)$, the terminal test $\one\{V_{a,b}\ge e^{Nc}\}$ satisfies
\begin{equation}\label{eq:ts-hoeffding}
 \alpha_{a,b}\le e^{-Nc},\qquad
 \lim_N -N^{-1}\log\beta_{a,b}(P,Q)=B_\rho(c;P,Q).
\end{equation}
Every sequence of tests satisfying the displayed null constraint has miss-probability limit superior exponent at most $B_\rho(c;P,Q)$.
\end{theorem}
The two KL orientations describe different efficiencies. For $a=b=n$, the optimal miss exponent per pair is
$-2\log\sum_i\sqrt{p_iq_i}$, whereas the attained-significance exponent is $2\JS(P,Q)$. The proofs in \cref{app:batch-proofs} use two independent multinomial type bounds and \eqref{eq:ts-code-sandwich}. These are the empirical-divergence principles underlying classical optimal multinomial tests \citep{hoeffding1965,gutman1989,unnikrishnan2012}.

An almost-sure growth theorem alone does not imply \eqref{eq:ts-stein}. Indeed, independently erase the reported wealth with probability $N^{-2}$ at each endpoint. The resulting process is bounded by the original e-process, retains its almost-sure growth rate by Borel--Cantelli, and has miss probability at least $N^{-2}$. The deterministic coding sandwich rules out this behavior. The present exponent theorem concerns the finite-alphabet code; a kernel or neural learner requires its own large-deviation control or a positive-weight coding component.

\subsection{Chernoff calibration}
For a common null distribution $R$ define the two-group Chernoff information per observation by
\begin{equation}\label{eq:ts-C-R}
 C_\rho(R;P,Q)=-\min_{0\le u\le1}
 \left[\rho\log\sum_i r_i^{1-u}p_i^u
 +(1-\rho)\log\sum_i r_i^{1-u}q_i^u\right].
\end{equation}
Endpoint values are one-sided limits on the common supports. The common $u$ reflects a single likelihood-ratio decision applied to both samples.

\begin{theorem}[Optimal symmetric-error exponent]\label{thm:ts-chernoff}
The smallest achievable sum of worst-case type-I error and miss probability at the fixed alternative has exponent
\begin{align}
 C_\rho^{\rm hom}(P,Q)
 &=\inf_R C_\rho(R;P,Q)\nonumber\\
 &=\inf_{S,T}\max\{G_\rho(S,T),\mathcal D_\rho(S,T;P,Q)\}.
 \label{eq:ts-C}
\end{align}
There is a unique $c_*\in(0,G_\rho(P,Q))$ with
$B_\rho(c_*;P,Q)=c_*=C_\rho^{\rm hom}(P,Q)$.
The terminal code test with threshold $e^{Nc_*}$ satisfies
\[
 \lim_N-N^{-1}\log\{\alpha_{a,b}+\beta_{a,b}(P,Q)\}
 =C_\rho^{\rm hom}(P,Q).
\]
No sequence of tests has a larger limit superior exponent for this error sum.
\end{theorem}
This is an alternative-dependent threshold choice in the classical Chernoff criterion \citep{chernoff1952}. A fixed significance threshold serves the level-constrained criterion in \cref{thm:ts-errors}.

\subsection{Pitman alternatives and the exact permutation test}
Suppose that the data-generating probabilities vary with $N$ as
\begin{equation}\label{eq:ts-local-model}
 P_N=r+u/\sqrt N+o(N^{-1/2}),\qquad
 Q_N=r+v/\sqrt N+o(N^{-1/2}),
\end{equation}
where $r_i>0$ and $\sum_i u_i=\sum_i v_i=0$. Put
\[
 \Lambda=\rho(1-\rho)\sum_i\frac{(u_i-v_i)^2}{r_i}.
\]
For a fixed endpoint, permute the source labels of the $N$ pooled observations, preserving the two sample sizes, and recompute $V_{a,b}$. Denote the exact upper-tail permutation p-value by $p^{\rm perm}_{a,b}$.

\begin{theorem}[Local efficiency and conditional calibration]\label{thm:ts-pitman}
The permutation p-value is the probability-ordered Fisher--Freeman--Halton p-value for the $2\times K$ table of source by outcome. It is finite-sample superuniform under every common law. Uniformly over all tables,
\begin{equation}\label{eq:ts-perm-sandwich}
 -\log p^{\rm perm}_{a,b}=T_{a,b}+O_K(\log(N+1)).
\end{equation}
It has the optimal Bahadur slope and satisfies \cref{thm:ts-errors,thm:ts-chernoff} when rejection is defined by
$p^{\rm perm}_{a,b}\le\alpha$, $e^{-Nc}$, or $e^{-Nc_*}$, respectively.
Under \eqref{eq:ts-local-model},
\begin{align}
 2T_{a,b}&\Rightarrow\chi^2_{K-1}(\Lambda),\label{eq:ts-local-lrt}\\
 \log V_{a,b}
 &=T_{a,b}-(K-1)\log N+C_K(r,\rho)+o_P(1),\label{eq:ts-local-code}\\
 C_K(r,\rho)
 &=2\log\Gamma(K)+(K-1)\log(2\pi)+\sum_i\log r_i
 -\tfrac{K-1}{2}\log\{\rho(1-\rho)\}.\nonumber
\end{align}
Consequently the direct terminal e-value test at fixed level has power tending to zero. The exact permutation test has limiting power
\begin{equation}\label{eq:ts-pitman-power}
 \Pr\{\chi^2_{K-1}(\Lambda)>\chi^2_{K-1,1-\alpha}\},
\end{equation}
the same as the classical homogeneity likelihood-ratio test.
\end{theorem}
At fixed pooled counts $C_i$, the null probability of the allocation $A$ is
\[
 \mathfrak h(A\mid C)=\frac{a!b!}{N!}\prod_i\frac{C_i!}{A_i!B_i!}.
\]
The code is proportional to its reciprocal. This gives the exact connection to classical conditional testing \citep{freeman1951}; the local power agrees with the likelihood-ratio limit \citep{wilks1938}. Pitman efficiency here is one relative to that omnibus test on a fixed finite alphabet. The exact conditional test is already a classical procedure. Sequential validity concerns the original code process; repeated monitoring of its separately permutation-calibrated p-values requires an additional correction. Exact enumeration is feasible for small tables. A fixed number of Monte Carlo permutations has a p-value floor, so the stated Bahadur slopes concern exact probabilities or an approximation accurate at the required exponential scale.

\subsection{Minority-scale error exponents under arbitrary imbalance}
The finite-alphabet exponent analysis also extends to the boundary allocations. Suppose $a/b\to\infty$, $b\to\infty$, and $P,Q$ remain strictly positive on the fixed alphabet. Use the capped identity-partition code
\[
 U_{a,b}=\kappa_{u,b}V_{u,b},\qquad
 u=\min\{a,bL(b)\},\quad L(b)=\lceil\log(b+1)\rceil\vee1.
\]
Rounding $bL(b)$ is unnecessary with this integer choice. The rule is the specialization of \cref{cor:imb-cap}; it is valid in the full sampling filtration and imposes no restriction on the original majority sample size.

\begin{theorem}[Boundary error exponents]\label{thm:ts-boundary-errors}
For the preceding capped code, the fixed-level miss exponent is
\begin{equation}\label{eq:ts-boundary-stein}
 \lim -b^{-1}\log\beta_{a,b}(P,Q)=\KL(P\|Q).
\end{equation}
For $0<c<\KL(Q\|P)$, a threshold $e^{bc}$ has null rejection probability at most $e^{-bc}$ and miss exponent
\begin{equation}\label{eq:ts-boundary-hoeffding}
 B_\infty(c;P,Q)=\inf_{S:\KL(S\|P)\le c}\KL(S\|Q).
\end{equation}
Both exponents are optimal among tests with the corresponding level constraint, including tests using all the uncapped data. An alternative-dependent threshold at the crossing $B_\infty(c)=c$ attains the optimal symmetric-error exponent
\[
 C(P,Q)=-\min_{0\le t\le1}\log\sum_i p_i^{1-t}q_i^t
\]
per minority observation. All statements remain true for a fixed positive-weight mixture containing $U$ and any other e-process valid in the same filtration. Interchanging the streams gives the reverse boundary regime.
\end{theorem}
Thus the one-sample limit has significance exponent $\KL(Q\|P)$ and miss exponent $\KL(P\|Q)$. Their opposite orientations persist under arbitrarily fast sample imbalance. The proof in \cref{app:batch-proofs} uses superexponential concentration of the retained reference empirical distribution on the minority scale. Interior exponents, boundary exponents, and local efficiency describe different asymptotic regimes.

\paragraph{Relation to earlier optimal tests.}
The finite-alphabet generalized likelihood-ratio statistic and its large-deviation optimality are classical \citep{gutman1989,unnikrishnan2012}. Kernel-MMD tests already attain optimal fixed-level two-sample miss exponents under the assumptions of \citet{zhu2021}; \citet{harsha2026,grootveld2026} give recent finite-alphabet refinements. These results coexist with smaller log-wealth growth for a particular bounded MMD betting strategy. The contribution here is the connection of the general-space, allocation-adaptive e-process to its optimal attained-significance rates and the explicit coding analysis under each calibration. The corresponding independence paper \citep{ramdas2026independence} develops the mutual-information counterpart; its classical finite-table error-exponent region is also covered by \citet{lapidoth2018} after exchanging the null and alternative conventions.

\section{Finite-sample simulations}\label{sec:experiments}
We compare finite implementations of the constructions with the kernel-MMD betting procedure of \citet{shekhar2024betting}, then study classifier betting and representation adaptation on images. The finite-alphabet and Gaussian experiments use seven alternatives, 400 independent replications per alternative, and a horizon of 1000 pairs. Four additional null experiments use 1000 replications each. All simulations use independent i.i.d. streams, and all methods receive the same observations within a replication. The additive methods also retain the exchangeable-pair validity discussed in \cref{sec:mmd}. Parameters are fixed across the alternatives and their corresponding nulls. The reported level is $\alpha=0.05$; results at $\alpha=0.01$ are included in the companion data.

\subsection{Methods and experimental design}
\paragraph{Kernel-MMD benchmarks.}
The primary benchmark, denoted MMD OGA--ONS, implements Definitions 5 and 8 of the corrected version arXiv:2112.09162v6: online gradient ascent on the RKHS ball of radius $1/2$, followed by online Newton-step betting with fractions in $[-1/2,1/2]$. MMD ERM--ONS replaces the OGA witness by the empirical mean-embedding difference normalized to RKHS norm $1/2$, while using the same ONS update. MMD full grid uses this empirical witness with a uniform wealth mixture over the 41 fractions $\{-1,-0.95,\ldots,0.95,1\}$. This third method isolates the effect of enlarging the betting-fraction range within the additive-payoff family.

\paragraph{Product-bet implementations.}
Kernel logistic uses \cref{sec:kernel} with six radii $R\in\{1,2,4,8,16,32\}$, step sizes $R/(2\sqrt t)$, and the augmented kernel $K=(1+k)/2$. Each learner has a uniform stake mixture over $\{0,1/4,1/2,1\}$, and the six learners have equal initial weights. Partition pair implements the simultaneous-pair Laplace predictor of \cref{sec:partition}. Partition code implements \eqref{eq:imb-code-e}. The finite-alphabet experiments use the identity partition. On the real line, both partition methods mix equally over five fixed standard-normal quantile partitions with $2,4,8,16,32$ cells. In ten dimensions they mix over the 30 coordinatewise quantile partitions having $2,4,8$ cells per coordinate. This marginal representation is suited to the sparse mean shift in that experiment. Joint refining partitions supply the broader approximation property in the universal theorem.

For real-valued vectors, the MMD methods use $k(x,y)=\exp\{-\|x-y\|^2/(2d)\}$, with the same fixed bandwidth choice as the Gaussian setting in \citet{shekhar2024betting}. For finite symbols they use $k(x,y)=\one\{x=y\}$. Kernel logistic uses the corresponding augmented kernel. Both observations in a pair are scored before either updates a learner. Kernel parameters and prior weights are fixed in advance. The finite banks retain anytime validity and the stated class-relative regret guarantees; universal approximation is supplied by the expanding banks in the theoretical constructions.

\paragraph{Distributions.}
The Bernoulli alternatives have success probabilities $(p,q)=(0.4,0.6)$ and $(0.2,0.8)$. The univariate Gaussian alternatives compare $N(0,1)$ with $N(0.5,1)$, $N(1,1)$, and $N(0,1.5^2)$. The ten-dimensional alternative compares $N(0,I_{10})$ with $N(0.5e_1,I_{10})$. The four-symbol alternative has probability vectors $(0.4,0.3,0.2,0.1)$ and $(0.1,0.2,0.3,0.4)$. The nulls use Bernoulli$(1/2)$, $N(0,1)$, $N(0,I_{10})$, and the first four-symbol distribution in both streams. The identity kernel and partition are unchanged by any relabeling of the symbols.

The sample-efficiency estimand is $\E[\min\{\tau_\alpha,1000\}]$, so a run that has not rejected contributes the horizon. Evidence growth is measured by $\log W_{1000}/1000$. Monte Carlo uncertainty is computed across independent replications, and comparisons between methods use paired differences. All methods are run to the horizon even after crossing, in order to measure their subsequent wealth growth. This is a comparison at equal observation budgets; the multi-learner kernel implementation uses more arithmetic per pair than one MMD learner.

\subsection{Balanced alternatives and null calibration}
\begin{table}[tbp]
\centering\small\setlength{\tabcolsep}{3pt}
\begin{tabular}{lrrrrrr}
\toprule
&\shortstack{MMD\\OGA--ONS}&\shortstack{MMD\\ERM--ONS}&\shortstack{MMD\\full grid}&\shortstack{Kernel\\logistic}&\shortstack{Partition\\pair}&\shortstack{Partition\\code}\\
\midrule
Bernoulli $(0.4,0.6)$ & 150.1 (4.5) & 137.9 (4.2) & 133.3 (4.5) & 159.3 (4.8) & 123.4 (4.7) & 179.6 (5.9) \\
Bernoulli $(0.2,0.8)$ & 23.3 (0.4) & 23.0 (0.3) & 19.8 (0.3) & 17.2 (0.4) & 12.3 (0.3) & 14.3 (0.4) \\
Gaussian mean $0.5$ & 147.8 (3.2) & 136.6 (2.9) & 111.6 (3.0) & 130.3 (3.4) & 128.9 (4.6) & 220.4 (6.2) \\
Gaussian mean $1$ & 52.8 (0.8) & 50.5 (0.8) & 40.8 (0.7) & 37.7 (0.9) & 30.7 (0.9) & 53.3 (1.3) \\
Gaussian scale $1.5$ & 169.8 (3.0) & 159.2 (2.9) & 129.7 (2.9) & 142.5 (3.2) & 244.9 (8.6) & 674.0 (11.1) \\
Gaussian, $d=10$ & 325.0 (4.1) & 316.4 (3.9) & 238.8 (3.5) & 214.4 (4.6) & 181.6 (4.9) & 277.9 (6.2) \\
Four symbols & 53.6 (0.9) & 52.3 (0.8) & 43.0 (0.7) & 40.1 (0.9) & 28.7 (0.9) & 49.1 (1.3) \\
\bottomrule
\end{tabular}
\caption{Estimated $\E[\min\{\tau_{0.05},1000\}]$ in pairs, with Monte Carlo standard errors in parentheses, from 400 replications. Unrejected runs contribute 1000. The partition-code power is $89.75\%$ in the scale-change row; every other entry has empirical power $100\%$ at this horizon.}\label{tab:simulation-stopping}
\end{table}

\begin{table}[tbp]
\centering\small\setlength{\tabcolsep}{3pt}
\begin{tabular}{lrrrrrrr}
\toprule
&\shortstack{MMD\\OGA--ONS}&\shortstack{MMD\\ERM--ONS}&\shortstack{MMD\\full grid}&\shortstack{Kernel\\logistic}&\shortstack{Partition\\pair}&\shortstack{Partition\\code}&$2\JS$\\
\midrule
Bernoulli $(0.4,0.6)$ & 0.0338 & 0.0351 & 0.0362 & 0.0348 & 0.0370 & 0.0340 & 0.0403 \\
Bernoulli $(0.2,0.8)$ & 0.1737 & 0.1744 & 0.2894 & 0.3764 & 0.3824 & 0.3804 & 0.3855 \\
Gaussian mean $0.5$ & 0.0340 & 0.0356 & 0.0472 & 0.0491 & 0.0447 & 0.0361 & 0.0606 \\
Gaussian mean $1$ & 0.0839 & 0.0858 & 0.1447 & 0.1993 & 0.1919 & 0.1743 & 0.2228 \\
Gaussian scale $1.5$ & 0.0297 & 0.0314 & 0.0472 & 0.0528 & 0.0303 & 0.0115 & 0.0734 \\
Gaussian, $d=10$ & 0.0140 & 0.0146 & 0.0262 & 0.0387 & 0.0413 & 0.0327 & 0.0606 \\
Four symbols & 0.0887 & 0.0910 & 0.1539 & 0.1995 & 0.2044 & 0.1965 & 0.2129 \\
\bottomrule
\end{tabular}
\caption{Mean $\log W_{1000}/1000$ and the information-theoretic optimum per pair. Standard errors and the growth over pairs 501--1000 are supplied with the simulation data. Gaussian information values are computed by one-dimensional numerical quadrature; the sparse ten-dimensional shift has the same divergence as its affected coordinate.}\label{tab:simulation-growth}
\end{table}

The kernel logistic implementation has a smaller mean rejection time than MMD OGA--ONS in six of the seven settings in \cref{tab:simulation-stopping}. For the unit Gaussian mean shift, the estimates are $37.7$ and $52.8$ pairs, respectively. The paired difference is $-15.15$ pairs, with a Monte Carlo $95\%$ interval $[-16.21,-14.09]$. For the sparse ten-dimensional shift, the estimates are $214.4$ and $325.0$ pairs, with paired difference $-110.62$ and interval $[-116.77,-104.46]$. The mild Bernoulli alternative favors MMD OGA--ONS: kernel logistic takes $159.3$ pairs compared with $150.1$, with difference $9.28$ and interval $[5.44,13.12]$.

\begin{figure}[tbp]
\centering
\includegraphics[width=.86\linewidth]{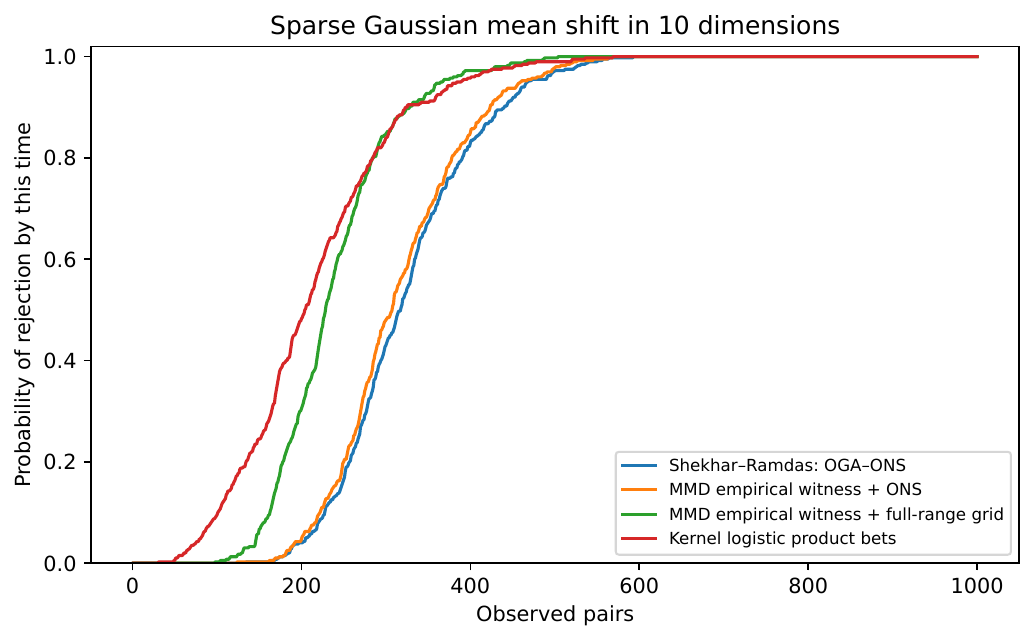}
\caption{Rejection probability by each pair count at level $0.05$, for the sparse ten-dimensional Gaussian mean shift. Curves use the same 400 replications for all methods.}\label{fig:simulation-power}
\end{figure}

The full-range additive grid reduces MMD rejection times in every simulated alternative. Kernel logistic stops sooner than this stronger benchmark in the strong Bernoulli, unit Gaussian shift, sparse ten-dimensional shift, and four-symbol settings. The full-range MMD grid stops sooner in the mild Bernoulli, smaller Gaussian shift, and Gaussian scale-change settings. Thus betting-fraction adaptation accounts for a substantial part of the improvement over the ONS implementation. The logarithmic-loss learner offers further gains when its additional information offsets its early learning cost. \Cref{fig:simulation-power} displays the rejection curves for the sparse shift.

\begin{figure}[tbp]
\centering
\includegraphics[width=.86\linewidth]{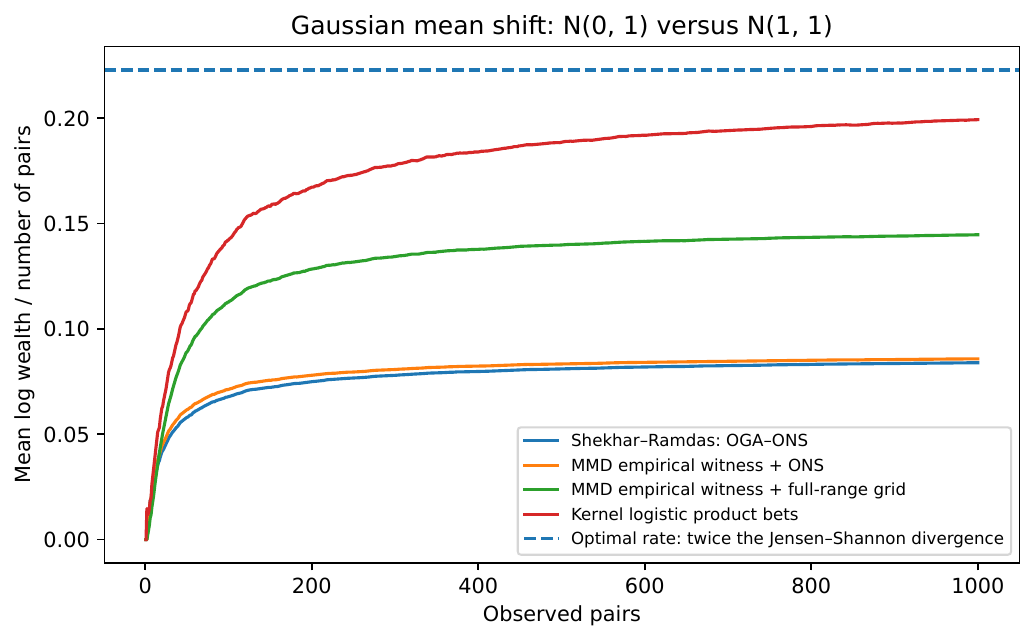}
\caption{Mean log wealth per pair for $N(0,1)$ versus $N(1,1)$. The dashed line is the information bound $2\JS(P,Q)=0.222843$. The estimated rates at 1000 pairs are reported in \cref{tab:simulation-growth}.}\label{fig:simulation-growth}
\end{figure}

Kernel logistic has the largest estimated growth among the four kernel methods in six settings; the full-range MMD grid has slightly larger growth in the mild Bernoulli case. For the unit Gaussian mean shift, the mean normalized log wealth is $0.1993$ for kernel logistic, $0.1447$ for MMD full grid, and $0.0839$ for MMD OGA--ONS, compared with the optimum $0.2228$. \Cref{fig:simulation-growth} shows the growth trajectories. The partition predictor performs particularly well on the small finite alphabets, where it estimates a small number of cell probabilities. The fixed quantile partitions learn the Gaussian scale change more slowly; partition coding also incurs the redundancy associated with its cell counts.

\begin{table}[tbp]
\centering\small\setlength{\tabcolsep}{3pt}
\begin{tabular}{lrrrrrr}
\toprule
&\shortstack{MMD\\OGA--ONS}&\shortstack{MMD\\ERM--ONS}&\shortstack{MMD\\full grid}&\shortstack{Kernel\\logistic}&\shortstack{Partition\\pair}&\shortstack{Partition\\code}\\
\midrule
Bernoulli & 23 & 32 & 31 & 10 & 34 & 10 \\
Gaussian, $d=1$ & 11 & 30 & 26 & 5 & 24 & 1 \\
Gaussian, $d=10$ & 3 & 2 & 18 & 7 & 47 & 5 \\
Four symbols & 21 & 29 & 27 & 13 & 35 & 0 \\
\bottomrule
\end{tabular}
\caption{Numbers of null rejections by 1000 pairs in 1000 independent replications, at level $0.05$. All methods use the same streams within each replication. Wilson intervals for each rate are included in the companion data.}\label{tab:simulation-null}
\end{table}

The null rejection counts in \cref{tab:simulation-null} range from zero to 47 among 1000 replications at nominal level $0.05$. Kernel logistic has counts 10, 5, 7, and 13 in the Bernoulli, univariate Gaussian, ten-dimensional Gaussian, and four-symbol experiments, respectively. The observed conservatism varies with the learner and its aggregation scheme. Anytime validity follows from the constructions; the null simulations check their numerical implementation.

\subsection{An allocation-adaptive implementation under imbalance}
For Bernoulli$(0.4)$ versus Bernoulli$(0.6)$, we additionally use 400 independent stream pairs and evaluate $n_2\in\{100,1000,10000\}$ and $n_1/n_2\in\{1,4,16,64\}$. Each replication uses the same underlying streams at all these endpoints. We compare the all-data identity-partition code, the allocation-adaptive mixture \eqref{eq:adaptive-code-mixture}, the capped and size-protected identity-partition code, and the capped and size-protected kernel pilot of \cref{cor:imb-cap}. Every method uses a single rule across the endpoints.

The cap is $m\lceil\log(m+1)\rceil$. The pilot uses $s=\lfloor\sqrt m\rfloor$ balanced pairs, the training rule \eqref{eq:imb-ogd}, and clipping levels $1,\ldots,\lfloor\log(m+1)\rfloor$, with the omitted clipping weight placed in cash as described above. The kernel $K=(1+\one\{x=y\})/2$ has an exact two-dimensional feature map, so training and evaluation use that finite representation. All clipping weights and sample-size penalties from the theorem are included.

\begin{table}[tbp]
\centering\small\setlength{\tabcolsep}{4pt}
\begin{tabular}{rrrrrrr}
\toprule
$n_2$&$n_1/n_2$&\shortstack{All-data\\code}&\shortstack{All-ratio\\code mixture}&\shortstack{Capped\\code}&\shortstack{Capped\\kernel pilot}&$I_{n_1,n_2}/n_2$\\
\midrule
1,000 & 1 & 0.0349 & 0.0342 & 0.0073 & -0.0096 & 0.0403 \\
1,000 & 4 & 0.0585 & 0.0579 & 0.0281 & -0.0010 & 0.0646 \\
1,000 & 16 & 0.0692 & 0.0685 & 0.0329 & 0.0007 & 0.0762 \\
1,000 & 64 & 0.0721 & 0.0714 & 0.0329 & 0.0007 & 0.0798 \\
\addlinespace
10,000 & 1 & 0.0396 & 0.0395 & 0.0359 & 0.0265 & 0.0403 \\
10,000 & 4 & 0.0637 & 0.0636 & 0.0597 & 0.0445 & 0.0646 \\
10,000 & 16 & 0.0752 & 0.0751 & 0.0684 & 0.0510 & 0.0762 \\
10,000 & 64 & 0.0787 & 0.0786 & 0.0684 & 0.0510 & 0.0798 \\
\bottomrule
\end{tabular}
\caption{Mean log wealth per minority observation for Bernoulli$(0.4)$ versus Bernoulli$(0.6)$, from 400 shared-stream replications. The all-ratio mixture is \eqref{eq:adaptive-code-mixture}. Capped methods include their sample-size correction. The directed KL limit is $D(Q\|P)=0.081093$.}\label{tab:simulation-imbalance}
\end{table}

\begin{figure}[tbp]
\centering
\includegraphics[width=.86\linewidth]{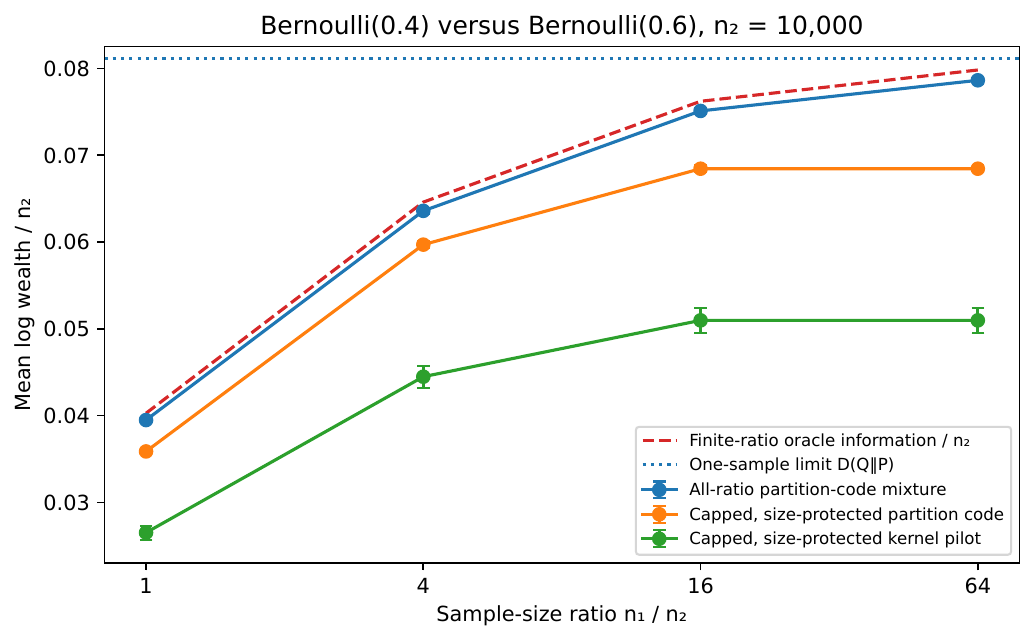}
\caption{Minority-scale growth at $n_2=10000$ for Bernoulli$(0.4)$ versus Bernoulli$(0.6)$. Error bars are pointwise Monte Carlo $95\%$ intervals. The allocation-adaptive coding mixture retains the all-data rate up to a constant wealth factor. Capping, endpoint protection, and the small pilot have visible finite-sample costs.}\label{fig:simulation-imbalance}
\end{figure}

At $(n_1,n_2)=(640000,10000)$, the allocation-adaptive coding mixture has mean normalized log wealth $0.0786$, close to the finite-ratio oracle information $0.0798$ and the limiting divergence $D(Q\|P)=0.0811$. Its normalized log wealth is at least the all-data value minus $\log2/n_2$, and it inherits the arbitrary-imbalance guarantee from its capped component. The capped code alone has mean $0.0684$, and the capped kernel pilot has mean $0.0510$. Their cap retains 100000 reference observations, so the two largest original ratios give identical retained data. These finite-sample effects are visible in \cref{fig:simulation-imbalance}.

The small pilot and the generic sample-size correction produce appreciable overhead at $n_2=1000$, especially for the kernel procedure. At $n_2=10000$ their costs are smaller on the minority scale. The all-data/capped coding mixture illustrates how a practical component can preserve finite-sample efficiency while a universal component supplies the limiting guarantee. Effective kernel implementations call for further work on training allocation, complexity adaptation, and data reuse.

\FloatBarrier
\subsection{Image classification and sequential fine-tuning}\label{sec:image-results}
We apply \cref{sec:finetuning} to the 1797-image UCI optical-digits subset distributed with scikit-learn \citep{uci1998,sklearn-digits}. A fixed stratified split reserves 900 images for ten-digit CNN pretraining and 897 for independently resampled testing streams. The $8\times8$ images are resized to $28\times28$; noisy settings add independent Gaussian pixel noise with standard deviation $0.5$ and clip intensities to $[0,1]$. The experiments therefore concern controlled distributions of real digit images with 784 pixel coordinates. Pretraining observations are external to the reported streaming sample sizes. Monte Carlo uncertainty is conditional on the fixed image population and pretrained checkpoint. 
Architectures, preprocessing, source distributions, and optimization settings are given in Appendix~\ref{app:image-experiments}.

\paragraph{Classifier betting and MMD on a shared representation.}
The first image study compares pixel MMD, MMD on the pretrained CNN features, a logistic source classifier on those features, and a raw-pixel neural classifier. It uses 100 replications for each of four alternatives and 200 for each of two nulls, with horizon 1000 pairs. The source heads are refitted every 50 pairs. In the weaker noisy digit-frequency shift, mean truncated rejection times at $\alpha=0.05$ are $976.6$ pairs for pixel MMD OGA--ONS, $248.7$ for full-range MMD on CNN features, and $235.8$ for the CNN product bet. The paired difference between the last two is $-12.9$ pairs with interval $[-32.0,6.2]$. Learned features account for much of the improvement over pixel MMD. Product betting has smaller mean rejection times than feature MMD in the other three settings, with paired intervals excluding zero. \Cref{tab:image-stopping,tab:image-paired} also compare product and additive bets using identical source predictions.

\paragraph{Fine-tuning on the two streams.}
The second study uses the same pretrained CNN and image population, with four restricted-class alternatives: noisy 3 versus 8; noisy 4 versus 9; 3 versus 8 after rotating both streams by $90$ degrees; and the rotated digits mixed in proportions $(0.7,0.3)$ under $P$ and $(0.3,0.7)$ under $Q$. The null uses the same rotated $(0.5,0.5)$ mixture in both streams. During adaptation, only source labels enter training.

The frozen and fine-tuned arms use the same $L^2$-regularized logistic-head fitting procedure. After 50 neutral warm-up pairs, both update every 25 pairs. The fine-tuned arm takes ten balanced replay gradient steps through both convolutional layers and the feature layer, then refits its logistic head. Each arm mixes uniformly over temperatures $\{0.5,1,2,4\}$ and fractions $\{0,1/4,1/2,1\}$. There are 100 paired replications per alternative and 200 null replications, each continued to 500 pairs. This experiment compares encoder adaptation while holding the head-fitting rule and wealth aggregation fixed.

\begin{table}[tbp]
\centering\small\setlength{\tabcolsep}{5pt}
\begin{tabular}{@{}lrrrrr@{}}\toprule
&\multicolumn{2}{c}{Mean $\log W_{500}/500$}&\multicolumn{3}{c}{Mean $\tau_{0.05}\wedge500$}\\
\cmidrule(lr){2-3}\cmidrule(l){4-6}
Setting & Frozen & Tuned & Frozen & Tuned & Mixture\\\midrule
Noisy 3 versus 8 & 1.1018 & 1.1308 & 54.97 & 54.88 & 54.93\\
Noisy 4 versus 9 & 1.2196 & 1.2278 & 54.05 & 54.01 & 54.03\\
Rotated 3 versus 8 & 0.7894 & 1.0293 & 56.82 & 56.14 & 56.33\\
Rotated mixtures & 0.0830 & 0.0896 & 96.87 & 100.96 & 96.80\\
\bottomrule\end{tabular}
\caption{Sequential fine-tuning on UCI digit images, with 100 replications per alternative. Both arms start from the same pretrained encoder and use the same logistic-head fit. All streaming costs, including the 50-pair warm-up, are counted. Mixture denotes \eqref{eq:frozen-tuned-mixture}. Paired Monte Carlo intervals appear in \cref{tab:finetune-paired}.}\label{tab:finetune-main}
\end{table}

\begin{figure}[tbp]
\centering\includegraphics[width=.86\linewidth]{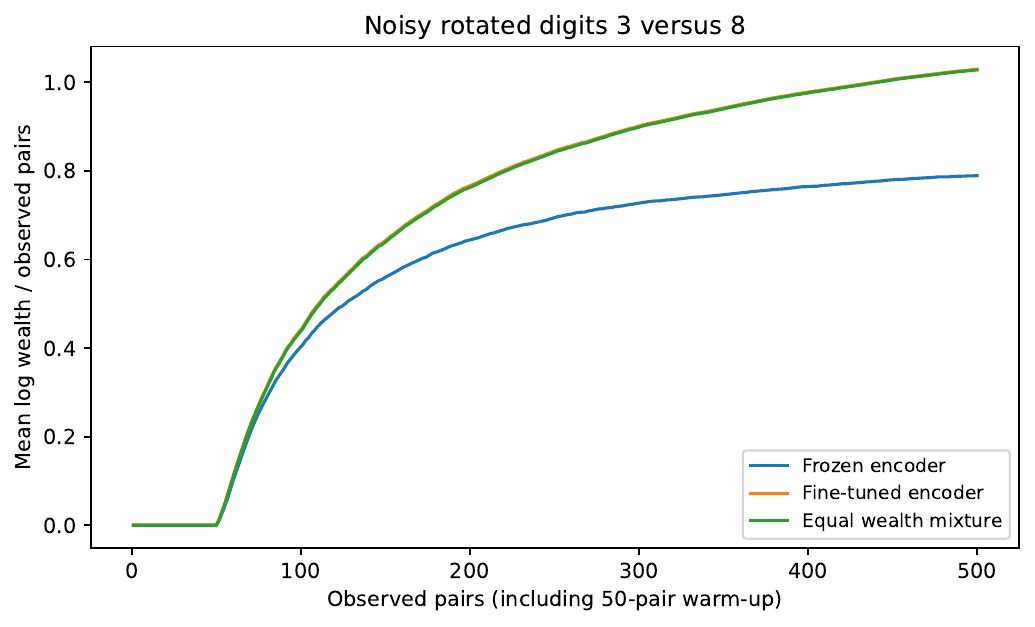}
\caption{Mean log wealth per pair for noisy, rotated digit 3 versus digit 8. Both source classifiers learn from earlier observations; only the fine-tuned arm updates the pretrained encoder. Trials continue to 500 pairs after any initial rejection. The 50-pair neutral warm-up is included.}\label{fig:finetune-growth}
\end{figure}

The pretrained features already separate the unrotated pure classes well, so their rejection times are close to the warm-up boundary. For rotated 3 versus 8, fine-tuning raises mean normalized log wealth from $0.7894$ to $1.0293$, an increase of about $30\%$; the paired growth difference has $95\%$ Monte Carlo interval $[0.2337,0.2462]$ nats per pair. \Cref{fig:finetune-growth} shows this adaptation to the shared rotation shift. The interval is descriptive and unadjusted for the other comparisons.

For the weaker rotated mixture, fine-tuning increases terminal growth but takes $4.09$ more pairs to reject on average, with paired interval $[0.78,7.40]$. The equal wealth mixture has mean rejection time $96.80$, compared with $96.87$ for the frozen arm, while retaining the pathwise guarantee \eqref{eq:frozen-tuned-mixture}. The role of temperature and stake adaptation is examined separately in Appendix~\ref{app:finetune-calibration}.

All three methods reject in all 100 replications of each fine-tuning alternative by 500 pairs. Under the null, the frozen, fine-tuned, and mixture methods reject in $1$, $3$, and $3$ of 200 replications at nominal level $0.05$. Across the two studies, the image experiments contain 1400 trials. Appendix~\ref{app:image-experiments} gives their full protocols, comparison tables, additional rejection curves, and null checks.
\FloatBarrier

\paragraph{Reproducibility.}
The ancillary package contains the experiment programs, all per-replication outcomes, aggregate curves, Monte Carlo intervals, and deterministic random seeds. Main alternative, null, and imbalance seeds are 20261105, 20261106, and 20261107, respectively; trial and scenario identifiers spawn separate NumPy random streams. An independent ten-replication pilot uses seed 731947. Parameters and the seven alternative scenarios are unchanged between the pilot and the main comparison. Optimized routines are checked against the manuscript's reference implementations, including frozen-pair predictions, direct finite-feature MMD updates, and the kernel pilot. Results at the additional significance level and intermediate imbalance endpoints are included with the data. 

\FloatBarrier

\subsection{Fixed-sample calibration: local power and attained significance}\label{sec:batch-numerics}
We evaluate the identity-partition code and its exact probability-ordered permutation calibration on Bernoulli observations. Write $a=rb$, with $r\in\{1,4\}$. Under the local alternatives the success probabilities are
\[
 p_{a,b}=\tfrac12-\tfrac12\delta_{a,b},\qquad
 q_{a,b}=\tfrac12+\tfrac12\delta_{a,b},\qquad
 \delta_{a,b}=2\sqrt{\tfrac14(a^{-1}+b^{-1})}.
\]
The noncentrality in \cref{thm:ts-pitman} is four, so the classical local-power limit at level $0.05$ is $0.516005$. Each configuration uses 8000 independent count pairs. The exact permutation p-value sums the full hypergeometric support in logarithmic coordinates. A likelihood-ratio comparator uses the $\chi^2_1$ cutoff. All three methods receive the same count pairs.

\begin{table}[tbp]\centering\small
\begin{tabular}{rrrrr}\toprule
$a/b$ & $b$ & Direct code (\%) & Exact permutation (\%) & Likelihood ratio (\%)\\\midrule
1 & 100 & 3.79 & 48.48 & 54.41 \\
1 & 400 & 1.79 & 49.73 & 52.73 \\
1 & 1,600 & 0.79 & 50.44 & 51.74 \\
1 & 6,400 & 0.31 & 51.61 & 52.30 \\
4 & 100 & 2.36 & 48.38 & 51.91 \\
4 & 400 & 1.04 & 48.95 & 51.32 \\
4 & 1,600 & 0.41 & 50.29 & 51.19 \\
4 & 6,400 & 0.14 & 51.08 & 51.58 \\
\bottomrule\end{tabular}
\caption{Fixed-endpoint rejection proportions under local alternatives, from 8000 replications per row. The common asymptotic target for the exact conditional and likelihood-ratio tests is $51.60\%$. Monte Carlo standard errors for these powers are at most $0.56$ percentage points. The direct code's local power vanishes.}\label{tab:ts-local-power}
\end{table}

The fixed-alternative experiment uses Bernoulli$(0.4)$ versus Bernoulli$(0.6)$, 1000 replications per configuration, $b\in\{100,1000,10000\}$, and the same two ratios. At $(a,b)=(10000,10000)$, the mean values of $-\log p/N$ are $0.019697$ for reciprocal-code calibration and $0.020295$ for exact permutation calibration; the target is $G_{1/2}=0.020136$. At $(40000,10000)$ they are $0.012738$ and $0.012999$, with target $G_{4/5}=0.012918$. These two calibrations have the same fixed-alternative exponent but markedly different local power. Full summaries, seed 2026111207, and the exact-enumeration implementation accompany the paper. Ten deterministic checks cover the coding sandwich, null expectation, hypergeometric ordering, local constant, and KL orientations.

\section{Discussion}
For independent pairs, the oracle multiplier is a product of two class probabilities, and its log-growth shortfall is exactly an excess logarithmic prediction loss. Partition-based Laplace prediction and kernel logistic learning recover the Bayes risk. Their pathwise regret bounds yield almost-sure and expected log-growth optimality, as well as first-order efficiency of rejection times. A characteristic kernel with a constant feature provides the required approximation on general measurable data spaces.

The batch calibrations separate three regimes: attained significance is governed by forward information, miss probabilities by a reverse projection, and local power by the calibration after removing logarithmic learning cost. On finite alphabets the exact conditional version recovers a classical test, while the direct code supplies sequential validity. These connections inherit the established empirical-divergence benchmarks.

For unequal samples, the weighted information bound interpolates between the two-sample and one-sample problems. When the reference sample is much larger, the attainable rate per test observation is $\KL(Q\|P)$. Universal coding and pilot-trained kernel scores achieve this limit, including infinite divergence. A slowly growing cap on the reference sample preserves the rate and limits both computation and sample-size penalties.

Neural source classifiers extend the predictive interface to learned image representations. Predict-then-train updates preserve validity while allowing both the encoder and head to adapt. Their improvement in raw log wealth equals their reduction in cumulative predictive logarithmic loss; representation and head errors are separated by \cref{prop:representation-gap}. In the digit experiments, fine-tuning yields larger terminal wealth under a shared rotation shift, while the weaker mixture reveals an early learning cost. A frozen/tuned wealth mixture retains either arm's evidence within a constant factor.

The learning rule, model mixture, and betting fraction provide complementary forms of adaptation. Norm and kernel mixtures enrich the score class; universal portfolios compete with the best fixed fraction at logarithmic cost. The exchangeable-pair construction gives a second choice of null model, with its own information bound. Exact kernel learning and continuous portfolios have quadratic cumulative cost, while finite fraction grids containing full betting preserve the main asymptotic rates at lower computational cost. Sparse kernel approximations with sublinear excess log loss, effective complexity adaptation, and data-adaptive allocation between the two streams are natural directions for further work. The general-space and all-allocation guarantees provide a foundation for this algorithmic development. The simulations show that larger long-run growth can coexist with additional early learning cost, so practical design benefits from attention to both wealth growth and rejection times.

\appendix
\section{A kernel construction attaining the exchangeable-pair optimum}\label{app:swap}
For applications requiring validity for conditionally exchangeable pairs, replace the product multiplier by
\begin{equation}\label{eq:ranking-factor}
E_t^{\mathrm{sw}}=2\sig(h_t(X_t)-h_t(Y_t)).
\end{equation}
It obeys $E_t^{\mathrm{sw}}(x,y)+E_t^{\mathrm{sw}}(y,x)=2$ for every score and is consequently conditionally fair under every exchangeable pair law. Learn the score with the convex ranking loss
\[
\ell_t^{\mathrm{sw}}(h)=\spplus(h(Y_t)-h(X_t)),
\quad
\nabla\ell_t^{\mathrm{sw}}(h)
=\{\sig(h(X_t)-h(Y_t))-1\}\{K(X_t,\cdot)-K(Y_t,\cdot)\}.
\]
Its gradient norm is at most two. The proof of \cref{thm:kernel-regret} therefore gives $6R\sqrt m$ regret on a radius-$R$ ball, now with $\log2$ replacing $\log4$ in the wealth identity.

\begin{proposition}\label{prop:swap-kernel}
Under the kernel assumptions of \cref{thm:kernel-optimal}, mixtures of the ranking learners over unbounded radii, with the same activation and stake options, satisfy
\[
n^{-1}\log W_n^{\mathrm{sw}}\longrightarrow G_{\mathrm{sw}}(P,Q)\as
\]
under independent i.i.d. pairs $P\otimes Q$, while remaining anytime valid for the conditional exchangeable-pair null.
\end{proposition}
\begin{proof}
It remains to show that RKHS scores approximate the oracle ranking risk. Let $h^*=\log(p/q)$ with extended-real values and $h_B^*=\clip(h^*,-B,B)$. Under $P\otimes Q$, $h^*(X)$ is never $-\infty$ and $h^*(Y)$ is never $+\infty$, so their difference is well-defined almost surely. It is the log odds of the oracle orientation probability in \eqref{eq:swap-oracle}. The difference of clipped scores converges to it, possibly to $+\infty$.

If this oracle difference is nonnegative, the log multiplier obtained by clipping is between zero and $\log2$. If it is negative, monotonicity and contraction of clipping place the clipped difference between the original difference and zero. The clipped log multiplier is then between the integrable oracle log multiplier and zero. Dominated convergence proves that the expected log multiplier from $h_B^*$ tends to $G_{\mathrm{sw}}$. For bounded measurable scores $h,g$, the 1-Lipschitz property of $v\mapsto\log\sig(v)$ gives
\[
\left|\E\log E^{\mathrm{sw}}(h)-\E\log E^{\mathrm{sw}}(g)\right|
\le\E_P|h-g|+\E_Q|h-g|
\le2\norm{h-g}_{L^2(M)}.
\]
RKHS density therefore permits approximation of each clipped score. The regret, activation, and countable-comparator argument of \cref{thm:kernel-optimal} gives the desired lower bound, and \cref{thm:swap} supplies the upper bound. Validity was established directly from the swap identity.
\end{proof}
The ranking learner minimizes pair-orientation log loss and uses a nonlinear swap-normalized payoff to attain the exchangeable-pair information bound.

\section{A finite-space illustration of the two benchmarks}\label{app:finite}
Let $\X=\{0,1\}$, $P=\operatorname{Bernoulli}(p)$ and $Q=\operatorname{Bernoulli}(q)$. Then
\[
G(P,Q)=2H\left(\frac{p+q}{2}\right)-H(p)-H(q).
\]
Writing $a=p(1-q)$ and $b=(1-p)q$, the exchangeable-pair benchmark is
\[
G_{\mathrm{sw}}(P,Q)
=a\log\frac{2a}{a+b}+b\log\frac{2b}{a+b},
\]
with zero terms interpreted by continuity. Equal-coordinate pairs contribute zero to the latter expression because their orientation is uninformative. At $p=1,q=0$, these formulas give $\log4$ and $\log2$, respectively. At $p=q$, both are zero. 

\section{Adaptive representations and neural growth guarantees}\label{app:adaptive-proofs}
This appendix develops the learning guarantees used in \cref{sec:finetuning}. Let $\X$ and the feature space be standard Borel, let $M=(P+Q)/2$, and write $\eta^*=dP/d(P+Q)$. The predictor at time $t$ is fixed conditionally on $\F_{t-1}$, including its representation, normalization, and any algorithmic randomness. Independent pretraining data may be included in the initial sigma-field. Under the i.i.d. alternative, future testing observations are independent of this initial data and the preceding pairs.

\subsection{Representation and prediction error}
For a fixed measurable feature map $\phi$, denote its pushforward laws by $P_\phi,Q_\phi$ and set $M_\phi=(P_\phi+Q_\phi)/2$. Define a measurable feature posterior by
\[
 \bar\eta=\frac{dP_\phi}{d(P_\phi+Q_\phi)}.
\]
Then $\bar\eta(\phi(Z))=\E_M[\eta^*(Z)\mid\sigma(\phi(Z))]$ almost surely. The representation can be learned from preceding observations; all statements below then hold conditionally on that past for the resulting fixed map.

\begin{proposition}[Representation and head decomposition]\label{prop:representation-gap}
For any measurable $g$ taking values in $(0,1)$ and predictor $\eta=g\circ\phi$,
\begin{align}
 &G(P,Q)-\E_{P\otimes Q}\log\{4\eta(X)(1-\eta(Y))\}\label{eq:representation-decomp}\\
 &\quad=G(P,Q)-G(P_\phi,Q_\phi)
       +2\int\kl\{\bar\eta(\phi(z))\|g(\phi(z))\}\,dM(z).\nonumber
\end{align}
The representation loss satisfies
\begin{equation}\label{eq:representation-information}
 G(P,Q)-G(P_\phi,Q_\phi)
 =2\int\kl\{\eta^*(z)\|\bar\eta(\phi(z))\}\,dM(z)\ge0.
\end{equation}
It is zero exactly when $\eta^*$ is $\sigma(\phi)$-measurable up to $M$-null sets. The last term in \eqref{eq:representation-decomp} may be infinite.
\end{proposition}
\begin{proof}
Introduce a source variable $C\in\{0,1\}$ with equal priors and conditional image law $P$ when $C=1$ and $Q$ when $C=0$. The posterior given $Z$ is $\eta^*(Z)$, and the posterior given $\phi(Z)$ is $\bar\eta(\phi(Z))$. Let $H(u)=-u\log u-(1-u)\log(1-u)$. Conditioning on $\phi(Z)$ in the nonnegative cross-entropy terms gives
\begin{align*}
 \E_M\kl(\eta^*\|g\circ\phi)
 &=\E_M\bigl[-\bar\eta\log(g\circ\phi)
       -(1-\bar\eta)\log(1-g\circ\phi)\bigr]-\E_M H(\eta^*)\\
 &=\E_M\kl(\bar\eta\|g\circ\phi)
       +\E_M H(\bar\eta)-\E_M H(\eta^*),
\end{align*}
where $\bar\eta$ in these expressions is evaluated at $\phi(Z)$. The entropy terms are bounded by $\log2$, so the calculation is valid even if the cross-entropy is infinite. The same conditioning argument gives
\[
 \E_M\kl(\eta^*\|\bar\eta\circ\phi)
 =\E_M H(\bar\eta)-\E_M H(\eta^*).
\]
For posterior values zero or one, the corresponding original posterior agrees almost surely on that feature event, and the usual zero-term conventions apply. The Bayes-risk identity \eqref{eq:bayes-risk} for $P,Q$ and their feature laws now proves \eqref{eq:representation-information}. Combining the first identity with \cref{prop:calibration} proves \eqref{eq:representation-decomp}. Nonnegativity and the equality condition follow from the corresponding facts for Bernoulli relative entropy.
\end{proof}

A pretrained representation can discard information useful for a particular source distinction. For example, let $Z=(U,V)$, let $U$ have the same law in both streams, and let $V=1$ under $P$ and $V=0$ under $Q$. The feature map $\phi(Z)=U$ has $G(P_\phi,Q_\phi)=0$, whereas a feature map retaining $V$ has rate $\log4$. A trainable representation can express the latter distinction. Whether an optimizer learns it from the available source observations is measured by its prequential logarithmic loss. When the pretrained features already preserve the relevant posterior, fine-tuning instead changes the finite-class head approximation, estimation error, and calibration.

\Needspace{15\baselineskip}
\subsection{A sufficient accumulated-loss condition}
The following condition concerns the predictive losses of a chosen network and optimizer and permits nonconvergent network parameters.

\begin{proposition}[Optimal growth from vanishing average predictive shortfall]\label{prop:neural-growth}
Let $\eta_t$ be any predictable source predictor satisfying
\[
 \varepsilon_t\le\eta_t(z)\le1-\varepsilon_t
 \quad\text{for all }z,
 \qquad \varepsilon_t=\frac1{t+2}.
\]
Under fixed independent i.i.d. streams with laws $P,Q$, put
\[
 r_t=2\int\kl\{\eta^*(z)\|\eta_t(z)\}\,dM(z).
\]
If $n^{-1}\sum_{t=1}^n r_t\to0$ almost surely, then the product wealth satisfies
\[
 \frac{\log W_n^\eta}{n}\longrightarrow G(P,Q)
 \qquad\text{almost surely}.
\]
If additionally $n^{-1}\sum_{t=1}^n\E r_t\to0$, then the normalized expected log wealth converges to $G(P,Q)$ as well. The null test-supermartingale guarantee holds for every such learning rule regardless of these learning conditions.
\end{proposition}
\begin{proof}
Let $Z_t=\log\{4\eta_t(X_t)(1-\eta_t(Y_t))\}$ and
$D_t=Z_t-\E[Z_t\mid\F_{t-1}]$. The clipping bounds place $Z_t$ in an interval of length
\[
 2\log\frac{1-\varepsilon_t}{\varepsilon_t}=2\log(t+1).
\]
Hence $\E D_t^2\le4\log^2(t+1)$, and
$\sum_t\E D_t^2/t^2<\infty$. The martingale $\sum_{t=1}^nD_t/t$ is bounded in $L^2$ and converges almost surely. Kronecker's lemma gives $n^{-1}\sum_{t=1}^nD_t\to0$ almost surely. By \cref{prop:calibration},
\[
 \frac{\log W_n^\eta}{n}
 =G(P,Q)-\frac1n\sum_{t=1}^n r_t
       +\frac1n\sum_{t=1}^nD_t,
\]
which proves almost-sure convergence. Taking expectations in the same identity proves the expected-log assertion under its stated condition. The null argument is the conditional calculation in \cref{sec:finetuning}.
\end{proof}

The clipping condition can be imposed after any source learner. Clipping a prediction to $[\varepsilon_t,1-\varepsilon_t]$ increases the loss for either label by at most $-\log(1-\varepsilon_t)$. Its added paired loss is therefore at most $2\log\{(t+2)/(t+1)\}$, whose sum through time $n$ is $2\log\{(n+2)/2\}=o(n)$. Thus clipping preserves a vanishing average excess-log-loss guarantee of an unclipped learner. In contrast, a fixed clipping threshold retains a fixed approximation constraint when source posteriors approach zero or one.

\subsection{Mixtures, data reuse, and unequal blocks}
Let $W^F,W^T$ be the processes obtained by frozen and trainable encoders, each possibly using a fixed temperature/fraction mixture. Conditional Tonelli gives validity of \eqref{eq:frozen-tuned-mixture}; independence of the two learning algorithms is unnecessary. Since $W_t^{\rm mix}\ge W_t^F/2$ and $W_t^{\rm mix}\ge W_t^T/2$, its rejection time obeys
\[
 \tau_\alpha^{\rm mix}
 \le\min\{\tau_{\alpha/2}^F,\tau_{\alpha/2}^T\}.
\]
This is the stopping-time consequence of the factor-two wealth guarantee. It compares with constituent thresholds $2/\alpha$. The mixture may cross the threshold $1/\alpha$ later than either constituent crosses that same threshold.

A fixed mixture $(1-\epsilon)W^{\rm practical}+\epsilon W^{\rm universal}$, with $0<\epsilon<1$, is at least $(1-\epsilon)W^{\rm practical}$ and $\epsilon W^{\rm universal}$. When the processes share the same filtration and null, the universal lower growth bound and the information upper bound prove universal optimality of the mixture. For example, $\epsilon=0.01$ costs at most $-\log(0.99)=0.010051$ nats relative to the practical component. Maintaining the universal component also has its own computational cost; delayed activation gives the finite-at-each-time construction of \cref{sec:activation}.

For unequal prospective blocks, let a balanced-prior source predictor $\eta$ be frozen before a block containing $a$ fresh $P$ observations and $b$ fresh $Q$ observations. Put $w=a/(a+b)$ and $v=1-w$, and define
\[
 \eta_w(z)=\frac{w\eta(z)}{w\eta(z)+v\{1-\eta(z)\}}.
\]
The block factor is
\[
 \prod_{i=1}^a\frac{\eta_w(X_i)}w
 \prod_{j=1}^b\frac{1-\eta_w(Y_j)}v.
\]
Under a common law $R$, its conditional expectation is $(c/w)^a\{(1-c)/v\}^b\le1$, where $c=\E_R\eta_w$. This follows by weighted AM--GM and is the classifier form of \cref{thm:imb-oracle}. Block sizes and the model are fixed prospectively; the model may be updated after the block is scored. Multiplying these factors supplies validity at the block-release times. The sample-size-protected constructions of \cref{sec:imbalance} provide the corresponding endpoint-selection guarantees. The image experiments use balanced pairwise monitoring.

\FloatBarrier
\section{Image experiments: protocols and additional results}\label{app:image-experiments}
The two image studies in \cref{sec:image-results} address different algorithmic choices. The first compares image representations and payoff rules, including MMD on the same features as the source classifier. The second compares frozen and adapted encoders with a common source-head fit and betting scheme. Their update schedules and alternatives differ, and their outcomes are reported separately.

\subsection{Shared image population and independent pretraining}\label{app:image-data}
The experiments use the 1797-image optical-digits subset bundled with scikit-learn \citep{uci1998,sklearn-digits}. The original images have $8\times8$ integer-valued pixels in $\{0,\ldots,16\}$. A stratified split with seed 90173 reserves 900 images for representation learning and 897 for the resampling population. The two index sets are disjoint. Every testing stream samples independently with replacement from the specified held-out digit pools. Simulation uncertainty is conditional on this fixed population and one pretrained checkpoint.

Pixel intensities are divided by 16 and bilinearly resized to $28\times28$. Resizing preserves the original image information. In noisy experiments, fresh independent $N(0,0.5^2)$ pixel perturbations are added to every sampled image, followed by clipping to $[0,1]$. The resulting observations have 784 coordinates with independent high-dimensional nuisance noise. These are UCI optical-digit experiments. The implementation also supports a supplied MNIST archive, but all results reported here use the specified UCI population.

The pretrained encoder has two convolution--ReLU--pool stages with 8 and 16 channels, followed by a 32-dimensional ReLU feature layer and a ten-class head. It is trained on the 900 reserved images and their digit labels for 35 epochs, using AdamW with learning rate $0.002$, weight decay $0.001$, and batches of at most 128. Training samples noise levels from $\{0,0.25,0.5\}$. The final checkpoint is used without selection on source-testing performance. Its diagnostic digit accuracy on the held-out population is $95.76\%$ for clean images and $90.86\%$ in one noisy evaluation. These diagnostics and the held-out digit labels do not enter source-classifier fitting.

Feature means and standard deviations are determined from reserved pretraining images; standard deviations are bounded below by $0.1$. These centering and scaling constants are fixed throughout both studies, including after encoder fine-tuning. The 900 semantic-label pretraining images are additional external data and are not included in the reported streaming sample counts. Only source labels are supplied to the subsequent classifiers.

\begin{figure}[tbp]
\centering\includegraphics[width=.94\linewidth]{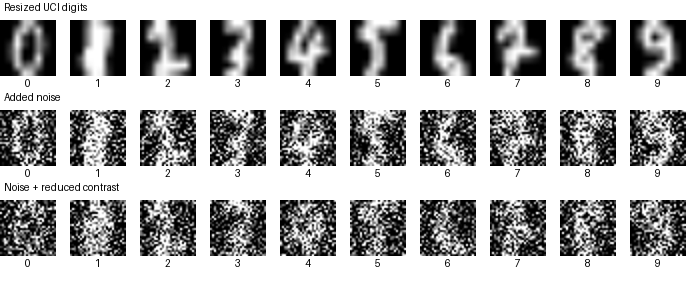}
\caption{Held-out optical-digit examples after resizing (top), independent pixel noise (middle), and contrast reduction followed by noise (bottom). All image experiments use this UCI subset and its fixed train/population split.}\label{fig:image-examples}
\end{figure}

\subsection{Fixed-representation classifier and MMD comparisons}\label{app:image-fixed}
\paragraph{Distributions.}
A digit is sampled uniformly within an even or odd group, and an image is then sampled uniformly from that digit's held-out pool. Stream $P$ assigns probability $0.5$ to the even group. The four alternatives are: clean images with $Q$'s even probability $0.75$; noisy images with the same probability change; noisy images with $Q$'s even probability $0.90$; and noisy images with unchanged digit frequencies but a $20\%$ reduction in $Q$'s contrast before noise is added. The two nulls use identical clean or noisy distributions.

\paragraph{Learners.}
A logistic source head with regularization parameter $C=0.1$ is refitted every 50 pairs on all preceding source-labeled CNN features. The first 50 pairs receive neutral classifier bets. A second source learner is a raw-pixel $784$--$64$--$1$ ReLU network with no external digit pretraining. At each 50-pair update it takes 20 balanced AdamW steps, using 64 earlier images per source per step, learning rate $0.001$, and weight decay $0.01$.

The product-bet arms mix uniformly over temperatures $\{0.5,1,2,4\}$ and fractions $\{0,1/4,1/2,1\}$. An additive control uses the identical CNN source logits and temperatures, with multipliers
\[
 1+\lambda\{\eta_{t,T}(X_t)-\eta_{t,T}(Y_t)\},
 \qquad\lambda\in\{-1,-0.95,\ldots,0.95,1\},
\]
and a uniform wealth mixture over temperatures and fractions. This comparison keeps the predictions fixed while changing the payoff and its stated fraction grid.

The pixel-MMD OGA--ONS method implements the online-gradient witness and online-Newton betting rule described in \cref{sec:experiments}. The stronger full-range MMD comparator uses the empirical mean-embedding witness and the same 41-point fraction grid in $[-1,1]$. A feature-MMD comparator applies that full-range method to the same 32-dimensional CNN features used by the logistic classifier. Each MMD procedure mixes equally over Gaussian-kernel bandwidths at $0.5$, $1$, and $2$ times a reference bandwidth determined from reserved pretraining images. Precomputed Gram matrices are used only for arithmetic efficiency: each witness accesses previously observed indices and fixed kernel parameters.

\paragraph{Replication.}
There are 100 paired replications per alternative and 200 per null, each continued to 1000 pairs. Every method receives identical streams within a replication. Rejection times include the initial training period, and nonrejections contribute 1000 to the truncated-time mean. Results at $\alpha=0.05$ are displayed below; $\alpha=0.01$ outcomes are included in the ancillary data. The independent two-replication-per-scenario pilot uses seed 731947 and the main run seed 2026110601. Configurations are unchanged between pilot and main runs.

\begin{table}[tbp]
\centering\small\setlength{\tabcolsep}{5pt}
\begin{tabular}{@{}lrrrr@{}}\toprule
Method & Clean 0.75 & Noisy 0.75 & Noisy 0.90 & Contrast\\\midrule
Pixel MMD OGA--ONS & 410.2 & 976.6 & 586.3 & 281.6\\
Pixel MMD, full fractions & 295.9 & 836.3 & 465.8 & 239.4\\
Feature MMD, full fractions & 244.8 & 248.7 & 133.1 & 89.9\\
CNN source head, product & 220.9 & 235.8 & 111.6 & 60.4\\
Same classifier, additive & 253.1 & 275.6 & 131.2 & 71.3\\
Raw-pixel network, product & 279.0 & 577.5 & 165.0 & 197.1\\
\bottomrule\end{tabular}
\caption{Mean $\tau_{0.05}\wedge1000$ in pairs for the first image study. The first three columns change $Q$'s even-digit probability from $P$'s value $0.5$; the last changes contrast at equal digit frequencies. Each alternative has 100 replications. Monte Carlo standard errors and both significance levels are retained in the ancillary summaries.}\label{tab:image-stopping}
\end{table}

Every method rejects in every alternative replication except in the weaker noisy frequency shift. There, pixel MMD--ONS rejects in $28\%$ of trials, full-range pixel MMD in $92\%$, and the raw-pixel source network in $64\%$. Feature MMD and both CNN-classifier bets reject in all 100 trials. Representation learning accounts for much of the improvement over pixel MMD. The paired comparisons in \cref{tab:image-paired} give the remaining differences at a common feature representation or a common source predictor.

\begin{table}[tbp]\centering\small
\begin{tabular}{@{}llrr@{}}\toprule
Scenario & Reference & Difference & Paired $95\%$ interval\\\midrule
Clean 0.75 & Feature MMD & $-24.0$ & $[-39.7,-8.2]$\\
Clean 0.75 & Same classifier, additive & $-32.3$ & $[-42.8,-21.7]$\\
Noisy 0.75 & Feature MMD & $-12.9$ & $[-32.0,6.2]$\\
Noisy 0.75 & Same classifier, additive & $-39.8$ & $[-52.5,-27.2]$\\
Noisy 0.90 & Feature MMD & $-21.5$ & $[-27.8,-15.2]$\\
Noisy 0.90 & Same classifier, additive & $-19.6$ & $[-23.2,-16.0]$\\
Contrast & Feature MMD & $-29.5$ & $[-31.4,-27.6]$\\
Contrast & Same classifier, additive & $-10.9$ & $[-11.5,-10.2]$\\
\bottomrule\end{tabular}
\caption{CNN product bet minus reference in mean truncated rejection time. Negative differences favor product betting. Intervals use paired Student-$t$ Monte Carlo calculations conditional on the fixed image population and pretrained model; they are descriptive and unadjusted for multiple comparisons.}\label{tab:image-paired}
\end{table}

\begin{figure}[tbp]\centering
\includegraphics[width=.86\linewidth]{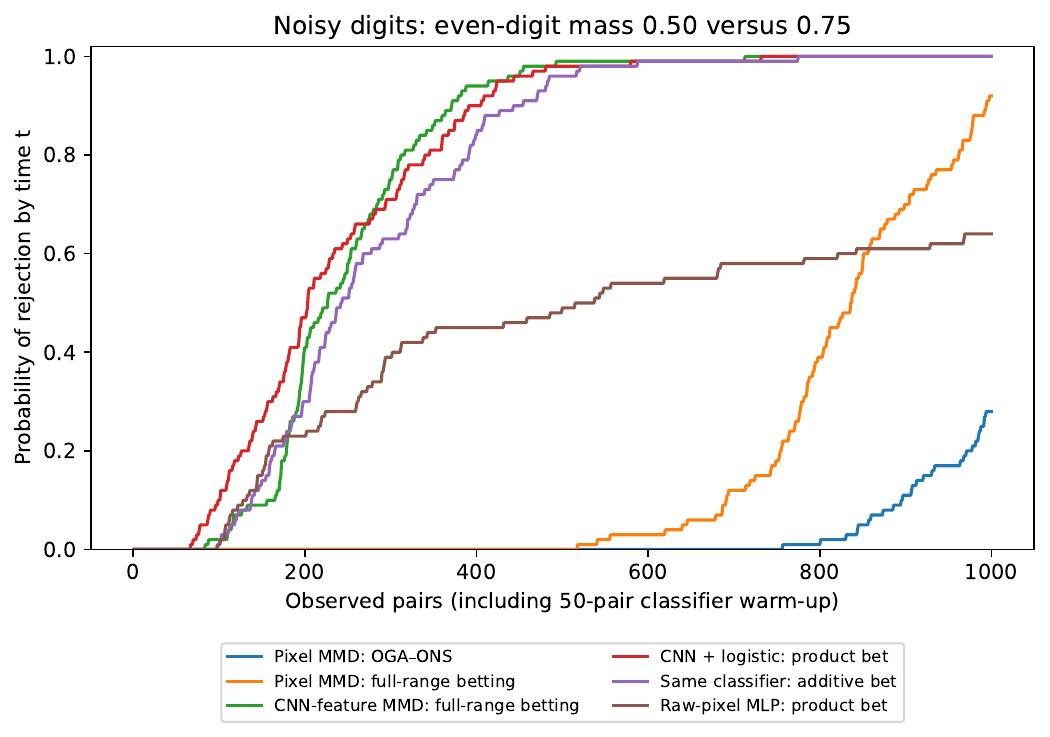}
\caption{Rejection probabilities in the weaker noisy digit-frequency shift. CNN-feature MMD and CNN product betting perform similarly; both reject substantially earlier than pixel MMD in this setting.}\label{fig:image-frequency-power}
\end{figure}

At 1000 pairs in the contrast experiment, mean normalized log wealth is $0.657$ for CNN product betting, $0.410$ for the identical-classifier additive bet, $0.081$ for feature MMD, and $0.029$ for full-range pixel MMD. In the clean frequency shift, the product and identical-classifier additive rates are $0.0324$ and $0.0328$, respectively. Predictive approximation, calibration, and the fraction mixture all affect these finite-horizon growth comparisons.

\begin{table}[tbp]\centering\small
\begin{tabular}{@{}lrr@{}}\toprule
Method & Clean null & Noisy null\\\midrule
Pixel MMD OGA--ONS & $0/200\;[0.00,1.83]$ & $0/200\;[0.00,1.83]$\\
Pixel MMD, full fractions & $4/200\;[0.55,5.04]$ & $0/200\;[0.00,1.83]$\\
Feature MMD, full fractions & $3/200\;[0.31,4.32]$ & $4/200\;[0.55,5.04]$\\
CNN source head, product & $3/200\;[0.31,4.32]$ & $2/200\;[0.12,3.57]$\\
Same classifier, additive & $3/200\;[0.31,4.32]$ & $0/200\;[0.00,1.83]$\\
Raw-pixel network, product & $2/200\;[0.12,3.57]$ & $4/200\;[0.55,5.04]$\\
\bottomrule\end{tabular}
\caption{Null rejections by 1000 pairs in the first image study at level $0.05$, followed by exact $95\%$ binomial intervals in percent. These are numerical checks of the analytically valid procedures.}\label{tab:image-null}
\end{table}

\FloatBarrier
\subsection{Fine-tuning protocol and paired comparisons}\label{app:image-finetuning}
The fine-tuning study reuses the identical split, pretrained encoder, and feature standardization. Four alternatives are specified in \cref{sec:image-results}. Every pure-class stream samples uniformly from its named digit pool. Mixture streams first sample digit 3 or 8 with the specified probabilities, then sample an image uniformly within that pool. The shared $90$-degree rotation precedes independent pixel noise. Semantic digit labels determine the simulation laws and are not supplied during source adaptation.

Both arms use 50 neutral warm-up pairs and update after every further 25 pairs. At each checkpoint, the frozen arm refits an $L^2$-regularized logistic head with $C=0.1$ on all preceding source-labeled features. The fine-tuned arm first fits the same head on its current representation, then takes ten AdamW steps using balanced replay minibatches of 32 earlier images per source. Both convolutional layers, the feature layer, and a temporary source head are updated. The encoder learning rate is $0.0003$, the temporary-head learning rate is $0.003$, weight decay is $0.001$, and the gradient norm is clipped at 5. Network weights are warm-started between checkpoints; optimizer state is reset. After these steps, the same $C=0.1$ logistic-head fit is repeated on the updated historical features. Future pairs are scored in evaluation mode. Previously recorded predictions and wealth increments remain unchanged.

Both arms use the same temperature/fraction mixture specified in \cref{sec:image-results}; an equal wealth mixture of the arms is also evaluated. The 600 trials comprise 100 replications for each alternative and 200 for the rotated-mixture null, with horizon 500 pairs. All trials continue to the horizon after any rejection to measure later growth. The main seed is 2026111001 and the independent two-replication pilot uses seed 183754. Learning configurations are unchanged between the pilot and main run. The trial comparisons are paired on the same sampled images.

\begin{table}[tbp]\centering\small\setlength{\tabcolsep}{4pt}
\begin{tabular}{@{}lrrrr@{}}\toprule
&\multicolumn{2}{c}{Normalized log-wealth difference}&\multicolumn{2}{c}{Rejection-time difference}\\
\cmidrule(lr){2-3}\cmidrule(l){4-5}
Setting & Mean & $95\%$ interval & Mean & $95\%$ interval\\\midrule
Noisy 3 versus 8 & $0.0290$ & $[0.0274,0.0306]$ & $-0.09$ & $[-0.15,-0.03]$\\
Noisy 4 versus 9 & $0.0082$ & $[0.0076,0.0087]$ & $-0.04$ & $[-0.08,-0.00]$\\
Rotated 3 versus 8 & $0.2399$ & $[0.2337,0.2462]$ & $-0.68$ & $[-0.92,-0.44]$\\
Rotated mixtures & $0.0066$ & $[0.0047,0.0084]$ & $4.09$ & $[0.78,7.40]$\\
\bottomrule\end{tabular}
\caption{Fine-tuned minus frozen in mean $\log W_{500}/500$ and mean $\tau_{0.05}\wedge500$. There are 100 paired replications in each row. Intervals are descriptive paired Student-$t$ Monte Carlo intervals, conditional on the pretrained checkpoint and image pool, and are not multiplicity-adjusted. The final endpoint in the second time interval is negative before rounding.}\label{tab:finetune-paired}
\end{table}

Every method rejects in all alternative trials by 500 pairs at both reported levels. At $\alpha=0.01$, mean rejection times for the frozen, tuned, and mixture methods are respectively $(56.72,56.52,56.64)$ for noisy 3 versus 8, $(55.37,55.29,55.31)$ for noisy 4 versus 9, $(58.91,58.24,58.40)$ for rotated 3 versus 8, and $(117.33,120.81,117.10)$ for the rotated mixture. As at level $0.05$, the weaker mixture shows an early adaptation cost despite stronger terminal growth.

\begin{figure}[tbp]\centering
\includegraphics[width=.86\linewidth]{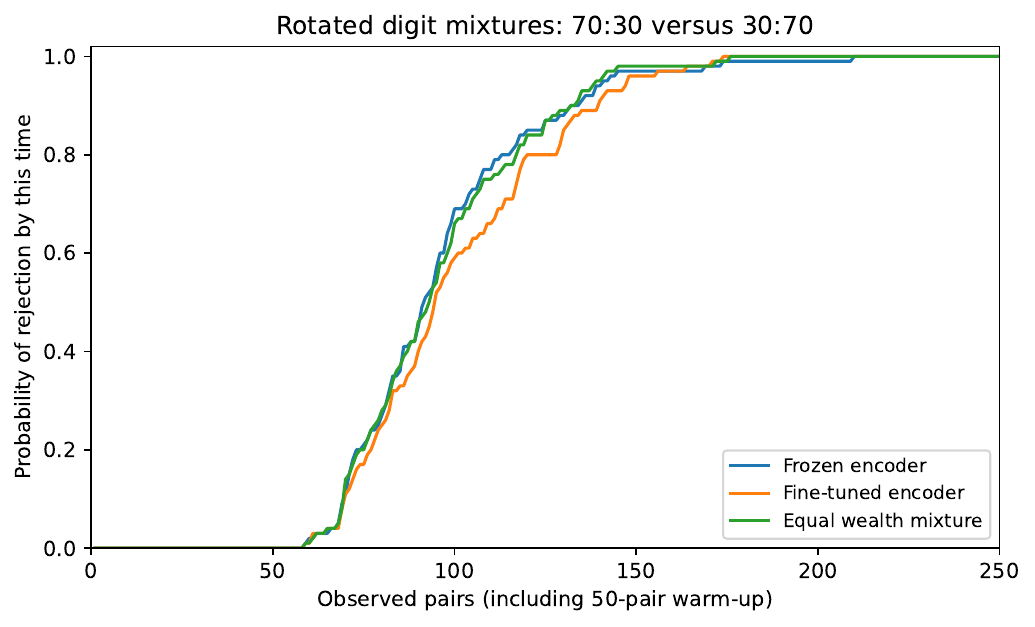}
\caption{Rejection probabilities at level $0.05$ for the rotated-digit mixtures. The source distributions have digit-3 probabilities $0.7$ and $0.3$. The frozen/tuned mixture retains the wealth guarantee \eqref{eq:frozen-tuned-mixture}; its rejection-time curve also reflects the transient learning behavior of the two arms.}\label{fig:finetune-mixture-power}
\end{figure}

\begin{table}[tbp]\centering\small
\begin{tabular}{@{}lrr@{}}\toprule
Method & Rejections at $\alpha=0.05$ & Rejections at $\alpha=0.01$\\\midrule
Frozen encoder & $1/200\;[0.01,2.75]$ & $0/200\;[0.00,1.83]$\\
Fine-tuned encoder & $3/200\;[0.31,4.32]$ & $0/200\;[0.00,1.83]$\\
Equal wealth mixture & $3/200\;[0.31,4.32]$ & $0/200\;[0.00,1.83]$\\
\bottomrule\end{tabular}
\caption{Null rejections by 500 pairs under identical rotated-digit mixtures, followed by exact $95\%$ binomial intervals in percent.}\label{tab:finetune-null}
\end{table}

\subsection{Prequential loss and probability adaptation}\label{app:finetune-calibration}
The exact identity \eqref{eq:neural-loss-comparison} concerns the unscaled full product. The reported procedures additionally mix temperatures and fractions. Retaining the prequential logits allows the two comparisons to be separated.

\begin{table}[tbp]\centering\small
\begin{tabular}{@{}lrrr@{}}\toprule
Setting & Raw product gain & Monte Carlo SE & Mixture gain\\\midrule
Noisy 3 versus 8 & $34.24$ & $0.27$ & $14.52$\\
Noisy 4 versus 9 & $15.00$ & $0.12$ & $4.08$\\
Rotated 3 versus 8 & $132.16$ & $1.11$ & $119.97$\\
Rotated mixtures & $-0.38$ & $0.68$ & $3.29$\\
\bottomrule\end{tabular}
\caption{Fine-tuned minus frozen in mean terminal log wealth, in nats over 500 pairs. Raw product uses temperature one and full betting. Mixture gain compares the temperature/fraction mixtures of the two arms, rather than their equal outer mixture. Raw-product gains equal the corresponding reduction in cumulative prequential logarithmic loss.}\label{tab:finetune-calibration}
\end{table}

For the rotated pure-class alternative, fine-tuning improves both the raw predictive loss and the adapted betting wealth. For the weaker rotated mixture, the raw-product gain is $-0.38$ nats with Monte Carlo standard error $0.68$, while temperature/fraction adaptation produces a gain of $3.29$ nats. Thus this latter terminal improvement reflects the combined representation and probability-adaptation procedure. In the null, the mean raw-product fine-tuning difference is $-9.00$ nats with standard error $0.37$, consistent with additional estimation cost when the source labels carry no population information.

\section{Proofs for fixed-sample efficiency}\label{app:batch-proofs}
Throughout this appendix $K$ is fixed. Except where stated otherwise, the alternative distributions $P,Q$ have positive probabilities on all $K$ symbols and $a/N\to\rho\in(0,1)$. The proofs allow randomized tests by adjoining independent randomization to the observation space.

\subsection{Two-sample types and the reverse projection}
For a type $S$ with denominator $a$, its multinomial probability under $P$ satisfies
\begin{equation}\label{eq:ts-types}
 (a+1)^{-K}e^{-a\KL(S\|P)}
 \le\PP_P(\widehat P_a=S)\le e^{-a\KL(S\|P)}.
\end{equation}
The corresponding bounds for $\widehat Q_b=T$ multiply, since the samples are independent. There are at most $(a+1)^K(b+1)^K$ type pairs. These inequalities follow from the standard type-class bounds
$(a+1)^{-K}e^{aH(S)}\le |\mathcal T_a(S)|\le e^{aH(S)}$.
One direct proof of the lower bound observes that the type $S$ is a mode of the multinomial distribution with parameter $S$, and there are at most $(a+1)^K$ types. Every string of that type has probability $e^{-aH(S)}$ under $S$. Zero cells can be omitted. Changing the sampling law from $S$ to $P$ gives \eqref{eq:ts-types}.

Compactness and continuity on the finite simplices now give the large-deviation rate
$\mathcal D_\rho(S,T;P,Q)$ for empirical pairs. More explicitly, the upper bound for any set follows by summing \eqref{eq:ts-types}, and a lower bound follows by choosing types converging to any interior point of that set. The functions $G_{\rho_N}$ and $\mathcal D_{\rho_N}$ converge uniformly to their $\rho$ versions, since $P,Q$ are positive and $\rho_N$ stays in a compact subinterval of $(0,1)$.

The reverse projection in \eqref{eq:ts-J} has the explicit minimizer
\[
 r_i^*=\frac{p_i^\rho q_i^{1-\rho}}{Z_\rho},\qquad
 Z_\rho=\sum_i p_i^\rho q_i^{1-\rho}.
\]
Indeed,
\[
 \rho\KL(R\|P)+(1-\rho)\KL(R\|Q)
 =\KL(R\|R^*)-\log Z_\rho.
\]
For $P\ne Q$, strict H\"older inequality gives $Z_\rho<1$, so $J_\rho>0$.

\subsection{Proof of the Stein exponent and its converse}
By \eqref{eq:ts-code-sandwich}, acceptance at fixed level implies
\[
 G_{\rho_N}(\widehat P_a,\widehat Q_b)
 \le \frac{\log(1/\alpha)+d_{a,b}}N\longrightarrow0.
\]
The type bound and compactness imply that its probability under $P^a\otimes Q^b$ has limit inferior exponent at least
\[
 \inf_{S,T:G_\rho(S,T)=0}\mathcal D_\rho(S,T;P,Q)
 =\inf_R\mathcal D_\rho(R,R;P,Q)=J_\rho.
\]
Here $G_\rho(S,T)=0$ if and only if $S=T$.

For the converse, fix a common law $R$ and a test with level at most $\alpha$. Let $A_N$ be its acceptance event, including any test randomization, and let
$L_N=dR^N/d(P^a\otimes Q^b)$. Under $R$, streamwise strong laws give
\[
 N^{-1}\log L_N\longrightarrow
 \rho\KL(R\|P)+(1-\rho)\KL(R\|Q)=:D_R.
\]
For any $\varepsilon>0$,
\[
 \beta_{a,b}(P,Q)\ge
 e^{-N(D_R+\varepsilon)}
 R^N\{A_N,\log L_N\le N(D_R+\varepsilon)\}.
\]
The final probability has limit inferior at least $1-\alpha$. Thus the limit superior miss exponent is at most $D_R$; minimizing over $R$ proves the converse and equality in \eqref{eq:ts-stein}.

\subsection{Proof of the Hoeffding tradeoff}
For threshold $e^{Nc}$, the coding sandwich places the acceptance region between empirical-divergence sublevel sets at $c$ and $c+d_{a,b}/N$. The type bounds give a limit inferior exponent of at least $B_\rho(c;P,Q)$. A type sequence approaching any $(S,T)$ with $G_\rho(S,T)<c$ is eventually accepted, proving the reverse bound with the open constraint.

The open and closed infima agree for $c>0$. If $G_\rho(S,T)=c$, let $M=\rho S+(1-\rho)T$ and set
\[
 S_\delta=(1-\delta)S+\delta M,\qquad
 T_\delta=(1-\delta)T+\delta M.
\]
The pooled distribution remains $M$, and convexity of relative entropy gives
$G_\rho(S_\delta,T_\delta)\le(1-\delta)G_\rho(S,T)<c$.
Meanwhile $\mathcal D_\rho(S_\delta,T_\delta;P,Q)\to\mathcal D_\rho(S,T;P,Q)$. This proves the asserted exponent for the code.

For the converse, fix $(S,T)$ with $G_\rho(S,T)<c$ and the common null $R=\rho S+(1-\rho)T$. Under $S^a\otimes T^b$, the normalized log likelihood ratio against $R^N$ tends to $G_\rho(S,T)$. A level bound $e^{-Nc}$ and a likelihood-ratio split therefore imply
\[
 \PP_{S^a\otimes T^b}(\text{reject})
 \le o(1)+e^{-N\{c-G_\rho(S,T)-\varepsilon\}}
 \longrightarrow0
\]
for sufficiently small positive $\varepsilon$. A second change of measure, from $S^a\otimes T^b$ to $P^a\otimes Q^b$, then bounds the limit superior miss exponent by $\mathcal D_\rho(S,T;P,Q)$. Minimize over the open constraint and use the preceding continuity argument. This proves \cref{thm:ts-errors}.

\subsection{Proof of the Chernoff characterization}
For a fixed $R$, the minimum sum of errors for the simple pair
$(R^a\otimes R^b,P^a\otimes Q^b)$ is the sum of the pointwise minimum of the two likelihoods. Grouping by the two sample types and applying \eqref{eq:ts-types} shows that its exponent is
\begin{equation}\label{eq:ts-C-variational}
 \inf_{S,T}\max\{\mathcal D_\rho(S,T;R,R),
                         \mathcal D_\rho(S,T;P,Q)\}.
\end{equation}
Writing the maximum as a supremum over a single weight $u\in[0,1]$ and applying finite-dimensional convex minimax proves equality with \eqref{eq:ts-C-R}. For fixed $u$, minimization over $S$ and $T$ separates, giving the two logarithms of exponential affinities. If $R$ has zero cells, restrict to their common support; one-sided endpoint limits give the same identity.

For fixed $(S,T)$, the first divergence in \eqref{eq:ts-C-variational} is minimized by $R=\rho S+(1-\rho)T$, with minimum $G_\rho(S,T)$. Interchanging the infima gives \eqref{eq:ts-C}.

The function $B_\rho$ is continuous and nonincreasing on $[0,G_\rho(P,Q)]$. Right continuity, including at zero, follows by compactness of minimizing pairs. Left continuity at positive arguments follows by the pooled-mixture perturbation from the preceding proof. It has $B_\rho(0)=J_\rho>0$ and $B_\rho(G_\rho(P,Q))=0$. Hence the strictly decreasing continuous function $B_\rho(c)-c$ has a unique zero $c_*$. For every pair with $G_\rho(S,T)\le c_*$ its alternative divergence is at least $c_*$, while every other pair has $G_\rho(S,T)>c_*$. A minimizer of $B_\rho(c_*)$ makes the maximum equal to $c_*$. This proves $C_\rho^{\rm hom}=c_*$.

The code test at this threshold has type-I error at most $e^{-Nc_*}$ and miss exponent $c_*$, by \cref{thm:ts-errors}. Its summed-error exponent is at least $c_*$. For any common $R$, the worst-case sum is at least the sum for the simple pair with null $R$. The preceding simple-hypothesis result bounds every test's exponent by $C_\rho(R;P,Q)$; minimize over $R$. This proves \cref{thm:ts-chernoff}.

\subsection{Permutation ordering, significance, and error exponents}
Condition on the pooled counts $C$ and on sample sizes $a,b$. Under the common-law null, the allocation counts have probability
\[
 \mathfrak h(A\mid C)=\frac{\prod_i\binom{C_i}{A_i}}{\binom Na}
 =\frac{a!b!\prod_i C_i!}{N!\prod_i A_i!B_i!}.
\]
The code's gamma normalizers and pooled null likelihood are fixed at these margins, so $V_{a,b}$ is proportional to $\prod_i A_i!B_i!$. Thus its upper-tail permutation probability sums exactly those allocations having probability no greater than the observed allocation. This proves the Fisher--Freeman--Halton identity and finite-sample conditional validity, including its conservative convention at ties.

There are at most $(N+1)^K$ possible allocation tables. Consequently
\[
 \mathfrak h(A^{\rm obs}\mid C)\le p^{\rm perm}_{a,b}
 \le (N+1)^K\mathfrak h(A^{\rm obs}\mid C).
\]
The bound $\log(n!)=n\log n-n+O(\log(n+1))$, uniformly for $0\le n\le N$, gives
\[
 -\log\mathfrak h(A\mid C)=T_{a,b}+O_K(\log(N+1)).
\]
This proves \eqref{eq:ts-perm-sandwich}, including tables with empty cells. The strong law gives the Bahadur slope. For a fixed significance threshold, acceptance implies $T_{a,b}=O_K(\log N)$; the Stein proof applies. For exponentially small significance levels, the permutation sandwich puts acceptance between the same divergence sublevel sets up to $O_K(\log N/N)$, so the Hoeffding type argument applies. Its exact validity supplies the required null bound. The Chernoff conclusion follows at $c_*$. This proves all the global-efficiency assertions in \cref{thm:ts-pitman}.

\subsection{Local likelihood-ratio and coding expansions}
Under \eqref{eq:ts-local-model}, independent triangular-array multinomial central limit theorems give
\[
 Z_N:=\sqrt{\frac{ab}{N}}(\widehat P_a-\widehat Q_b)
 \Rightarrow N\bigl(\sqrt{\rho(1-\rho)}(u-v),\Sigma_r\bigr),
 \quad \Sigma_r=\operatorname{diag}(r)-rr^\top.
\]
These central limit theorems follow, for example, by the second-order expansion of the characteristic functions of the bounded categorical summands. Taylor expansion of the two relative entropies around their common pooled distribution yields
\[
 2T_{a,b}=\sum_i Z_{N,i}^2/r_i+o_P(1).
\]
After standardization by $\sqrt{r_i}$, the covariance is the orthogonal projection onto the $(K-1)$-dimensional subspace perpendicular to $(\sqrt{r_i})_i$. The standardized mean lies in that subspace and has squared norm $\Lambda$. This proves \eqref{eq:ts-local-lrt}.

For counts $U$ of a sample of size $n$ with all empirical probabilities converging to positive limits, Stirling expansion gives
\[
 \log q_K(U)-\sum_i U_i\log(U_i/n)
 =-\frac{K-1}{2}\log n+\log\Gamma(K)
 +\frac{K-1}{2}\log(2\pi)+\frac12\sum_i\log(U_i/n)+o_P(1).
\]
Apply this to the two samples and use $\log(ab)=2\log N+\log\{\rho(1-\rho)\}+o(1)$. Subtracting the pooled null log likelihood gives \eqref{eq:ts-local-code}. Since $T_{a,b}=O_P(1)$, the direct code tends to zero in probability and its fixed-threshold power tends to zero.

To establish the conditional reference limit, consider deterministic pooled proportions $C/N\to r>0$. Under the multivariate hypergeometric law, $A$ has mean $aC/N$ and covariance
\[
 \operatorname{Cov}(A\mid C)=\frac{ab}{N-1}
 \{\operatorname{diag}(C/N)-(C/N)(C/N)^\top\}.
\]
Thus $\sqrt{N/(ab)}(A-aC/N)$ is tight. Uniform Stirling expansion on each bounded region in its $(K-1)$-dimensional lattice gives masses proportional to a centered Gaussian density with covariance $\Sigma_r$ and lattice spacing of order $N^{-1/2}$. Riemann sums and tightness identify the conditional Gaussian limit. The same quadratic expansion gives conditional convergence of $2T_{a,b}$ to $\chi^2_{K-1}$.

The additive correction in \eqref{eq:ts-local-code}, after replacing $r$ by $C/N$ and $\rho$ by $a/N$, depends only on the conditioned margins up to $o_P(1)$. Both the observed table under the local alternative and a random permutation concentrate in this expansion regime. The continuous chi-square limiting distribution therefore makes the permutation cutoff asymptotically the likelihood-ratio cutoff. The random pooled proportions under \eqref{eq:ts-local-model} converge to $r$, so the deterministic-margin statement applies in probability. Combining it with the noncentral limit proves \eqref{eq:ts-pitman-power}.

\subsection{Proof of the boundary exponent theorem}
Write $u=\min\{a,bL(b)\}$. Then $u/b\to\infty$ and
$\log u=o(b)$ regardless of the growth rate of the original $a$.
For each fixed $\delta>0$, Hoeffding's inequality gives
\[
 \PP_P(\|\widehat P_u-P\|_\infty>\delta)
 \le2K e^{-2u\delta^2}.
\]
Its decay is superexponential on the $b$ scale. On the complementary event, with $\delta<\min_i p_i/2$, Taylor expansion in $b/(u+b)$ is uniform over all possible $\widehat Q_b$ and gives
\[
 b^{-1}T_{u,b}=\KL(\widehat Q_b\|\widehat P_u)+O(b/u).
\]
Also, uniformly in the second argument's empirical probabilities,
$|\KL(\widehat Q_b\|\widehat P_u)-\KL(\widehat Q_b\|P)|\le C_P\delta$.
The code and endpoint-correction penalties are $O_K(\log u+\log b)=o(b)$. It follows that
\begin{equation}\label{eq:ts-boundary-equiv}
 b^{-1}\log U_{a,b}\quad\hbox{and}\quad\KL(\widehat Q_b\|P)
\end{equation}
are exponentially equivalent at speed $b$: for every positive tolerance, the probability that they differ by more than that tolerance has logarithm divided by $b$ tending to $-\infty$.

Applying the one-sample type bounds under $Q$ gives miss-exponent lower bounds $\KL(P\|Q)$ at a fixed level and $B_\infty(c;P,Q)$ at threshold $e^{bc}$. Open and closed infima at positive $c$ agree by mixing a boundary candidate $S$ with $P$. Null error is bounded by the inverse threshold because $U$ is an e-process.

For the converse choose the common null $R=P$. The $X$ observations have exactly the same law under this null and the alternative, so the likelihood ratios involve only the $b$ observations from the second stream. The Stein likelihood-ratio argument gives upper bound $\KL(P\|Q)$, even for tests using every original $X$ observation. For the exponential-level converse, choose any $S$ with $\KL(S\|P)<c$ and compare $P^a\otimes S^b$ first to $P^{a+b}$ and then to $P^a\otimes Q^b$. The two-change-of-measure argument gives upper bound $\KL(S\|Q)$; minimize and take the closure. These converses match the lower bounds from \eqref{eq:ts-boundary-equiv}.

The function $B_\infty$ is continuous and nonincreasing from $\KL(P\|Q)$ to zero on $[0,\KL(Q\|P)]$. Its crossing with $c$ equals
$\inf_S\max\{\KL(S\|P),\KL(S\|Q)\}=C(P,Q)$ by the simple-hypothesis Chernoff argument. The common null $R=P$ gives the corresponding converse for every test, while the capped code at this crossing attains it. A mixture with positive weight $\varepsilon$ on $U$ has miss probability no larger than that of $U$ at a threshold increased by $1/\varepsilon$. This constant has no effect on any exponent; the mixture's validity and the same converses complete \cref{thm:ts-boundary-errors}.

\paragraph{Acknowledgments.}
ChatGPT assisted with mathematical derivations, manuscript preparation, reference implementations, and execution and analysis of the experiments.

\begingroup\small
\setlength{\bibsep}{1.5pt}

\endgroup

\begin{thebibliography}{99}
\providecommand{\natexlab}[1]{#1}
\providecommand{\url}[1]{\texttt{#1}}
\expandafter\ifx\csname urlstyle\endcsname\relax
  \providecommand{\doi}[1]{doi: #1}\else
  \providecommand{\doi}{doi: \begingroup \urlstyle{rm}\Url}\fi


\bibitem[Ramdas(2026)]{ramdas2026independence}
Aaditya Ramdas.
\newblock Log-optimal independence testing by betting.
\newblock Stanford Digital Repository (Technical Report)  \href{https://purl.stanford.edu/pg979gk9581}{https://purl.stanford.edu/pg979gk9581}, 2026.

\bibitem[Alpaydin and Kaynak(1998)]{uci1998}
E.~Alpaydin and C.~Kaynak.
\newblock Optical recognition of handwritten digits.
\newblock UCI Machine Learning Repository, 1998.
\newblock \doi{10.24432/C50P49}.

\bibitem[Balsubramani and Ramdas(2016)]{balsubramani2016}
Akshay Balsubramani and Aaditya Ramdas.
\newblock Sequential nonparametric testing with the law of the iterated
  logarithm.
\newblock Proceedings of the Conference on Uncertainty in Artificial
  Intelligence, 2016.
\newblock URL \url{https://arxiv.org/abs/1506.03486}.

\bibitem[Cover(1991)]{cover1991}
Thomas~M. Cover.
\newblock Universal portfolios.
\newblock \emph{Mathematical Finance}, 1\penalty0 (1):\penalty0 1--29, 1991.
\newblock URL \url{https://www-isl.stanford.edu/~cover/papers/paper93.pdf}.

\bibitem[Gretton et~al.(2012)Gretton, Borgwardt, Rasch, Sch{\"o}lkopf, and
  Smola]{gretton2012kernel}
Arthur Gretton, Karsten~M. Borgwardt, Malte~J. Rasch, Bernhard Sch{\"o}lkopf,
  and Alexander Smola.
\newblock A kernel two-sample test.
\newblock \emph{Journal of Machine Learning Research}, 13\penalty0
  (25):\penalty0 723--773, 2012.
\newblock URL \url{https://jmlr.org/papers/v13/gretton12a.html}.

\bibitem[Gr{\"u}nwald et~al.(2024)Gr{\"u}nwald, de~Heide, and
  Koolen]{grunwald2024safe}
Peter Gr{\"u}nwald, Rianne de~Heide, and Wouter Koolen.
\newblock Safe testing.
\newblock \emph{Journal of the Royal Statistical Society Series B: Statistical
  Methodology}, 2024.
\newblock \doi{10.1093/jrsssb/qkae011}.
\newblock URL \url{https://doi.org/10.1093/jrsssb/qkae011}.

\bibitem[Kelly(1956)]{kelly1956}
John~L. Kelly, Jr.
\newblock A new interpretation of information rate.
\newblock \emph{Bell System Technical Journal}, 35\penalty0 (4):\penalty0
  917--926, 1956.
\newblock \doi{10.1002/j.1538-7305.1956.tb03809.x}.

\bibitem[Lh{\'e}ritier and Cazals(2018)]{lheritier2018}
Alix Lh{\'e}ritier and Fr{\'e}d{\'e}ric Cazals.
\newblock A sequential non-parametric multivariate two-sample test.
\newblock \emph{IEEE Transactions on Information Theory}, 64\penalty0
  (5):\penalty0 3361--3370, 2018.
\newblock \doi{10.1109/TIT.2018.2800658}.
\newblock URL \url{https://ieeexplore.ieee.org/document/8276312/}.

\bibitem[Lh{\'e}ritier and Cazals(2019)]{lheritier2019}
Alix Lh{\'e}ritier and Fr{\'e}d{\'e}ric Cazals.
\newblock Low-complexity nonparametric Bayesian online prediction with
  universal guarantees.
\newblock In \emph{Advances in Neural Information Processing Systems},
  volume~32, 2019.
\newblock URL
  \url{https://proceedings.neurips.cc/paper_files/paper/2019/hash/122e27d57ae8ecb37f3f1da67abb33cb-Abstract.html}.

\bibitem[Orabona(2026)]{orabona2019}
Francesco Orabona.
\newblock Online learning: A modern introduction using convex optimization.
\newblock arXiv:1912.13213v10, 2026.
\newblock URL \url{https://arxiv.org/abs/1912.13213v10}.

\bibitem[Pandeva et~al.(2024)Pandeva, Bakker, Naesseth, and
  Forr{\'e}]{pandeva2024}
Teodora Pandeva, Tim Bakker, Christian~A. Naesseth, and Patrick Forr{\'e}.
\newblock E-valuating classifier two-sample tests.
\newblock arXiv:2210.13027v2, 2024.
\newblock URL \url{https://arxiv.org/abs/2210.13027}.

\bibitem[Pandeva et~al.(2024)Pandeva, Forr{\'e}, Ramdas, and Shekhar]{pandeva2024deep}
Teodora Pandeva, Patrick Forr{\'e}, Aaditya Ramdas, and Shubhanshu Shekhar.
\newblock Deep anytime-valid hypothesis testing.
\newblock In \emph{Proceedings of the 27th International Conference on Artificial Intelligence and Statistics},
volume 238 of \emph{Proceedings of Machine Learning Research}, pages 622--630, 2024.
\newblock URL \url{https://proceedings.mlr.press/v238/pandeva24a.html}.

\bibitem[Podkopaev and Ramdas(2023)]{podkopaev2023predictive}
Aleksandr Podkopaev and Aaditya Ramdas.
\newblock Sequential predictive two-sample and independence testing.
\newblock arXiv:2305.00143v2, 2023.
\newblock URL \url{https://arxiv.org/abs/2305.00143}.

\bibitem[Ramdas et~al.(2023)Ramdas, Gr{\"u}nwald, Vovk, and
  Shafer]{ramdas2023savi}
Aaditya Ramdas, Peter Gr{\"u}nwald, Vladimir Vovk, and Glenn Shafer.
\newblock Game-theoretic statistics and safe anytime-valid inference.
\newblock \emph{Statistical Science}, 38\penalty0 (4), 2023.
\newblock \doi{10.1214/23-STS894}.
\newblock URL \url{https://arxiv.org/abs/2210.01948}.

\bibitem[scikit-learn developers(2026)]{sklearn-digits}
scikit-learn developers.
\newblock \texttt{sklearn.datasets.load\_digits}: dataset documentation.
\newblock Accessed November 10, 2026.
\newblock URL \url{https://scikit-learn.org/stable/modules/generated/sklearn.datasets.load_digits.html}.

\bibitem[Shekhar and Ramdas(2024)]{shekhar2024betting}
Shubhanshu Shekhar and Aaditya Ramdas.
\newblock Nonparametric two-sample testing by betting.
\newblock \emph{IEEE Transactions on Information Theory}, 70\penalty0
  (2):\penalty0 1178--1203, 2024.
\newblock \doi{10.1109/TIT.2023.3305867}.
\newblock URL \url{https://arxiv.org/abs/2112.09162}.

\bibitem[Shekhar(2025)]{shekhar2025optimal}
Shubhanshu Shekhar.
\newblock Optimal anytime-valid tests for composite nulls.
\newblock arXiv:2512.20039, 2025.
\newblock URL \url{https://arxiv.org/abs/2512.20039}.

\bibitem[Sriperumbudur et~al.(2011)Sriperumbudur, Fukumizu, and
  Lanckriet]{sriperumbudur2011}
Bharath~K. Sriperumbudur, Kenji Fukumizu, and Gert R.~G. Lanckriet.
\newblock Universality, characteristic kernels and {RKHS} embedding of
  measures.
\newblock \emph{Journal of Machine Learning Research}, 12\penalty0
  (70):\penalty0 2389--2410, 2011.
\newblock URL \url{https://www.jmlr.org/papers/v12/sriperumbudur11a.html}.

\bibitem[Turner et~al.(2024)Turner, Ly, and Gr{\"u}nwald]{turner2024generic}
Rosanne~J. Turner, Alexander Ly, and Peter Gr{\"u}nwald.
\newblock Generic e-variables for exact sequential $k$-sample tests that allow
  for optional stopping.
\newblock \emph{Journal of Statistical Planning and Inference}, 2024.
\newblock \doi{10.1016/j.jspi.2023.106116}.
\newblock URL \url{https://arxiv.org/abs/2106.02693}.

\bibitem[Wasserman et~al.(2020)Wasserman, Ramdas, and Balakrishnan]{wasserman2020universal}
Larry Wasserman, Aaditya Ramdas, and Sivaraman Balakrishnan.
\newblock Universal inference.
\newblock \emph{Proceedings of the National Academy of Sciences}, 117(29):16880--16890, 2020.
\newblock \doi{10.1073/pnas.1922664117}.
\newblock URL \url{https://arxiv.org/abs/1912.11436}.

\bibitem[Waudby-Smith et~al.(2026)Waudby-Smith, Sandoval, and
  Jordan]{waudbysmith2026}
Ian Waudby-Smith, Ricardo~J. Sandoval, and Michael~I. Jordan.
\newblock Instance-log-optimality of portfolio-based e-processes and their
  sequential hypothesis tests.
\newblock arXiv:2504.02818v2, 2026.
\newblock URL \url{https://arxiv.org/abs/2504.02818}.

\bibitem[Gutman(1989)]{gutman1989}
Michael Gutman.
\newblock Asymptotically optimal classification for multiple tests with empirically observed statistics.
\newblock \emph{IEEE Transactions on Information Theory}, 35(2):401--408, 1989.
\newblock \doi{10.1109/18.32134}.

\bibitem[Unnikrishnan(2012)]{unnikrishnan2012}
Jayakrishnan Unnikrishnan.
\newblock On optimal two sample homogeneity tests for finite alphabets.
\newblock In \emph{IEEE International Symposium on Information Theory}, pages 2027--2031, 2012.

\bibitem[Zhu et~al.(2021)Zhu, Chen, Chen, and Yang]{zhu2021}
Shengyu Zhu, Biao Chen, Zhitang Chen, and Pengfei Yang.
\newblock Asymptotically optimal one- and two-sample testing with kernels.
\newblock \emph{IEEE Transactions on Information Theory}, 67(4):2074--2092, 2021.
\newblock \doi{10.1109/TIT.2021.3059267}.
\newblock URL \url{https://arxiv.org/abs/1908.10037}.

\bibitem[Harsha et~al.(2026)Harsha, Ravi, and Koch]{harsha2026}
K.~V. Harsha, Jithin Ravi, and Tobias Koch.
\newblock Second-order asymptotics of two-sample tests.
\newblock arXiv:2601.09196, 2026.
\newblock URL \url{https://arxiv.org/abs/2601.09196}.

\bibitem[Grootveld et~al.(2026)Grootveld, Chen, and Gandikota]{grootveld2026}
Arick Grootveld, Biao Chen, and Venkata Gandikota.
\newblock Asymptotically optimal tests for one- and two-sample problems.
\newblock arXiv:2601.11727, 2026.
\newblock URL \url{https://arxiv.org/abs/2601.11727}.

\bibitem[Lapidoth and Pfister(2018)]{lapidoth2018}
Amos Lapidoth and Christoph Pfister.
\newblock Testing against independence and a R\'enyi information measure.
\newblock In \emph{IEEE Information Theory Workshop}, 2018.
\newblock URL \url{https://arxiv.org/abs/1805.11059}.


\bibitem[Bahadur(1960)]{bahadur1960}
R. R. Bahadur.
\newblock Stochastic comparison of tests.
\newblock \emph{Annals of Mathematical Statistics}, 31(2):276--295, 1960.

\bibitem[Chernoff(1952)]{chernoff1952}
Herman Chernoff.
\newblock A measure of asymptotic efficiency for tests of a hypothesis based on the sum of observations.
\newblock \emph{Annals of Mathematical Statistics}, 23(4):493--507, 1952.
\newblock \url{https://doi.org/10.1214/aoms/1177729330}.

\bibitem[Hoeffding(1965)]{hoeffding1965}
Wassily Hoeffding.
\newblock Asymptotically optimal tests for multinomial distributions.
\newblock \emph{Annals of Mathematical Statistics}, 36(2):369--401, 1965.

\bibitem[Freeman and Halton(1951)]{freeman1951}
G. H. Freeman and J. H. Halton.
\newblock Note on an exact treatment of contingency, goodness of fit and other problems of significance.
\newblock \emph{Biometrika}, 38(1--2):141--149, 1951.
\newblock \url{https://doi.org/10.1093/biomet/38.1-2.141}.

\bibitem[Wilks(1938)]{wilks1938}
Samuel S. Wilks.
\newblock The large-sample distribution of the likelihood ratio for testing composite hypotheses.
\newblock \emph{Annals of Mathematical Statistics}, 9(1):60--62, 1938.
\newblock \url{https://doi.org/10.1214/aoms/1177732360}.

\end{thebibliography}
\end{document}